\documentclass[a4paper,fleqn]{cas-sc-arxiv}

\usepackage[authoryear]{natbib}

\usepackage[utf8]{inputenc}
\usepackage{caption}
\usepackage{subcaption}
\usepackage{bm}
\usepackage{algpseudocode}
\usepackage{mathtools}

\begin{document}
\let\WriteBookmarks\relax
\def\floatpagepagefraction{1}
\def\textpagefraction{.001}
\shorttitle{Modified Dynamic Programming Algorithms for GLOSA Systems with Stochastic Signal Switching Times}
\shortauthors{P. Typaldos et~al.}

\title [mode = title]{Modified Dynamic Programming Algorithms for GLOSA Systems with Stochastic Signal Switching Times}

\author[1]{Panagiotis Typaldos}
\cormark[1]

\author[1,2]{Markos Papageorgiou}

\address[1]{Dynamic Systems and Simulation Laboratory, School of Production Engineering \& Management, Technical University of Crete, 73100, Chania, Greece}
\address[2]{Faculty of Maritime and Transportation, Ningbo University, Ningbo, China}

\cortext[cor1]{Corresponding author}

\begin{abstract}
A discrete-time stochastic optimal control problem was recently proposed to address the GLOSA (Green Light Optimal Speed Advisory) problem in cases where the next signal switching time is decided in real time and is therefore uncertain in advance. The corresponding numerical solution via SDP (Stochastic Dynamic Programming) calls for substantial computation time, which excludes problem solution in the vehicle's on-board computer in real time. To overcome the computation time bottleneck, as a first attempt, a modified version of Dynamic Programming, known as Discrete Differential Dynamic Programming (DDDP) was recently employed for the numerical solution of the stochastic optimal control problem. The DDDP algorithm was demonstrated to achieve results equivalent to those obtained with the ordinary SDP algorithm, albeit with significantly reduced computation times. The present work considers a different modified version of Dynamic Programming, known as Differential Dynamic Programming (DDP). For the stochastic GLOSA problem, it is demonstrated that DDP achieves quasi-instantaneous (extremely fast) solutions in terms of CPU times, which allows for the proposed approach to be readily executable online, in an MPC (Model Predictive Control) framework, in the vehicle's on-board computer. The approach is demonstrated by use of realistic examples. It should be noted that DDP does not require discretization of variables, hence the obtained solutions may be slightly superior to the standard SDP solutions.
\end{abstract}


\maketitle

\section{Introduction}
In an era of climate crisis, actions toward environment protection are due in Transportation systems, which have a significant impact on environmental pollution and energy resources shortage. Improved fuel efficiency in road vehicles implies reduced fuel consumption, emissions and expenses for the driver. Efficient intelligent transportation systems, including real-time traffic signals, contribute to this effort. Traffic lights were introduced, in the first place, as a measure to ensure safe crossing of antagonistic streams of vehicles, cyclists and pedestrians at urban junctions. However, enforcing safety via traffic lights affects the fuel consumption, as some vehicles are obliged to fully stop in front of a red light; and then accelerate, after the traffic light switching to green. To reduce the resulting vehicle delays and number of stops, several algorithms have been proposed and deployed over the past decades, aiming at optimizing the traffic signals operation \citep{hounsell2001urban, papageorgiou2003review, kosmatopoulos2006international}. In this context, fuel consumption is increasingly considered as an optimization or evaluation criterion, while developing and deploying signal control systems \citep{jamshidnejad2017sustainable}.

Fixed-time signal plans are derived offline based on historical constant demands; and are applied without deviations. This implies that switching times of the traffic lights are known in advance. In contrast, real-time signal control strategies make use of real-time measurements to compute in real time suitable signal settings. Depending on the employed signal control strategy, the control update period may range from one second to one signal cycle. Clearly, for real-time signals, the next switching time is not known before the switching decision has been actually made.

A common dilemma, for a vehicle approaching a red traffic light, is whether it should slow down or maintain its speed at risk to stop if the traffic light has not switched to green at arrival. This dilemma may be addressed by appropriately designed on-board systems, which receive the current state and timing of the traffic signal and guide the driver (or an automated vehicle), while approaching the traffic light, by giving speed advice that ensures that the vehicle will cross the traffic signal at green and with minimum fuel consumption and emissions. Systems (or apps) optimizing the vehicle approach to traffic lights are often referred to as Green Light Optimal Speed Advisory (GLOSA) systems \citep{stahlmann2016technical}.

In the case of fixed signals, and hence prior knowledge of the next switching time, a corresponding message is broadcasted by the signal controller well before switching. Under these conditions, the problem of how to optimize the approach to traffic signals has been addressed in different ways, including rule-based algorithms \citep{katsaros2011performance, sanchez2006predicting, ma2018hardware} and optimal control approaches, which consider explicitly the vehicle kinematics and are, by their nature, more efficient in producing fuel-optimal speed profiles \citep{lawitzky2013energy, typaldos2020vehicle}. Some of these works take into consideration additional factors, such as the car waiting queue. In \citep{yang2016eco}, an eco-CACC (Cooperative Adaptive Cruise Control) scheme was developed, where optimal speed trajectories were calculated by minimizing fuel-consumption and ensuring that the vehicles arrive at the intersection as soon as the last queued vehicle is discharged. An extension of this work for the case of multiple signalized intersections is presented in \citep{yang2020eco}. In the literature, most existing GLOSA problems are focusing on developing strategies for a single vehicle, however, there are few works that approach the problem in case of platoons of vehicles in mixed traffic environment. In \citep{zhao2018platoon}, a model predictive control (MPC) approach was proposed that minimizes the fuel-consumption of a platoon of mixed vehicles (automated and human-driven vehicles) while crossing an intersection during the green phase. 

On the other hand, the case of real-time signals with very short (e.g., second-by-second) control update periods is more complicated, as there is no exact prior knowledge of the next switching time. In this case, the best available knowledge can be presented as an estimate or as a probabilistic distribution for the next switching time within a short-term future time-window. A comparison of several signal control methods with respect to traffic efficiency, plan stability and the resulting GLOSA performance is presented in \citep{blokpoel2017signal}. In \citep{koukoumidis2011signalguru}, a system called SignalGuru was developed, which detects and predicts the traffic signal timing using a collection of mobile phone data. Instead of a single estimate, \cite{mahler2012reducing} produced a probability distribution for the next switching time based on past signal statistics. However, the available probability distribution is used heuristically (rather than optimally) to accordingly time-weigh the objective function within a deterministic optimal control problem that is solved via a dynamic programming (DP) algorithm. In \citep{sun2020optimal}, the uncertainty of traffic signal timing is considered as a data-driven chance-constrained robust optimization problem, where the red-light duration is defined as a random variable. A data-driven approach is adopted to formulate chance constraints, based on empirical sample data, and Dynamic Programming is employed to solve the optimization problem. In \citep{lawitzky2013energy}, a stochastic optimal control problem was proposed, which uses a probabilistic distribution of the next signal switching time, and the problem is solved using a discrete stochastic dynamic programming (SDP) algorithm.

In \citep{typaldos2020vehicle}, the problem of producing fuel-optimal vehicle trajectories for a vehicle approaching a traffic signal was considered for both cases of known and stochastic switching times. For the first case, an optimal control problem was formulated and solved analytically via PMP (Pontryagin's Maximum Principle). Subsequently, the case of stochastic switching time with known probability distribution was also addressed in the format of a stochastic optimal control problem, which was solved numerically using SDP. The proposed SDP algorithm may take several minutes to execute, which implies that the solution is not real-time feasible and can therefore not be obtained on-board the vehicle, but must be executed offline, on the infrastructure side, and be communicated to approaching vehicles according to their current states. To substantially reduce the computation time and memory requirements for the solution of the stochastic GLOSA problem and enable its solution on-board the vehicle, \cite{typaldos2021vehicle} employed a Discrete Differential Dynamic Programming (DDDP) algorithm \citep{heidari1971discrete}. DDDP is an iterative algorithm, each iteration receiving a feasible (but non-optimal) trajectory and transforming it to an enhanced one, to be used in the next iteration, and so forth, until convergence to the overall solution is achieved. Each iteration uses the conventional DP algorithm to solve the problem within a strongly reduced state domain around the received state trajectory of the previous iteration. The computational time required for each iteration is far lower compared to that of the one-shot full problem solution, due to the strongly reduced space domain considered. Thus, the problem may be solved by use of DDDP very fast (around 1 s CPU time or less) \citep{typaldos2021vehicle}. It should be noted that, in contrast to deterministic optimal control problems, general stochastic optimal control problems do not feature a solution trajectory, even for given initial states, due to the uncertainty of system evolution created by the included stochastic variables. However, in the specific stochastic GLOSA problem, such a trajectory is indeed present in the problem solution, and this allows for application of the DDDP algorithm.

As an attempt to further reduce the computation time for the solution of the stochastic GLOSA problem and enable its solution on-board the vehicle, a Differential Dynamic Programming (DDP) algorithm is employed in the current work. DDP was first introduced in \citep{jacobson1970differential}, while several extensions and applications were eventually conducted by \cite{murray1979constrained, murray1984differential}. DDP is also used in robotics \citep{tassa2014control, budhiraja2018differential} and unmanned aerial vehicles (UAV) \citep{xie2017differential}. DDP is also an iterative algorithm, which starts with a user-defined feasible trajectory to be used in the first iteration. Each iteration produces an improved trajectory, to be used in the subsequent iteration; until convergence to the solution of the original optimal control problem is achieved. Each iteration solves, for each time-step separately, a quadratic-linear approximation of the stochastic recursive Bellman equation; the approximation being computed around the initial trajectory of the iteration. It should be noted that, in contrast to the standard DP and DDDP algorithms, DDP requires no discretization of the state and control spaces, hence it may result in better solutions compared with discretised-value formulations. Moreover, DDP is found to deliver solutions quasi-instantaneously, which enables the approach to be readily implementable in an MPC (model predictive control) framework online, on the vehicle's computer. Note that communication between vehicles and traffic signals (V2I) is crucial in GLOSA problems \citep{lu2020next}, however in this work we focus on the methodological part of such systems.

The remainder of the paper is organized as follows: in Section \ref{sec:background}, the optimal control problems with known signal switching time and uncertain signal switching time, as proposed in \citep{typaldos2020vehicle} are briefly presented for completeness. Section \ref{sec:dddp} presents the DDDP algorithm, along the lines of \citep{typaldos2021vehicle}. Section \ref{sec:ddp} presents the DDP-based GLOSA problem algorithm, while demonstration results for the DDP algorithm's performance and comparison with the one-shot stochastic GLOSA and the DDDP results are presented in Section \ref{sec:results}. Finally, Section \ref{sec:concl} concludes this work, summarizing its contributions.

\section{Optimal Control-Based GLOSA with Known or Unknown Signal Switching Time} \label{sec:background}
This section outlines, for completeness, the deterministic and stochastic GLOSA optimal control problems, for the cases of known and unknown signal switching time, respectively, as well as the SDP solution algorithm for the latter, following \citep{typaldos2020vehicle}.

\subsection{Known Signal Switching Time - Deterministic GLOSA} \label{sec:back-determ}
\begin{figure}
	\centering
		\includegraphics[width=.5\linewidth]{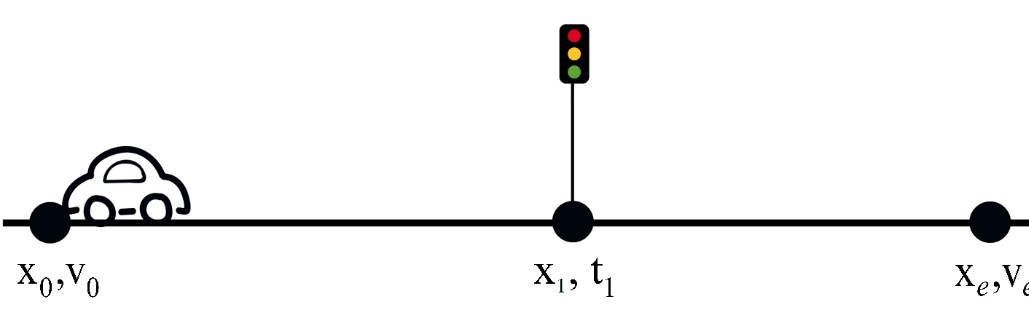}
	\caption{A vehicle starting from a given initial state ($x_0,v_0$) aims to reach a fixed final state ($x_e,v_e$) in presence of a traffic signal at a given position $x_1$ that switches to green at time $t_1$.}
	\label{fig:veh-example}
\end{figure}

Consider a vehicle traveling from an initial state $\bm{x}_0=[x_0,v_0 ]^T$, $x_0$ being a given initial position and $v_0$ a given initial speed, with the purpose to reach a fixed final state $\bm{x}_e=[x_e,v_e ]^T$ within a free (but weighted) time horizon $t_e$. Between the initial and final positions, there is a traffic signal, and hence the additional restriction that the vehicle must not pass the signal position $x_1$ before the known time $t_1$, when the traffic light turns green from red (see Fig. \ref{fig:veh-example}). The implicit assumption here is that a red light is active when the vehicle appears on the link (at time 0), but a generalisation, which includes the case where the vehicle appears when the traffic light is green, is given in \citep{typaldos2020vehicle}. The vehicle task is to appropriately adjust its acceleration (control variable), so as to minimize fuel consumption, while satisfying the initial and final conditions $\bm{x}_0$ and $\bm{x}_e$, as well as the intermediate (traffic signal) constraint.

The minimization problem outlined above is formulated as an optimal control problem, which accounts for the vehicle kinematics via the following state equations:
\begin{align}
\dot{x} &= v \\
\dot{v} &= a
\end{align}
where $a$ is the vehicle acceleration, which assumes the role of the control variable. The objective criterion to be minimized reads
\begin{equation}
J = w \cdot t_e + \dfrac{1}{2} \int_0^{t_e} a^2 dt \label{eq:det-cost}
\end{equation}

The utilized acceleration cost term $a^2$ in the criterion was demonstrated to be an excellent proxy for deriving fuel-minimizing vehicle trajectories \citep{typaldos2020minimization}. The final time $t_e$ is free but penalized with the weight $w$. For higher values of $w$, the resulting $t_e$ will be lower and vice versa. This, consequently, affects the acceleration cost, which, depending on higher or lower $w$ value, will also have increased or decreased values (for more details, see \citep{typaldos2020vehicle}). In addition, the green-light constraint, $t_S \geq t_1$, must be fulfilled, where $t_S$ is the time at which the vehicle crosses the signal position $x_1$, that is, $x(t_S ) = x_1$. If necessary, upper and lower bounds may be applied to speed $v$ and acceleration $a$.
 
The solution of this problem was addressed in \citep{typaldos2020vehicle} and can be obtained analytically using symbolic differentiation tools. Thus, the numerical solution of the deterministic GLOSA problem, for a specific problem instance, takes only fractions of a second of computation time and can be executed in real time on-board for each approaching vehicle. Note that it may be reasonable to continuously update the vehicle trajectory (in a model predictive control (MPC) loop) to account for possible deviations from the first calculated vehicle trajectory, which may occur because of a variety of disturbances, including a slower vehicle ahead.

For a given junction, the final state is the same for any initial vehicle state $\bm{x}_0$ and any switching time $t_1$. Therefore, the optimal value of criterion \eqref{eq:det-cost} of the deterministic GLOSA solution depends only on the initial state and the switching time and is denoted $J_{DG}^*(\bm{x}_0 ,t_1)$ for later use.

\subsection{Uncertain Signal Switching Time Problem – Stochastic GLOSA} \label{sec:back-stoch}
The traffic light switching time may be subject to short-term decisions in dependence of the prevailing traffic conditions in cases of real-time signals. In such cases, we typically have minimum and maximum admissible switching times; hence, based on statistics from past signal switching activity, we may derive a probability distribution of switching times within the admissible time-window of possible signal switching times. Thus, the problem can be cast in the format of a stochastic optimal control problem, which may be solved numerically using SDP techniques. To this end, the analytical solution of the deterministic GLOSA optimal control problem is used within the stochastic approach, as will be explained in this section (for more details, see \citep{typaldos2020vehicle}).

For the stochastic optimal control problem and the SDP algorithm \citep{bertsekas2012dynamic}, the discrete-time version of the vehicle kinematics, with time step $T$, is considered, as follows:
\begin{align}
x(k+1) &= x(k) + v(k)T + \dfrac{1}{2} a(k) T^2 \label{eq:dis-state1}\\
v(k+1) &= v(k) + a(k)T \label{eq:dis-state2}
\end{align}
where $x(k)$ and $v(k)$ correspond to the vehicle position and speed, respectively, at discrete times $k = 0,1,\dots$ (where $kT = t$), while the control variable $a(k)$ reflects the acceleration that remains constant over each time period $k$. The state and control variables are bounded within the following feasible regions
\begin{align}
\bm{x}(k) \in \bm{X} &= [\bm{x}_{min}, \bm{x}_{max}] \label{eq:bounds-x}\\
a(k) \in U &= [a_{min}, a_{max}] \label{eq:bounds-a}
\end{align}
with $\bm{x}_{min}, \bm{x}_{max}$ and $a_{min}, a_{max}$ being the lower and upper bounds of the states and acceleration, respectively. The traffic light's discrete switching time $k_1$ is not known, but it is assumed that a known range $k_{min} \leq k_1 \leq k_{max}$ of possible switching times exists, $k_{min}$ and $k_{max}$ being the minimum and maximum possible switching times, respectively.
 
For proper problem formulation, a virtual state $\tilde{x}(k)$ is introduced, that reflects formally the stochasticity of traffic light switching
\begin{equation} \label{eq:stoch-state}
\begin{aligned}
\tilde{x}(k+1) &= \tilde{x}(k) \cdot z(k) \\
\tilde{x}(0) &= 1
\end{aligned}
\end{equation}
where $z(k)$ is a binary stochastic variable defined as
\begin{equation} \label{eq:stoch-var}
z(k)=
\begin{cases}
0 & \text{if traffic light switches at time } k + 1 \\
1 & \text{else}
\end{cases}
\end{equation}
with \eqref{eq:stoch-state} and \eqref{eq:stoch-var}, the virtual state $\tilde{x}(k)$ is either equal to 1 if the traffic light has not yet switched until time $k - 1$; or equal to 0 if switching occurred at time $k$ or earlier. The virtual state $\tilde{x}$ is assumed measurable, which means that the system knows at each time $kT$ if switching has taken place or not within the last time period $[(k-1)T,kT]$. 

The stochastic variable $z(k)$ is independent of its previous values and takes values according to a time-dependent probability distribution $p(z|k)$. Based on the statistics of previous signal switching activity, availability of an a-priori discrete probability distribution $P(k), k_{min} \leq k_1 \leq k_{max}$, is assumed, for signal switching within the time-window, where $\sum_{k=k_{min}}^{k_{max}} P(k) = 1$. Since no switching takes place for $k \leq k_{min} - 1$, we have
\begin{equation}
p(0|k) = 0 \quad \text{for} \quad k < k_{min} - 1
\end{equation}

For $k \geq k_{min}$, the probability distribution $p(z|k)$, is obtained by use of "crop-and-scale", meaning that the a-priori probabilities of previous time steps, where switching did not take place, are distributed analogously to increase the probabilities of the remaining discrete times within the time-window \citep{lawitzky2013energy}. As shown in \citep{typaldos2020vehicle}, this update may be done by use of the following crop-and-scale formula that applies for $k_{min} \leq k \leq k_{max} - 1$ and for any a-priori distribution $P(k)$
\begin{equation}
p(0|k) = \dfrac{P(k + 1)}{\sum_{\kappa = k+1}^{k_{max}} P(\kappa)}
\end{equation}

The cost criterion of the stochastic problem is essentially the same as \eqref{eq:det-cost} in the deterministic case. However, in the stochastic case, the exact value of the criterion depends on the stochastic variable's realization, and therefore we consider minimisation of the expected value. As shown in \citep{typaldos2020vehicle}, the cost criterion in discrete time becomes
\begin{equation}
J = E \Bigg\{ \frac{1}{2} \sum_{k=0}^{k_1 - 1} a(k)^2 + J_{DG}^*[\bm{x}(k_1), k_1] \Bigg\}
\end{equation}
where the expectation refers to the stochastic variable $z(k), k = 0,\dots,k_{max} - 1$. This criterion captures the quadratic acceleration cost up until the stochastic switching time $k_1$, when the vehicle is at state $\bm{x}(k_1)$. After the signal switching, the problem instantly becomes a deterministic GLOSA problem (Section \ref{sec:back-determ}), and the corresponding optimal cost-to-go is $J_{DG}^*[\bm{x}(k_1), k_1 ]$, which will be denoted as the "escape cost".

To obtain a formally proper cost criterion, the stochastic variable $z(k)$ and the virtual variable $\tilde{x}(k)$, introduced earlier, are used, and, as shown in \citep{typaldos2020vehicle}, this yields the objective function in the required form, as follows
\begin{equation} \label{eq:stoch-cost}
\begin{aligned}
J = E \Bigg\{ \tilde{x}(k) & \sum_{k=0}^{k_{max} - 1} \bigg[ \frac{1}{2} a(k)^2 + [1-z(k)] J_{DG}^* [\bm{x}(k), a(k), k + 1] \bigg] \Bigg\}
\end{aligned}
\end{equation}
where, for notational simplicity, $J_{DG}^*[\bm{x}(k),a(k),k+1] = J_{DG}^* [\bm{x}(k+1),k+1]$ when $\bm{x}(k+1)$ is replaced via the state equations \eqref{eq:dis-state1}, \eqref{eq:dis-state2} as a function of $\bm{x}(k),a(k)$. 

Equations \eqref{eq:dis-state1}-\eqref{eq:stoch-cost} constitute an ordinary stochastic optimal control problem \citep{bertsekas2012dynamic}. Denoting the corresponding optimal cost-to-go function by $V[\bm{x}(k),\tilde{x}(k),k]$, the recursive Bellman equation for $0 \leq k \leq k_{max} - 1$ reads
\begin{align} \label{eq:rec-belm-eq}
V&[\bm{x}(k), \tilde{x}(k), k] = \min_{a(k) \in U} \bigg\{ E \Big\{\frac{1}{2} a(k)^2 + [1-z(k)] J_{DG}^* [\bm{x}(k), a(k), k + 1] \nonumber \\
& + V[\bm{x}(k+1), \tilde{x}(k)z(k), k+1] \Big\} \bigg\} = \min_{a(k) \in U} \bigg\{ \frac{1}{2} a(k)^2 + p(0|k) \cdot J_{DG}^* [\bm{x}(k), a(k), k + 1]  \\
& + [1 - p(0|k)] \cdot V[\bm{x}(k+1), 1, k+1] \bigg\} \nonumber
\end{align}
with boundary condition $V[\bm{x}(k_{max}), 1, k_{max}] = 0$. Note that, the minimum is required with respect to the control of time $k$ only, i.e., $a(k) \in U$, as typical in dynamic programming, which facilitates the numerical solution.

In this formulation, it has been assumed that the decision on the traffic light switching is taken within the last time period before the actual switching. The generalization to the case of taking a switching decision $\kappa$ time periods in advance of the actual switching is easy and is mentioned in \citep{typaldos2020vehicle}.

\subsection{Discrete SDP Numerical Solution Algorithm} \label{sec:back-sol}
For the numerical solution of the stochastic problem, using the SDP algorithm, the state and control variables must be discretised. As the discretization level has a significant impact on computational time and memory requirements, but also on the accuracy of the computed solution, an appropriate trade-off must be specified between reasonable computation requirements versus achievable solution quality.
For the discretization, the discrete time interval $T$ is set equal to 1 s, which is a reasonable choice for the problem at hand. Then, a general discretization interval $\Delta$ for the problem variables is assumed, and the discretization interval of acceleration is set $\Delta a = \Delta$. From \eqref{eq:dis-state2}, the discretization interval of speed assumes the same value ($\Delta v = \Delta a = \Delta$). Likewise, in view of \eqref{eq:dis-state1}, the discretization interval for the position is 
\begin{equation} \label{eq:Delta-x}
\Delta x = \dfrac{1}{2} \Delta \cdot T^2
\end{equation}

Based on these settings, it was shown in \citep{typaldos2020vehicle} that, if $x(k),v(k),a(k)$ are discrete points, then $x(k+1)$ and $v(k+1)$ resulting from \eqref{eq:dis-state1} and \eqref{eq:dis-state2} are also discrete points.

It is now straightforward to apply the discrete SDP algorithm to obtain a globally optimal closed-loop control law $a^*(k) = R[\bm{x}(k),k]$, which, for any given vehicle state $\bm{x}(k) \in \bm{X}$ and time $k$, delivers the corresponding optimal acceleration $a^*(k)$. The SDP algorithmic steps are summarized below. Note that the algorithm needs to consider only the case $\tilde{x}(k) = 1$, therefore any arguments pertaining to $\tilde{x}(k)$ are suppressed for convenience.
 
The SDP algorithm is as follows:

\begin{algorithmic}[1]
\State $V[\bm{x}(k_{\max}),k_{\max}) ] \gets 0 \quad \forall \bm{x}(k_{\max}) \in \bm{X}$
\For{each $k = k_{\max} - 1, \dots , 0$}
		\For{each discrete state $\bm{x}(k) \in \bm{X}$}
				\For{each discrete control $a(k) \in U$}
						\State Calculate $x(k + 1), v(k + 1)$
						\If {$\bm{x}(k + 1) \not\in \bm{X}$}
								\State $\Phi[\bm{x}(k), a(k), k + 1] \gets \infty$
								\State \textbf{continue}
						\EndIf
				\State $\Phi[\bm{x}(k), a(k), k + 1] \gets \frac{1}{2} a(k)^2 + p(0|k) \cdot J_{DG}^*[\bm{x}(k), a(k), k + 1]$ \par
				\hskip\algorithmicindent \hskip\algorithmicindent $+ [1 - p(0|k)] \cdot V[\bm{x}(k+1), k+1]$
				\EndFor
				\State $V[\bm{x}(k), k] \gets \min{\Phi[\bm{x}(k), a(k), k + 1]} \quad \forall a(k) \in U$
				\State $R[\bm{x}(k), k] = a^*(k) \gets \text{arg}\min\limits_{a(k) \in U}\{\Phi[\bm{x}(k), a(k), k + 1]\}$,
				\State with $a^*(k)$ the optimal control of point $[\bm{x}(k), k]$
		\EndFor
\EndFor
\end{algorithmic}

As mentioned, this algorithm delivers a globally optimal control law $a^*(k) = R[\bm{x}(k),k]$ for the full state domain $\bm{X}$. Note that general stochastic optimal control problems do not possess a solution trajectory for specific initial states $\bm{x}_0$, because the state evolution is uncertain in presence of the stochastic variables. However, in the specific stochastic GLOSA problem addressed here, the evolution of the vehicle states \eqref{eq:dis-state1} and \eqref{eq:dis-state2} is not affected by the stochastic variable $z(k)$, which concerns only the signal switching time. Thus, for a given initial state, i.e., vehicle position and speed at time 0, the optimal control law may be used to produce an optimal vehicle trajectory that the vehicle should pursue; until the signal switching actually occurs, in which case the vehicle must instantly behave according to the deterministic GLOSA solution.

Note that the SDP algorithm may take several minutes to execute in an ordinary PC, as will be reported in Section \ref{sec:results}. This implies that the solution cannot be obtained on the vehicle side in real time.

\section{Discrete Differential Dynamic Programming (DDDP)} \label{sec:dddp}
In \citep{typaldos2021vehicle}, a modified DP algorithm called DDDP was employed for the stochastic GLOSA problem in order to address the major disadvantage of the discrete (S)DP algorithm, which is the high computation time required for the numerical solution of the optimal control problem, which increases exponentially with the number of state and control variables \citep{papageorgiou2015optimierung}. The DDDP algorithm was proposed in \citep{heidari1971discrete} for deterministic optimal control problems. The method can nevertheless be applied here, because an optimal vehicle trajectory may be derived for the stochastic GLOSA problem, despite its stochastic character, as mentioned in the previous section. The DDDP approach to the stochastic GLOSA problem is presented in this paper in more detail and with more comprehensive testing and comparison results, compared to the conference article \citep{typaldos2021vehicle}.

DDDP is an iterative algorithm, calling for a feasible starting state trajectory to be specified externally. Each iteration $l$ receives a feasible (but non-optimal) state trajectory $\bm{x}^{(l-1)}(k)$, and transforms it to an enhanced one $\bm{x}^{(l)}(k)$, to be used in the next iteration. To this end, each iteration solves a discrete SDP problem by use of the standard SDP algorithm presented above. What changes at each iteration $l$ is the considered state domain 
\begin{equation}
\bm{X}_C^{(l)} = \Big\{ \bm{x}(k) \: \big| \: |\bm{x}(k) - \bm{x}^{(l-1)}(k)| \leq \bm{\Delta}_C^{(l)} \wedge \bm{x}(k) \in \bm{X} \Big\}
\end{equation}
which is a strongly reduced subdomain of the original state domain $\bm{X}$ in \eqref{eq:bounds-x}. In other words, the discretized SDP problem is solved in each iteration $l$ within a corridor with width $\bm{\Delta}_C^{(l)}$ around the received state trajectory $\bm{x}^{(l-1)}(k)$ of the previous iteration, to produce a solution trajectory $\bm{x}^{(l)}(k)$ for use in the next iteration. The corridor width $\bm{\Delta}_C^{(l)}$, as well as the discretisation intervals $\Delta a^{(l)}, \Delta \bm{x}^{(l)}$, can vary in each iteration, typically at a decreasing rate. The iterations stops when a termination criterion is satisfied.

In the proposed GLOSA application \citep{typaldos2021vehicle}, the starting trajectory for the first iteration of DDDP is the optimal solution of the deterministic GLOSA problem, assuming the "pessimistic" case where the traffic light will switch from red to green at the latest possible time, that is, at $k_1 = k_{max}$, so as to cover the whole time range and be not too far from the stochastic optimal trajectory; see \citep{typaldos2020vehicle}. The discretisation intervals $\Delta a^{(l)}, \Delta \bm{x}^{(l)}$ for the first iteration are also given. The intervals are reduced by half, each time there is no solution improvement at two subsequent iterations. The corridor width is taken as 
\begin{equation} \label{eq:Delta-C}
\bm{\Delta}_C^{(l)} = \bm{C} \cdot \Delta a^{(l)}
\end{equation} 
where $\bm{C} = [C_x, C_v]$ are constant given values specifying the state corridor width. Thus, the corridor's initial width is defined through the chosen values for $\bm{C}$ and $\Delta a^{(l)}$. In the following iterations, whenever the discrete interval is reduced, there is an analogous reduction of the corridor width. Consequently, we have a constant number of feasible discrete points in all iterations, which facilitates the algorithm's fine-tuning, so as to improve its computational efficiency. Note that, $C_x$ and $C_v$ values may differ, as the magnitude of the two respective state variables differs.

The admissible control region $U$ (see \eqref{eq:bounds-a}) is kept the same at each iteration, although many related state transitions cannot be considered in view of the reduced state variables domain. The algorithm terminates whenever there is no improvement of the produced solution trajectory even after a reduction of the discretisation interval; or when the $\Delta a^{(l)}$ value becomes less than 0.125, which was found in \citep{typaldos2020vehicle} to lead to sufficiently accurate results.

It should be noted that there is no general guarantee that the DDDP iterations will actually converge to the full-domain SDP solution. In particular, if the state sub-domains $\bm{X}_C^{(l)}$ considered in the iterations are too small, the obtained DDDP solution may actually differ from the overall SDP solution. On the other hand, if the state sub-domains are selected large, the required number of iterations may decrease, but the computation time required to find the solution at each iteration increases accordingly. In conclusion, some fine-tuning regarding the size of sub-domains is necessary to ensure convergence to the overall SDP solution with minimum overall (all iterations) computation time.

\section{Differential Dynamic Programming (DDP)} \label{sec:ddp}
DDP is an iterative algorithm, which requires a given initial feasible trajectory (of states and controls) to be used in the first iteration. At each iteration and time step, a quadratic-linear approximation of the recursive Bellman equation \eqref{eq:rec-belm-eq}, around the initial trajectory of the iteration, is solved, leading to an improved trajectory. The improved trajectory is then used as initial trajectory in the next iteration, and the procedure continues until convergence to the solution of the original optimal control problem is achieved. The method was developed for deterministic optimal control problems; and it is not applicable to general stochastic problems, which do not feature solution trajectories \citep{bertsekas2012dynamic}. However, as explained earlier, the specific stochastic GLOSA problem admits feasible solution trajectories for specific initial states, hence DDP may be applied.

\subsection{General Constrained DDP} \label{sec:ddp-general}
The standard DDP method \citep{jacobson1970differential} does not account for constraints, so, for the needs of the problem presented in this work, extensions of DDP, proposed by \cite{murray1979constrained} for linear constraints and by \cite{yakowitz1986stagewise} for nonlinear constraints, is utilized. 

Consider a discrete-time optimal control problem with states $\bm{x}$ and controls $\bm{a}$ with the recursive Bellman equation $V(\bm{x}(k),k)=\min\{\phi(\bm{x},\bm{a},k) + V(\bm{f}(\bm{x},\bm{a},k),k+1)\}$, for $k=K-1,\dots,0$, where $V$ is the optimal cost function, $\phi$ is the objective function of the optimal control problem, $\bm{f}$ is the right-hand side of the state equation, and the minimisation must be carried out for all feasible controls $\bm{a}(k)$ that satisfy any inequality constraints present. The DDP procedure undertakes, at each time step of each iteration, a quadratic approximation of the term to be minimized in the recursive Bellman equation, i.e., of $\phi(\bm{x},\bm{a},k)+V(\bm{f}(\bm{x},\bm{a},k),k+1)$. The quadratic approximation is taken around nominal trajectories ($\bm{\bar{x}}(k), \bm{\bar{a}}(k)$), which are the initial trajectories of each iteration. Define $\delta \bm{x} = \bm{x}- \bm{\bar{x}}$ and $\delta \bm{a} = \bm{a} - \bm{\bar{a}}$. The aim is to find the optimal control law for $\delta \bm{a}$, which minimizes the quadratic approximation, subject to some inequality constraints. The procedure of the DDP algorithm for each iteration is described in brief as follows, see \citep{yakowitz1986stagewise, murray1979constrained} for more details. 

For $k = K-1, \dots, 0$, the quadratic approximation $Q(\bm{x},\bm{a},k)$ of the term $\phi(\bm{x},\bm{a},k)+V(\bm{f}(\bm{x},\bm{a},k),k+1)$ around a nominal point ($\bm{\bar{x}}(k),\bm{\bar{a}}(k)$), is required and can be expressed as follows (for simplicity, the indexes $k$ are omitted)
\begin{equation} \label{eq:quadratic}
Q(\bm{x},\bm{a}) = \dfrac{1}{2} \delta\bm{x}^T \bm{A} \delta\bm{x}+\delta\bm{x}^T \bm{B} \delta\bm{a} + \dfrac{1}{2} \delta\bm{a}^T \bm{C} \delta\bm{a} +\bm{D}^T\delta\bm{a}+\bm{E}^T\delta\bm{x}
\end{equation}
where the zero-order constant term is omitted, as it does not influence the minimisation; and $\bm{A}, \bm{B}, \bm{C}, \bm{D}, \bm{E}$ are matrices with appropriate dimensions and with elements that may depend on ($\bm{\bar{x}}(k), \bm{\bar{a}}(k)$).

Furthermore, we assume the general inequality constraints for the optimal control problem
\begin{equation}
\bm{c}\left(\bm{x},\bm{a}\right) \leq \bm{0} 
\end{equation}
and we denote as $\widetilde{\bm{c}}\left(\bm{x},\bm{a}\right)$ its linear approximation around the nominal trajectories $(\bm{\bar{x}},\bm{\bar{a}})$.

For the solution of the resulting Quadratic Programming (QP) problem, Fletcher's active set algorithm is used \citep{fletcher1971general}. The critical issue here is that the minimization in the recursive Bellman equation should be executed in a parametrized way, i.e., deliver $\delta \bm{a}$ as a function of $\delta \bm{x}$, which is not compatible with a numerical solution procedure. Therefore, in the QP problem formulation, we replace $\bm{x}$ with $\bm{\bar{x}}$ (i.e. we set $\delta \bm{x} = 0$) in \eqref{eq:quadratic} and $\tilde{c}(\bm{x}, \bm{a})$, and solve the following QP problem to obtain a numerical value for $\delta \bm{a}$
\begin{equation}
\begin{aligned}
&\min{Q\left({\bar{\bm{x}}}, \bar{\bm{a}}+\delta \bm{a}\right)}	\\
\text{s.t}&	\\
&\widetilde{\bm{c}}\left({\bar{\bm{x}}}, \bar{\bm{a}}+\delta \bm{a}\right)\le \bm{0}
\end{aligned}
\end{equation}

After solving the QP problem, the index set of the active constraints, denoted as $S$, leads to the linear system of equations
\begin{equation} \label{eq: activeConsSet}
{\widetilde{c}}_i\left({\bar{\bm{x}}},\bar{\bm{a}}+\delta \bm{a}\right)= 0,\ \ \ i\in S	
\end{equation}
and expressing \eqref{eq: activeConsSet} in matrix form, we have 
\begin{equation}
\bm{U} \delta\bm{a}-\bm{W}=\bm{0}	
\end{equation}
where $\bm{U}$ has rank equal to the number of indexes in $S$. This solution of the QP problem satisfies the necessary conditions of optimality for the Lagrangian \citep{fletcher1971general}
\begin{equation} \label{eq:lagrange}
L\left(\bm{a},\bm{\lambda}\right)=\frac{1}{2}\delta\bm{a}^T \bm{C} \delta\bm{a}+\bm{D}^T\delta\bm{a}+\bm{\lambda}^T\left(\bm{U}^T\delta\bm{a}-\bm{W}\right)
\end{equation}	
where $\bm{\lambda}$ is the Lagrange multiplier vector. These necessary conditions may be obtained by differentiating $L$ with respect to $\delta \bm{a}$ and $\bm{\lambda}$, which yields the following linear system of equations
\begin{equation} \label{eq:langSystem}
\left[\begin{matrix}\bm{C}&\bm{U}\\\bm{U}^T&\bm{0}\\\end{matrix}\right]\left[\begin{matrix}\delta\bm{a}\\\bm{\lambda}\\\end{matrix}\right]=\left[\begin{matrix}-\bm{D}\\\bm{W}\\\end{matrix}\right]	
\end{equation}

In order to re-introduce $\delta \bm{x}$, i.e. the variablity of state variables, we consider that the constraints that are active while we assumed $\bm{x}=\bar{\bm{x}}$, remain active for other $\bm{x}$, i.e.
\begin{equation}
{\widetilde{c}}_i\left(\bar{\bm{x}} + \delta \bm{x},\bar{\bm{a}} + \delta \bm{a}\right)=0,\ \ \ i\in S
\end{equation}	
which in matrix form reads
\begin{equation}
\bm{U}\delta\bm{a}-\bm{W}-\bm{X}\delta\bm{x}=\bm{0}	
\end{equation}

Then, we may derive a control law, i.e. $\delta \bm{a}$ as a function of $\delta \bm{x}$ by solving analytically the accordingly extended version of \eqref{eq:langSystem}, i.e.
\begin{equation}
\left[\begin{matrix}\bm{C}&\bm{U}\\\bm{U}^T&\bm{0}\\\end{matrix}\right]\left[\begin{matrix}\delta\bm{a}\\\bm{\lambda}\\\end{matrix}\right]=\left[\begin{matrix}-\bm{D}-\bm{B}\delta\bm{x}\\\bm{W}+\bm{X}\delta\bm{x}\\\end{matrix}\right]	
\end{equation}
which yields the following linear control law with readily computed vector $\bm{\alpha}$ and matrix $\bm{\beta}$ 
\begin{equation} \label{eq:ctr-law}
\delta\bm{a}\left(\bm{x}\right)=\bm{\alpha}+\bm{\beta}\delta\bm{x}	
\end{equation}

This minimisation outcome is substituted in the quadratic approximation $Q$ to obtain, from the Bellman equation, the approximate optimal cost function
\begin{equation}
V\left(\bm{x}(k),k\right)=Q\left(\bm{x}(k),\bm{\alpha}(k)+\bm{\beta}(k)\delta\bm{x}(k), k\right)	
\end{equation}
Which must be stored for use at the next time step $k-1$ within the same iteration. Note that this approximate cost function is quadratic by construction, which facilitates the quadratic approximation in the Bellman equation at the next time step $k-1$.

Given $\delta \bm{a}(\delta \bm{x},k)$, from \eqref{eq:ctr-law}, the nominal trajectories for the next iteration are calculated as follows for $k=0,\dots, K-1$
\begin{align}
\bm{a}(k) &= \bar{\bm{a}}(k) + \epsilon[\bm{\alpha}(k) + \bm{\beta}(k)\delta\bm{x}(k)] \label{eq:forward-1}\\
\bm{x}(k+1) &= f(\bm{x}(k),\bm{a}(k),k) \label{eq:forward-2}\\
\bm{x}(0) &= \bm{x}_0 \label{eq:forward-3}
\end{align}
where $ 0 < \epsilon \leq 1$ is used to ensure the convergence of the algorithm. Specifically, $\epsilon$ values smaller than 1 may facilitate convergence if needed; but may decrease the convergence rate if selected too small. The algorithm terminates whenever the following convergence test is satisfied
\begin{equation}\label{eq:terminal}
\bigg(\sum_{k=0}^{K-1}\big[\big(\bm{a}(k) - \bar{\bm{a}}(k)\big)^2\big]\bigg)^{1/2} < \epsilon_1
\end{equation}
where $\epsilon_1$ is a small positive number.

In contrast to the SDP and DDDP algorithms, DDP does not apply any discrete DP procedures, hence no discretization of the state and control domains is necessary. As a result, the required computation time is independent of any discretization intervals, and the obtained solution trajectories are value-continuous, hence potentially better than those obtained by the other algorithms, even when using fine discretization.

\subsection{Constrained DDP for the GLOSA Problem} \label{sec:ddp-glosa}

To apply the DDP algorithm to the GLOSA stochastic optimal control problem described in Section \ref{sec:back-stoch}, the quadratic approximation of the right-hand side of \eqref{eq:rec-belm-eq} around a nominal trajectory $(\bm{\bar{x}}(k), \bar{a}(k))$ is first required. Note that only the $J_{DG}^*$ term needs to be approximated, as the rest terms are already quadratic. For the approximation, the LOESS (locally estimated scatterplot smoothing) method was utilized \citep{cleveland2017local}. Specifically, a number of $(x,v)$ points are specified around the nominal trajectory at each iteration, and the corresponding values of $J_{DG}^*(\bm{x}(k),kT)$ are computed. It was found that some 4 such points per time step are sufficient for a good fit. Then, LOESS fits the following smooth function to the data set
\begin{equation} \label{eq:fit}
p(x, v) = p_1 \delta x^2 + p_2 \delta v^2 + p_3 \delta x \delta v + p_4 \delta x + p_5 \delta v + p_6
\end{equation}
where $p_1, \dots, p_6$ are the coefficients of fitting. The resulting fitted function \eqref{eq:fit}, along with the other quadratic terms of \eqref{eq:rec-belm-eq}, can be easily expressed in the form of \eqref{eq:quadratic}. It was verified that the resulting quadratic function \eqref{eq:quadratic} is always convex, something that was expected from the physical background of the problem.

The nominal trajectory of the first iteration is taken to be the optimal solution of the deterministic GLOSA problem, assuming the "pessimistic" case where the traffic light switches from red to green at the latest possible time, i.e., $k_1 = k_{max}$. 

For the state constraints \eqref{eq:bounds-x}, we may replace the state equations \eqref{eq:dis-state1}, \eqref{eq:dis-state2} therein, hence, for the admissible region of the SDP \eqref{eq:bounds-x} and \eqref{eq:bounds-a}, we have the equivalent constraints
\begin{equation} \label{eq:ctr-constr}
\begin{aligned}
\dfrac{2(x_{min}-x(k)-v(k)T)}{T^2} &\leq a(k) \leq \dfrac{2(x_{max}-x(k)-v(k)T)}{T^2} \\
\dfrac{v_{min}-v(k)}{T} &\leq a(k) \leq \dfrac{v_{max}-v(k)}{T} \\
a_{min} &\leq a(k) \leq a_{max}
\end{aligned}
\end{equation}
which can be written in terms of $\delta x$ and $\delta a$ as follows
\begin{equation} \label{eq:ctr-constr-mat}
\begin{bmatrix*}[r]
1 \\ -1 \\ 1 \\ -1 \\ 1 \\ -1
\end{bmatrix*} 
\delta a -
\begin{bmatrix}
\dfrac{2(x_{max}-\bar{x}(k+1))}{T^2} \\ 
\dfrac{-2(x_{min}-\bar{x}(k+1))}{T^2} \\
\dfrac{(v_{max}-\bar{v}(k+1))}{T} \\
\dfrac{-(v_{min}-\bar{v}(k+1))}{T} \\
(a_{max} - \bar{a}(k)) \\
-(a_{min} - \bar{a}(k))
\end{bmatrix}
- \delta \bm{x} \leq
\begin{bmatrix}
0 \\ 0 \\ 0 \\ 0 \\ 0 \\ 0
\end{bmatrix}
\end{equation}

Having the quadratic approximation of the recursive Bellman equation and the linear inequality constraints \eqref{eq:ctr-constr}, the resulting QP is as described in the previous subsection. In fact, the resulting QP problem for the stochastic GLOSA application is simpler than the general case of Section \ref{sec:ddp-general}. After setting $\delta \bm{x} = 0$ in \eqref{eq:ctr-constr-mat}, we have only one decision variable, $\delta a$, with constant upper and lower bounds. Thus, for the solution of the QP problem, it suffices to first compute analytically the value of $\delta a$ that minimises the quadratic function. If this value does not violate any of the upper or lower bounds, then the QP solution is given by that value; else, the QP solution equals the activated bound. All further procedures are as explained in Section \ref{sec:ddp-general}.

\section{Demonstration Results and Comparison} \label{sec:results}

In this section, the performance of the proposed approaches is tested, demonstrated and compared. Specifically, the DDDP and DDP algorithms are applied to realistic scenarios and evaluated regarding the quality of their solution compared to the standard SDP, as well as the computational time needed for each algorithm to converge to the optimal solution.

\subsection{Scenarios Setup} \label{sec:results-setup}

For the evaluation, several scenarios have been considered in various situations, such as cases where the vehicle needs to accelerate or decelerate before or after the traffic signal; or cases where the vehicle is forced to fully stop and wait until the traffic light switches from red to green. For this paper, three scenarios have been chosen for demonstration with full illustration of state and control trajectories, according to Table 1. For consistency, the chosen scenarios are the same as in previous related works \citep{typaldos2020vehicle, typaldos2021vehicle}. All scenarios concern mixed acceleration-deceleration cases. However, in Scenario 1, the vehicle starts with low initial speed and needs to mostly accelerate in order to optimally pass the traffic signal. On the other hand, in Scenarios 2 and 3, the vehicle starts at higher speed, hence it needs to mostly decelerate to avoid a full stop at the traffic signal. Compared to Scenarios 1 and 2, in Scenario 3 the vehicle starts closer to the traffic light, which leads to a full stop for few time steps before the traffic signal switches from red to green.

\begin{table}[!htb]
\centering
\caption{Investigated scenarios.}
\begin{tabular}{cccccc}
\hline
			& $x_0$ (m) & $v_0$ (m/s) 	& $x_e$ (m) & $v_e$ (m/s) 	& $x_1$ (m) \\\hline
Scenario 1  & 0			& 5   			& 220 		& 11			& 150 		\\
Scenario 2  & 0			& $\bm{11}$		& 220 		& 11 			& 150 		\\
Scenario 3  & $\bm{50}$	& $\bm{11}$		& 220 		& 11 			& 150 		\\
\hline
\end{tabular}
\label{tb:scenarios}
\end{table}

For all evaluated scenarios, the weight parameter $w$ in \eqref{eq:det-cost} is equal to 0.1. The bounds for the states and control are set to $[x_{min}, x_{max}]=[0,150]$ m, $[v_{min}, v_{max}]=[0,16]$ m/s and $[a_{min}, a_{max}]=[-3,3]$ m/s$^2$, respectively. The time step $T$ is 1 s and the switching time window for the traffic signal is $[k_{min}, k_{max}]=[10,30]$ s with uniform a-priori probability distribution. Regarding DDDP, for the initial discretization, $\Delta a = \Delta v = 0.5$ is used, which leads to $\Delta x = 0.25$ m according to \eqref{eq:Delta-x}; and the initial corridor width in \eqref{eq:Delta-C} is set $\bm{C} = [20,4]$, i.e., we have $C_x = 5 \cdot C_v$. The choice of the initial discretization and corridor width value will be justified later in this section. For the DDP algorithm, parameters $\epsilon$ in \eqref{eq:forward-1} and $\epsilon_1$ in \eqref{eq:terminal} are set equal to 1 and 10$^{-4}$, respectively.

\subsection{DDDP algorithm results} \label{sec:results-dddp}

\begin{figure}[h!t]
\centering
\begin{minipage}{\linewidth}
\centering
	\begin{subfigure}{0.30\linewidth} 
		\includegraphics[width=\textwidth]{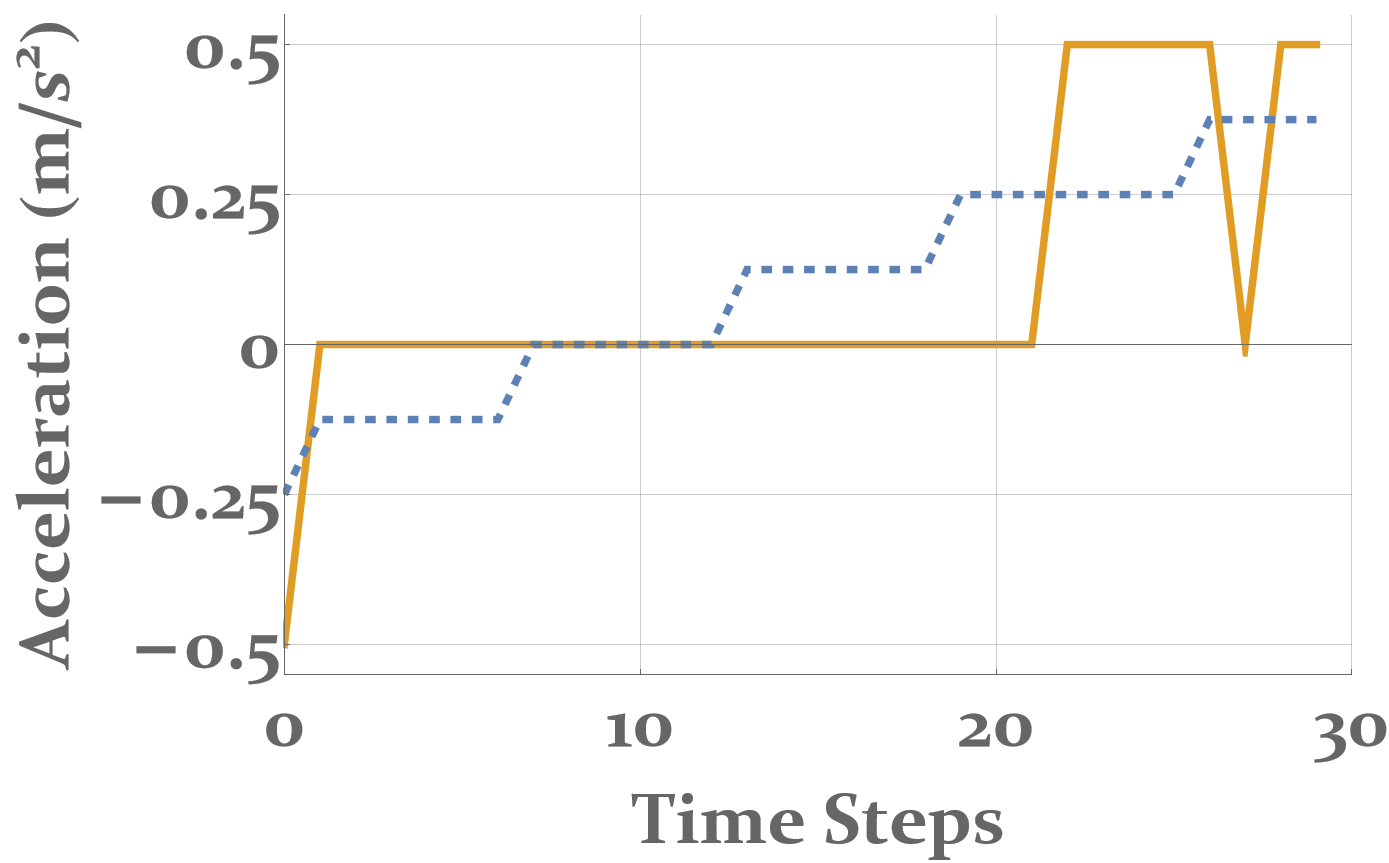}
		\caption*{Iteration 1}
	\end{subfigure}%
	\begin{subfigure}{0.30\linewidth} 
		\includegraphics[width=\textwidth]{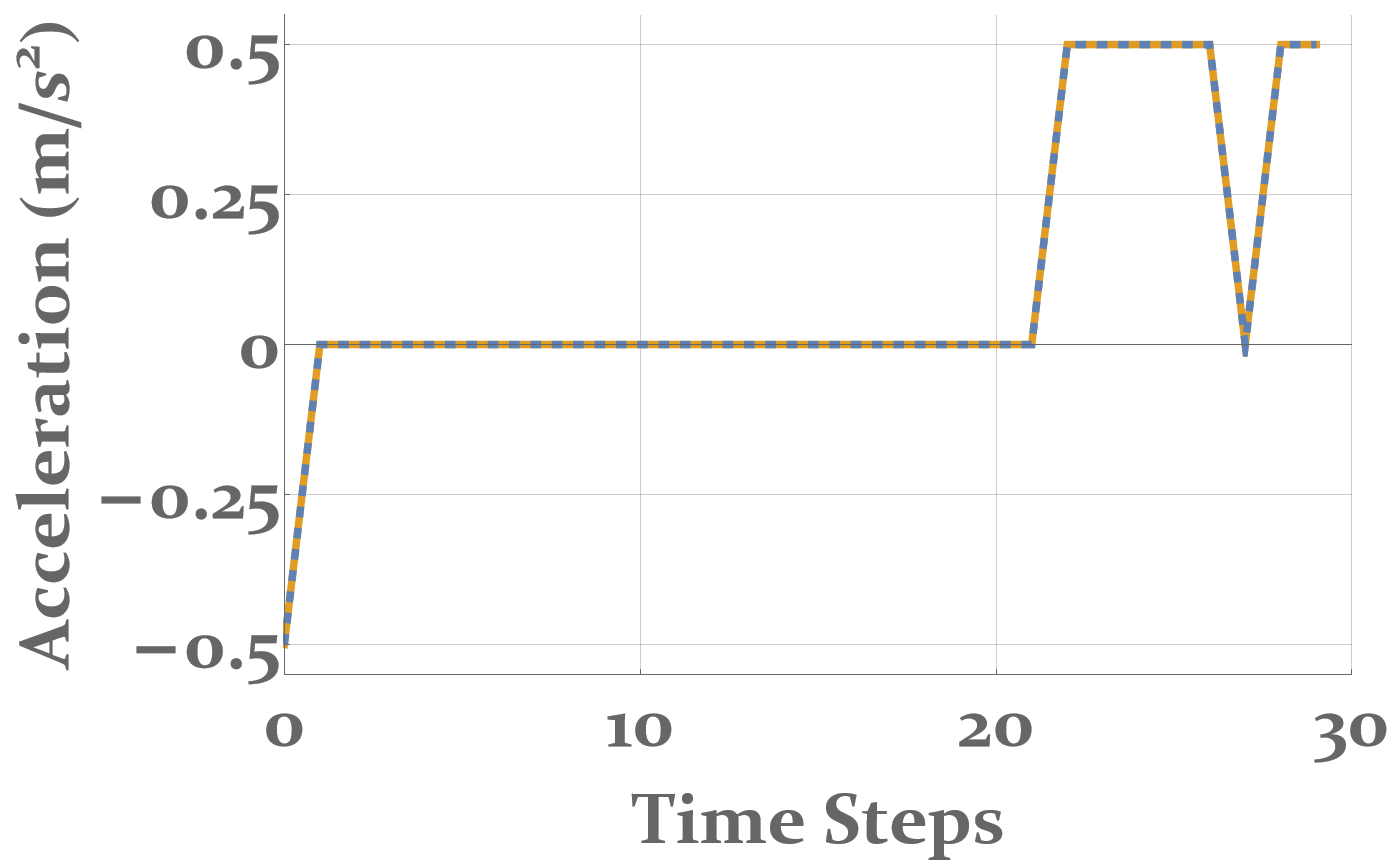}
		\caption*{Iteration 2}
	\end{subfigure}%
	\begin{subfigure}{0.30\linewidth}
		\includegraphics[width=\textwidth]{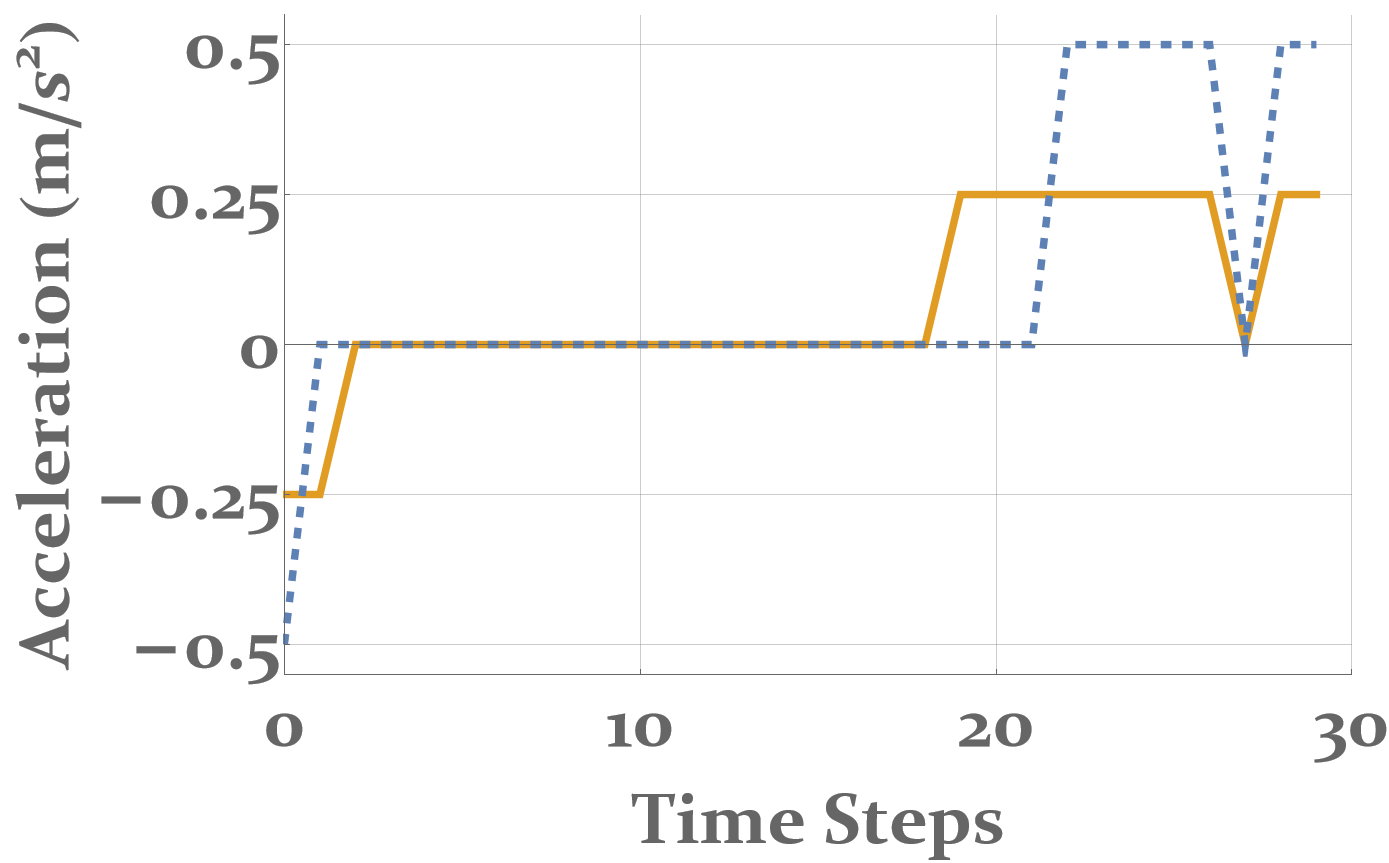}
		\caption*{Iteration 3}
	\end{subfigure}
	\begin{subfigure}{0.30\linewidth}
		\includegraphics[width=\textwidth]{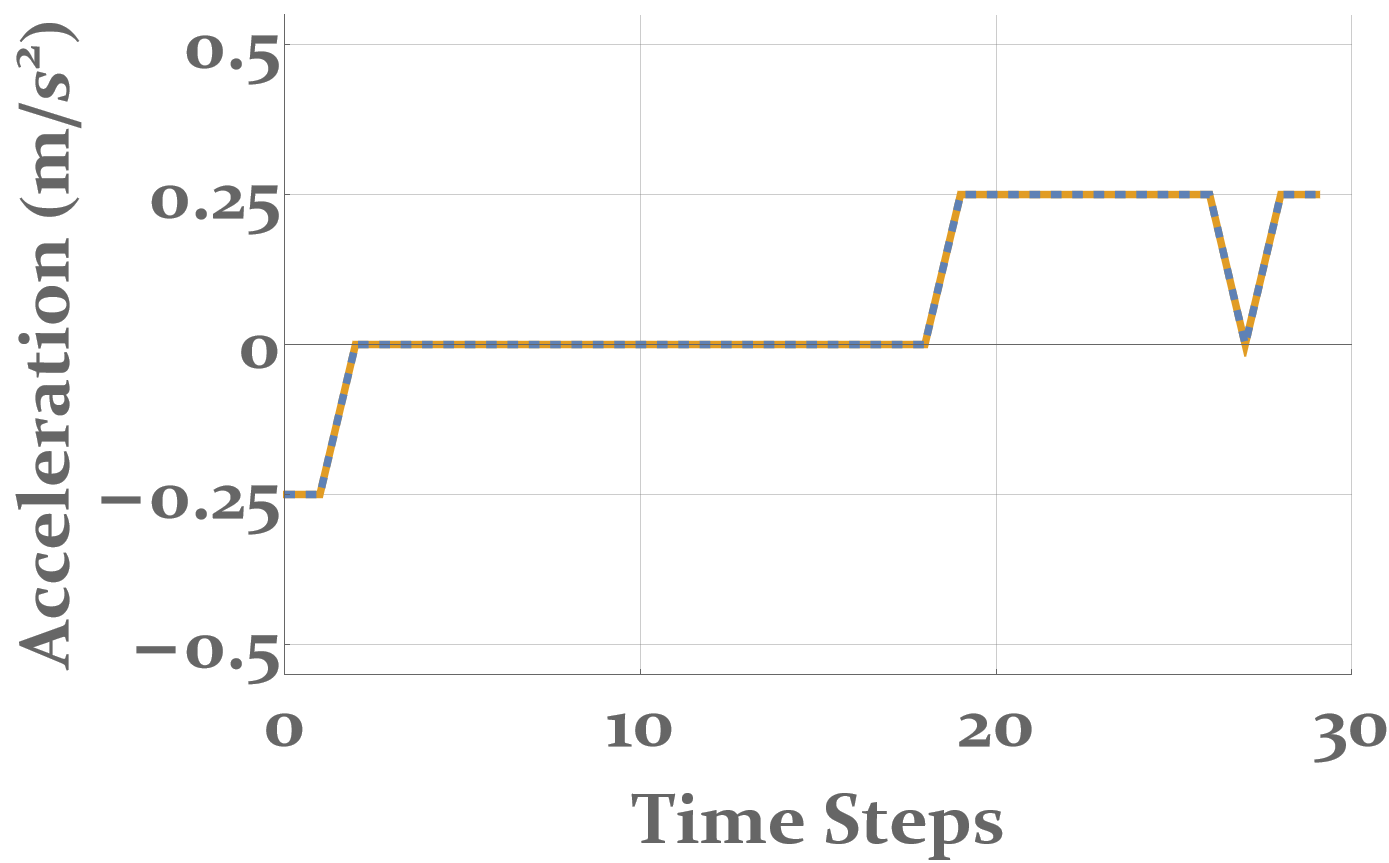}
		\caption*{Iteration 4}
	\end{subfigure}%
	\begin{subfigure}{0.30\linewidth}
		\includegraphics[width=\textwidth]{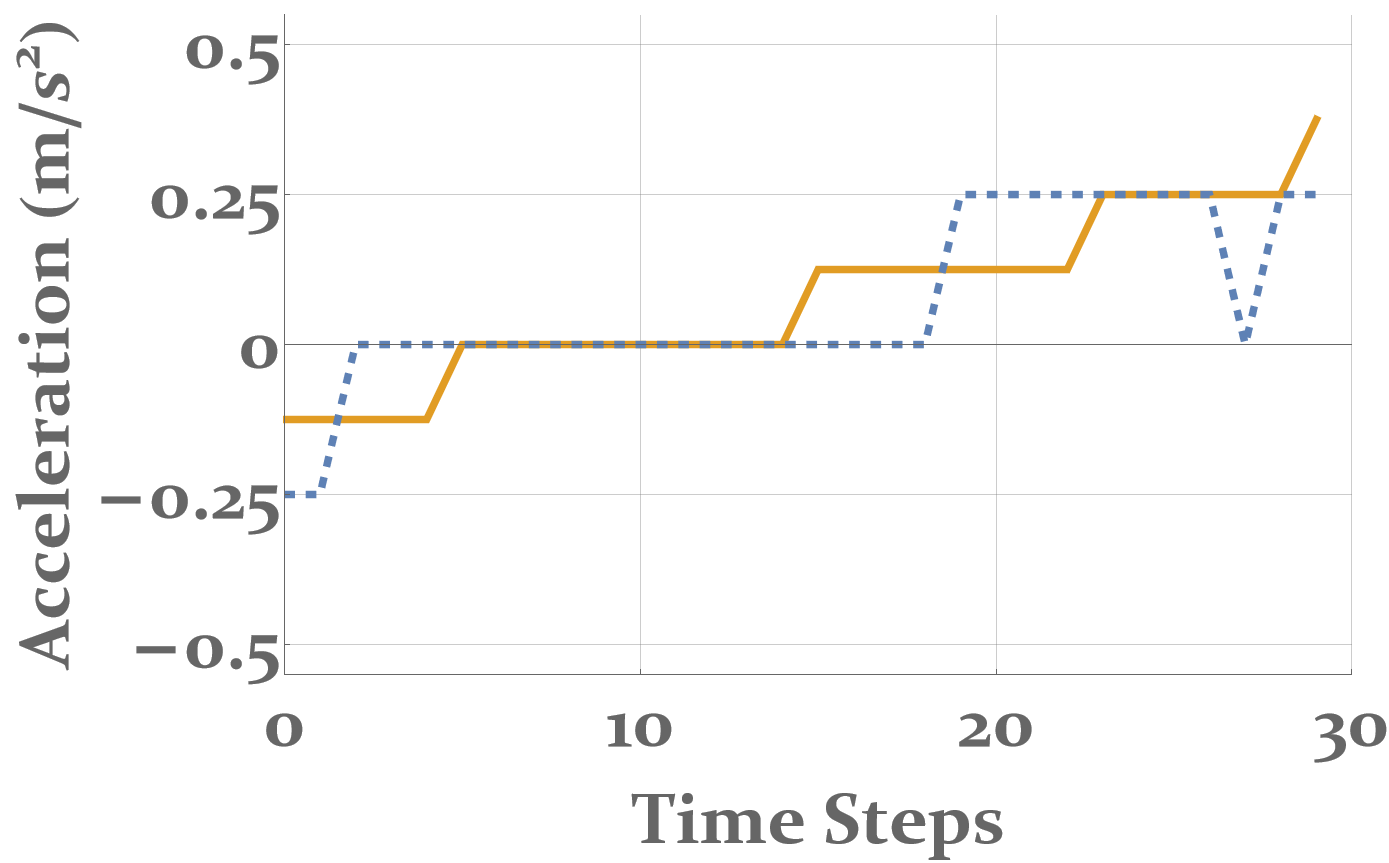}
		\caption*{Iteration 5}
	\end{subfigure}%
	\begin{subfigure}{0.30\linewidth}
		\includegraphics[width=\textwidth]{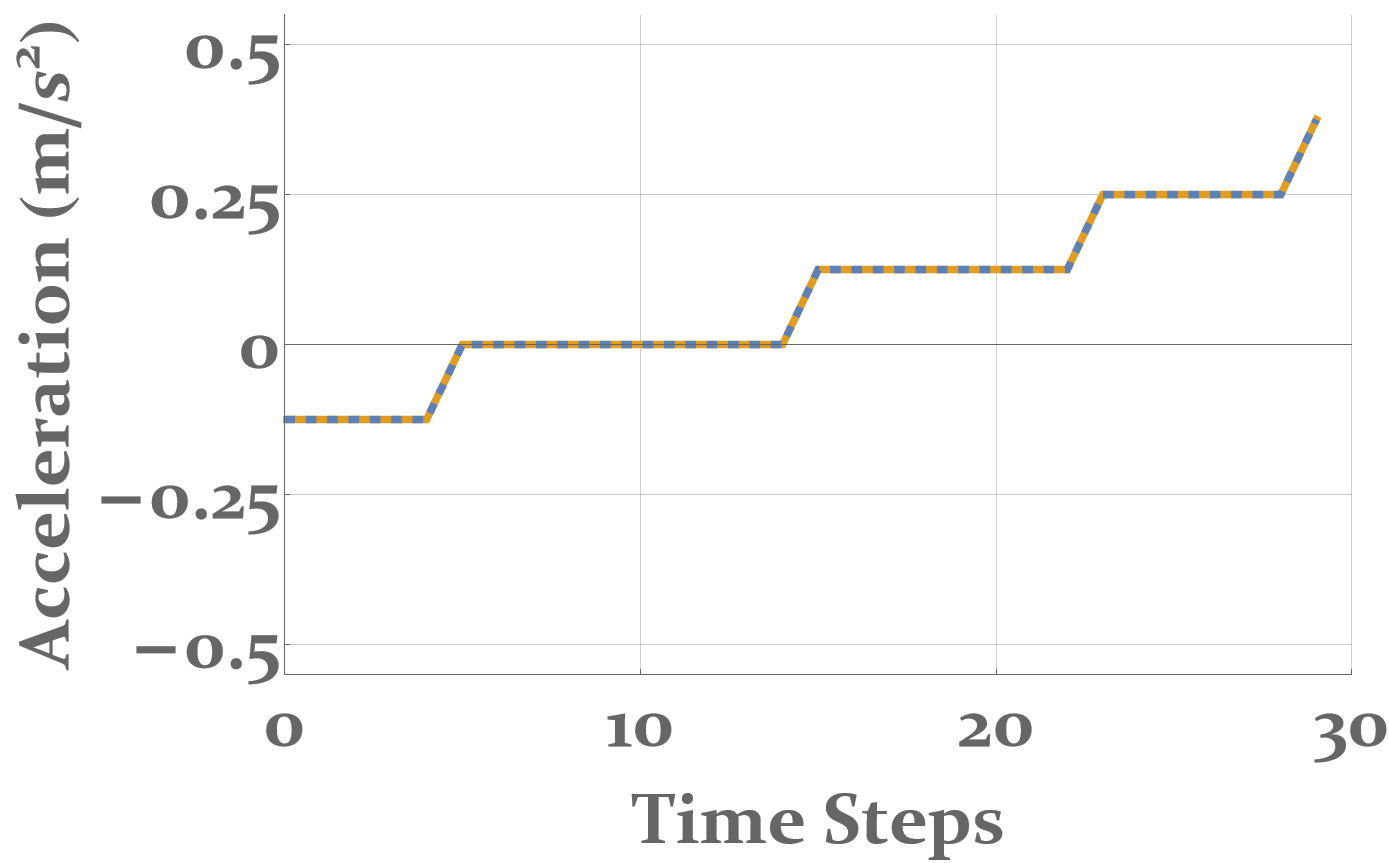}
		\caption*{Iteration 6}
	\end{subfigure}
	\subcaption{}
\end{minipage}
\begin{minipage}{\linewidth}
\centering
	\begin{subfigure}{0.30\linewidth} 
		\includegraphics[width=\textwidth]{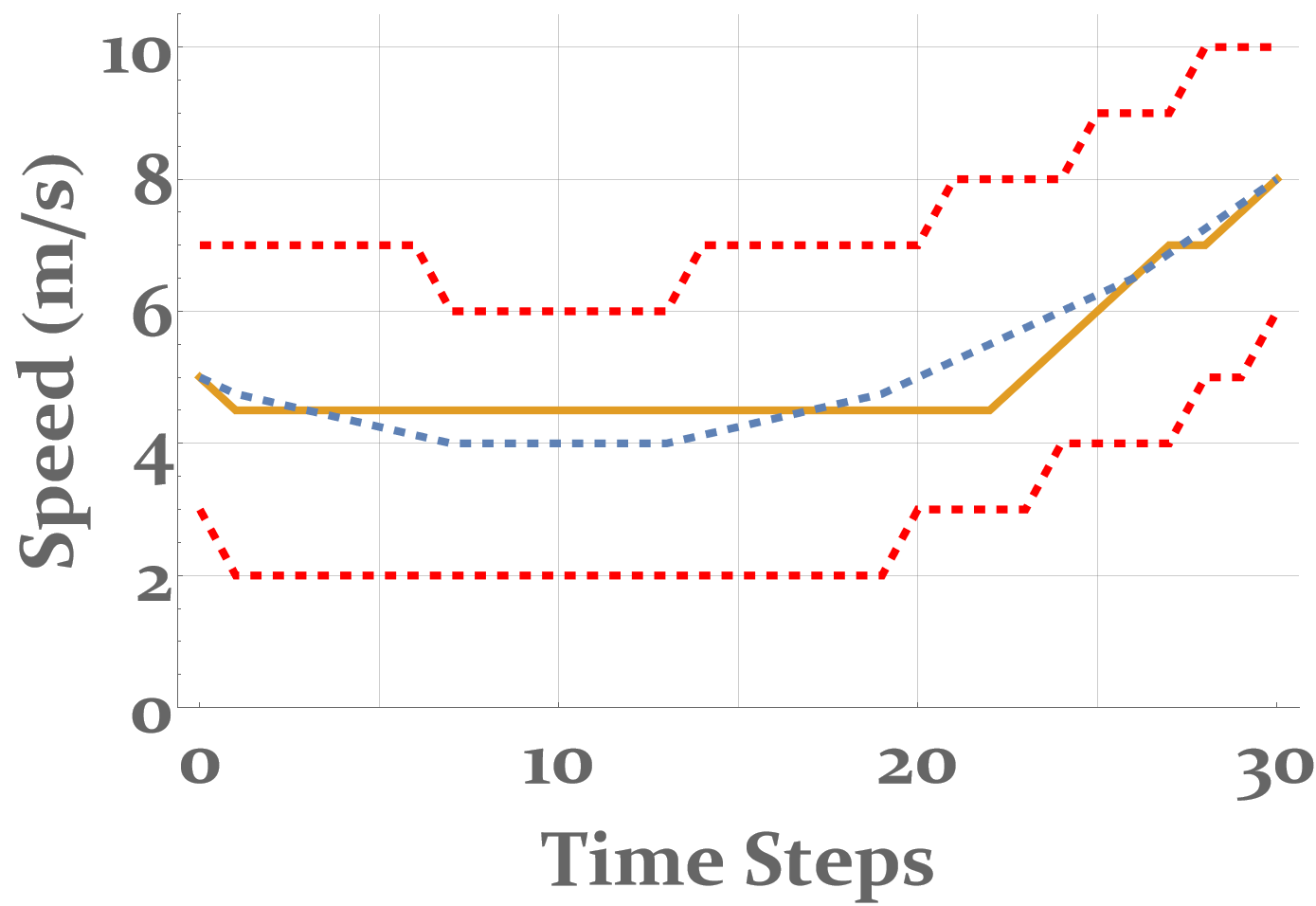}
		\caption*{Iteration 1}
	\end{subfigure}%
	\begin{subfigure}{0.30\linewidth} 
		\includegraphics[width=\textwidth]{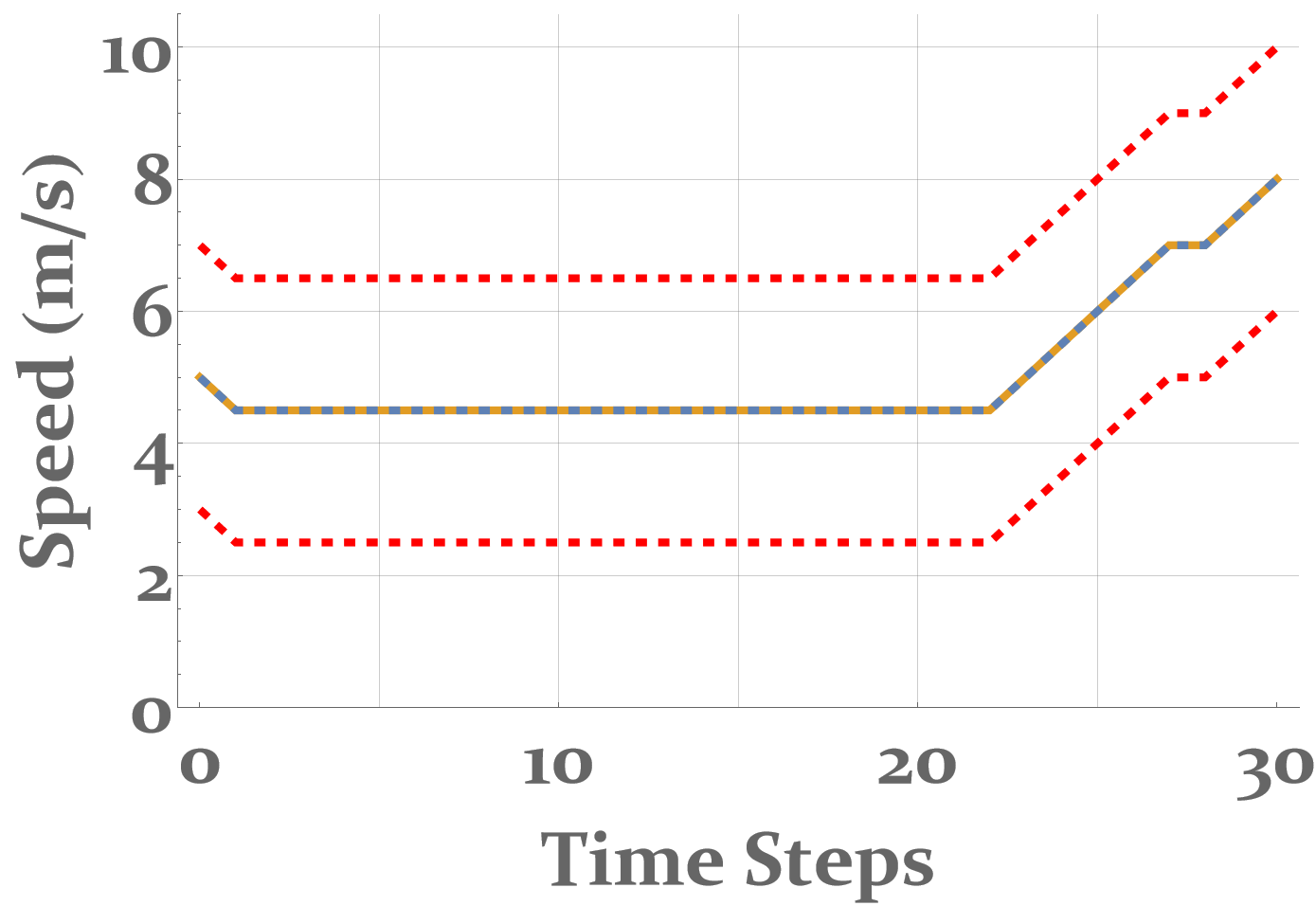}
		\caption*{Iteration 2}
	\end{subfigure}%
	\begin{subfigure}{0.30\linewidth}
		\includegraphics[width=\textwidth]{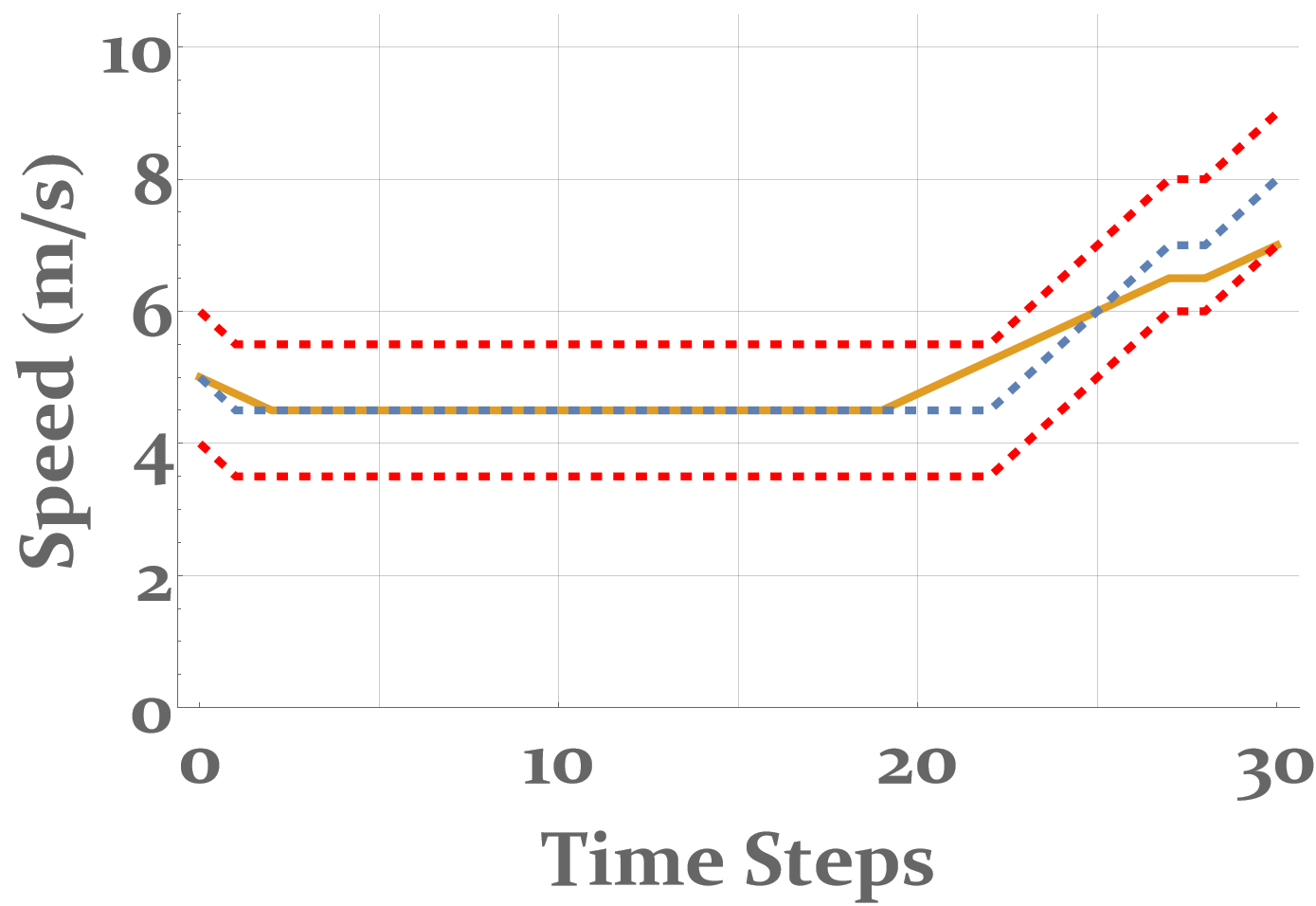}
		\caption*{Iteration 3}
	\end{subfigure}
	\begin{subfigure}{0.30\linewidth}
		\includegraphics[width=\textwidth]{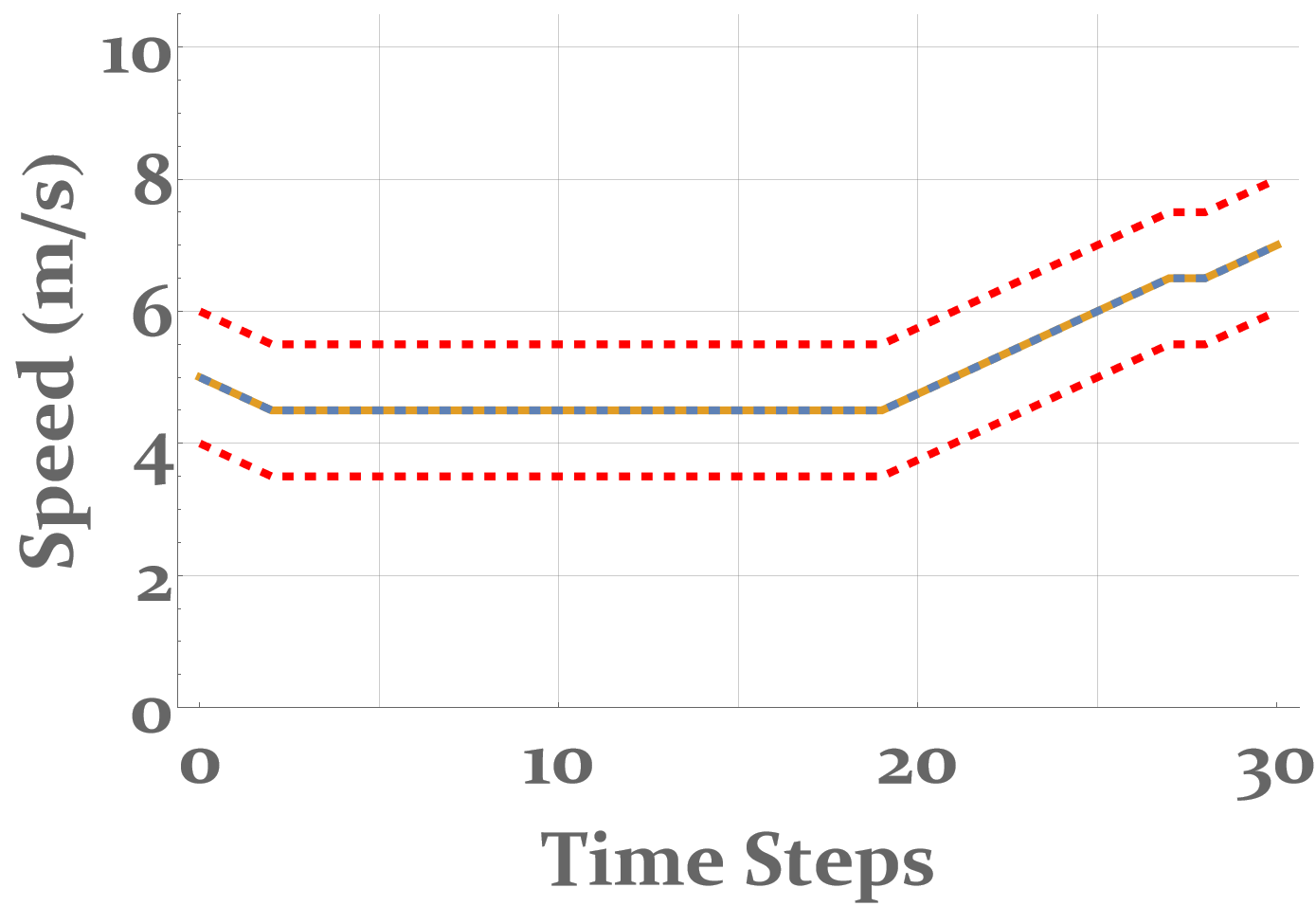}
		\caption*{Iteration 4}
	\end{subfigure}%
	\begin{subfigure}{0.30\linewidth}
		\includegraphics[width=\textwidth]{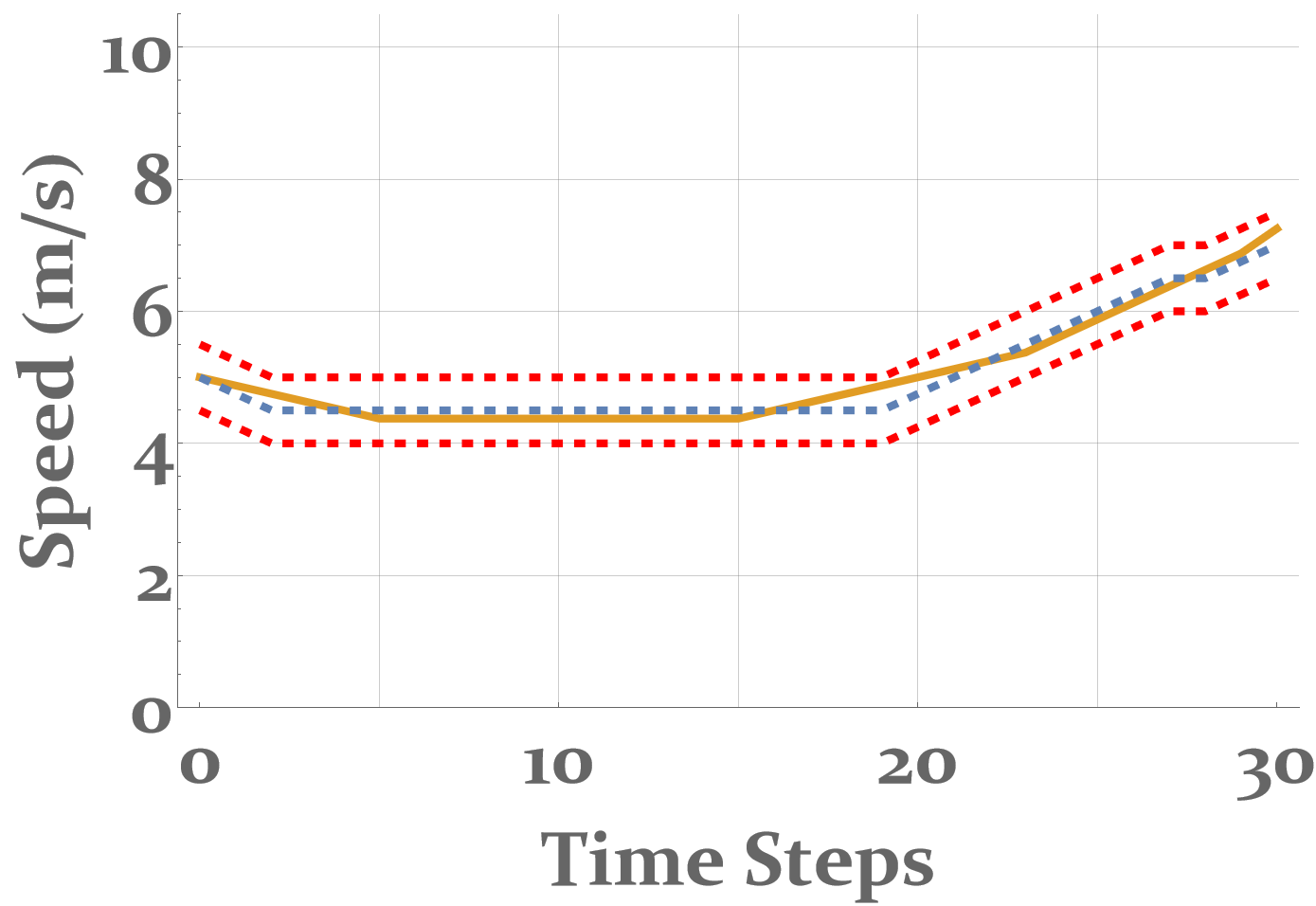}
		\caption*{Iteration 5}
	\end{subfigure}%
	\begin{subfigure}{0.30\linewidth}
		\includegraphics[width=\textwidth]{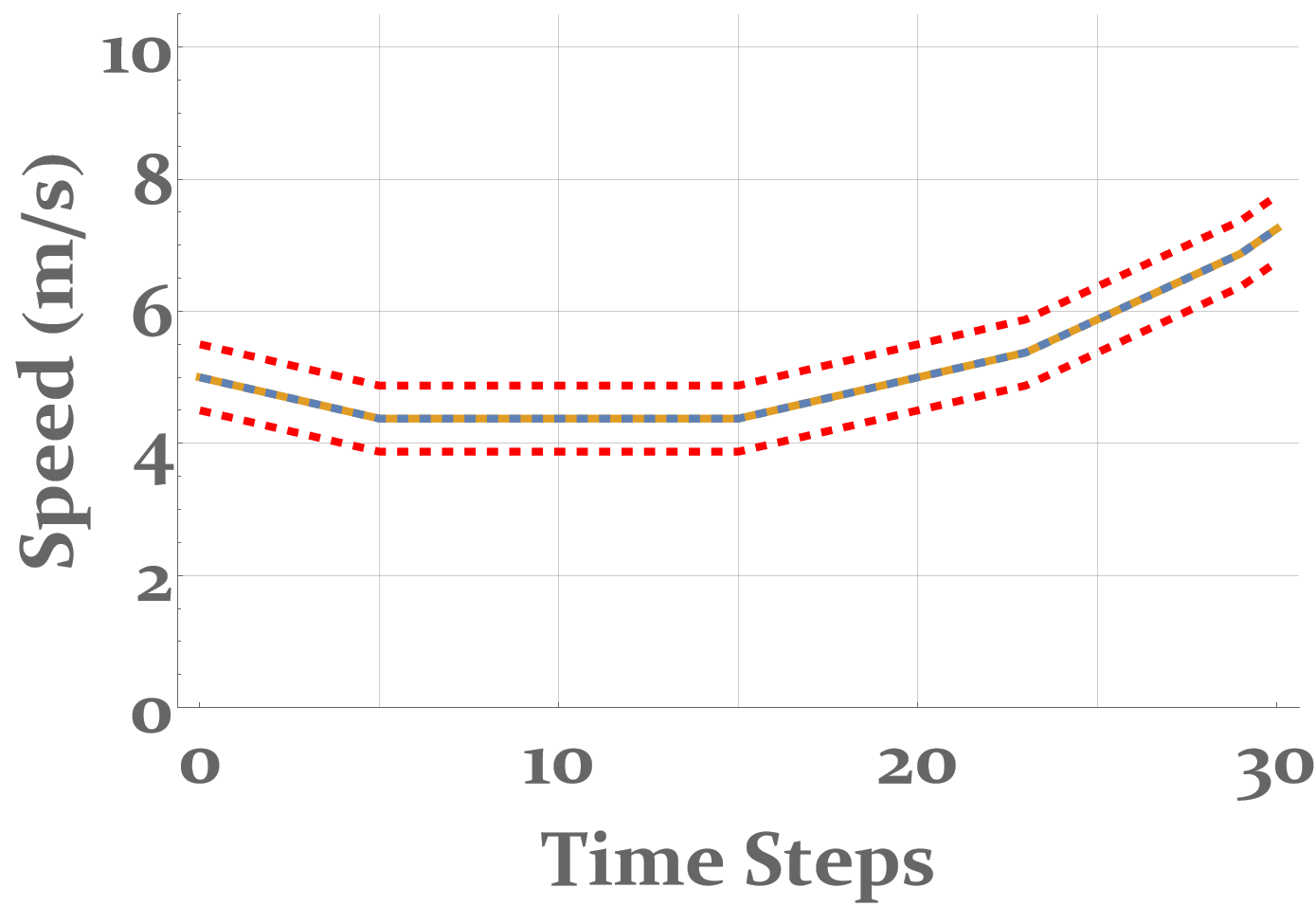}
		\caption*{Iteration 6}
	\end{subfigure}
	\subcaption{}
\end{minipage}
\caption{Received (blue dashed line) and optimal (orange line) (a) acceleration and (b) speed trajectories of DDDP algorithm in each iteration for Scenario 1.} 
\label{fig:DDDP-s-1}
\end{figure}

Fig. \ref{fig:DDDP-s-1} (Scenario 1) and Figs. \ref{fig:DDDP-s-2} and \ref{fig:DDDP-s-3} (for Scenarios 2 and 3, respectively, see Appendix A) display the optimal state (speed) and control (acceleration) trajectories over each iteration of the DDDP algorithm. Note that, position trajectories are not included, as the differences in those trajectories, at each iteration, are small and barely visible. In each iteration, we consider a corridor $\Delta_C^{(l)} = [- \bm{C} \cdot \Delta a^{(l)}, \bm{C} \cdot \Delta a^{(l)}]$ around the respective received state trajectories, which cannot of course extend out of the full state bounds.

In all Figs \ref{fig:DDDP-s-1}, \ref{fig:DDDP-s-2} and \ref{fig:DDDP-s-3}, the dashed blue lines represent, for each iteration, the received trajectory, the solid orange lines represent the optimal trajectories derived, and the red dashed lines reflect the corresponding corridor bounds. For Scenario 1 (Fig. \ref{fig:DDDP-s-1}), it can be observed that, starting with the initial chosen discretization, the first DDDP iteration improves the initial trajectory, leading to a better solution, which is optimal within the considered sub-domain. In the second iteration, no further improvement can be achieved within the new sub-domain bounds, which means that, with the current discretization values, the best achievable solution has been reached. By reducing the discretization at the 3rd iteration, a reduction in the corridor width is observed, but the number of discrete points remains the same. This reduction enables further improvement of the solution, which remains the same in iteration 4; hence the discretisation is further reduced, leading to further improvement in iteration 5; and the procedure continues with iteration 6, where no further improvement is reached and the algorithm terminates, as the discretisation has reached the minimum required value. Similar behaviour is observed for Scenarios 2 and 3, where 10 iterations are required for convergence, and the respective Figs. \ref{fig:DDDP-s-2} and \ref{fig:DDDP-s-3} display a subset thereof. Table \ref{tb:DDDP_cost_da} contains the values of the objective function for each iteration of the DDDP algorithm for all scenarios, assuming $\Delta a = 0.5$ and $\bm{C} = [20,4]$, where the choice of these values will be explained in the following.

The DDDP has been evaluated for different corridor $C_v$ and initial discretization $\Delta a$ values, for the three investigated scenarios (see Table \ref{tb:dddpC-Da} in Appendix A). It can be seen that, for all presented initial discretization values, the obtained optimal solution at convergence reaches the same globally optimal solution, if the corridor width is sufficiently large, specifically for $C_v \geq 4$ (or lower, in some cases). This is because, for very small values of $C_v$, DDDP does not have sufficient flexibility to improve the solution at each iteration due to the extremely limited state space, something that may lead to a local optimum of the overall problem at convergence. At the same time, a reduction of the number of DDDP iterations is obtained for larger corridors. This behaviour is expected, as the bigger corridors, i.e., bigger admissible sub-regions, lead to potentially better solutions at each iteration, and hence to fewer iterations to converge to the optimal solution. On the other hand, despite the decrease on the number of iterations, the overall computation time may increase due to the higher computation time required at each iteration, which in turn, is due to more feasible discrete state points included in each iteration's sub-region.  Based on these observations, the selection of the values $\bm{C} = [20,4]$ and $\Delta a = 0.5$ seems to be a reasonable choice.

\begin{table}[tb]
\begin{center}
\caption{Optimal cost evolution and descretization change in each iteration of DDDP algorithm.}\label{tb:DDDP_cost_da}
\begin{tabular}{cccc||ccc||ccc}
\hline
& \multicolumn{3}{c||}{Scenario 1} & \multicolumn{3}{c||}{Scenario 2}  & \multicolumn{3}{c}{Scenario 3} \\\hline
Iter. & $\Delta a = \Delta v$ & $\Delta x$ & Cost & $\Delta a = \Delta v$ & $\Delta x$ & Cost & $\Delta a = \Delta v$ & $\Delta x$ & Cost \\\hline
1  & 0.5	&		0.25   & 1.357 		& 0.5	& 0.25   & 4.311		& 0.5	& 0.25   & 6.6218 \\
2  & 0.5 	&		0.25   & 1.357 		& 0.5 	& 0.25   & 4.231		& 0.5 	& 0.25   & 6.6218 \\
3  & 0.25 	&		0.125  & 1.223 		& 0.5 	& 0.25   & 4.231		& 0.25 	& 0.125   & 6.4532 \\
4  & 0.25 	&		0.125  & 1.223 		& 0.25 	& 0.125  & 4.012		& 0.25 	& 0.125  & 6.4346 \\
5  & 0.125	&		0.0625 & 1.175		& 0.25	& 0.125  & 3.988		& 0.25	& 0.125  & 6.4208 \\
6  & 0.125	&		0.0625 & 1.175		& 0.25	& 0.125  & 3.958		& 0.25	& 0.125  & 6.4091 \\
7  &		&			   &			& 0.25	& 0.125  & 3.940		& 0.25	& 0.125  & 6.4091 \\
8  &		&			   &			& 0.25	& 0.125  & 3.940		& 0.125	& 0.0625  & 6.3583 \\
9  &		&			   &			& 0.125	& 0.0625 & 3.906		& 0.125	& 0.0625 & 6.3582 \\
10 &		&			   &			& 0.125	& 0.0625 & 3.906		& 0.125	& 0.0625 & 6.3582 \\
\hline
\end{tabular}%
\end{center}
\end{table}

The accuracy of the DDDP algorithm compared to the full one-shot SDP solution may be assessed by contrasting the optimal cost and the optimal states and control trajectories, obtained in the last respective DDDP iterations of the three scenarios, to the corresponding optimal trajectories derived from the full-range SDP with a discretisation of $\Delta a=0.125$. As a matter of fact, DDDP manages, in all cases, to converge to the exact same optimal solution as the full-range SDP. The obtained optimal costs of both approaches, for the three scenarios, are 1.175, 3.906 and 6.3582, respectively, while the computation time differences are remarkable, as the one-shot SDP needs approximately 613 s to obtain the solution (same for all scenarios), while DDDP needs only 0.69, 1.24 and 1.12 s for the respective three scenarios. More importantly, this big reduction in computation time enables the DDDP algorithm to be executable in real time, even in an MPC mode, on the vehicle side, similarly to the deterministic GLOSA.


\subsection{DDP algorithm results} \label{sec:results-ddp}

\begin{figure}[t]
\centering
\begin{minipage}[b]{\textwidth}
	\begin{subfigure}{0.33\textwidth} 
		\includegraphics[width=\textwidth]{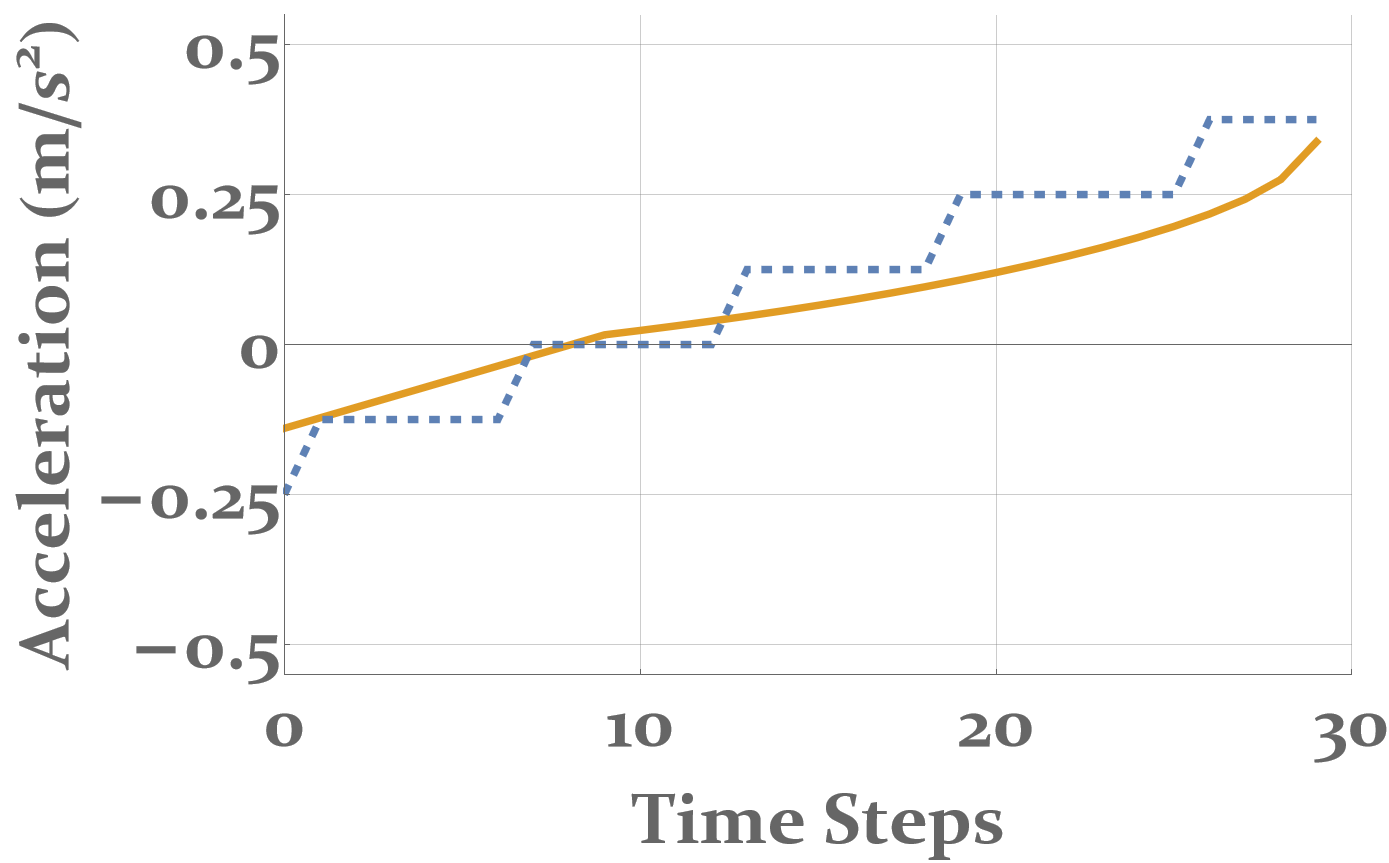}
		\caption*{Iteration 1}
	\end{subfigure}%
	\begin{subfigure}{0.33\textwidth} 
		\includegraphics[width=\textwidth]{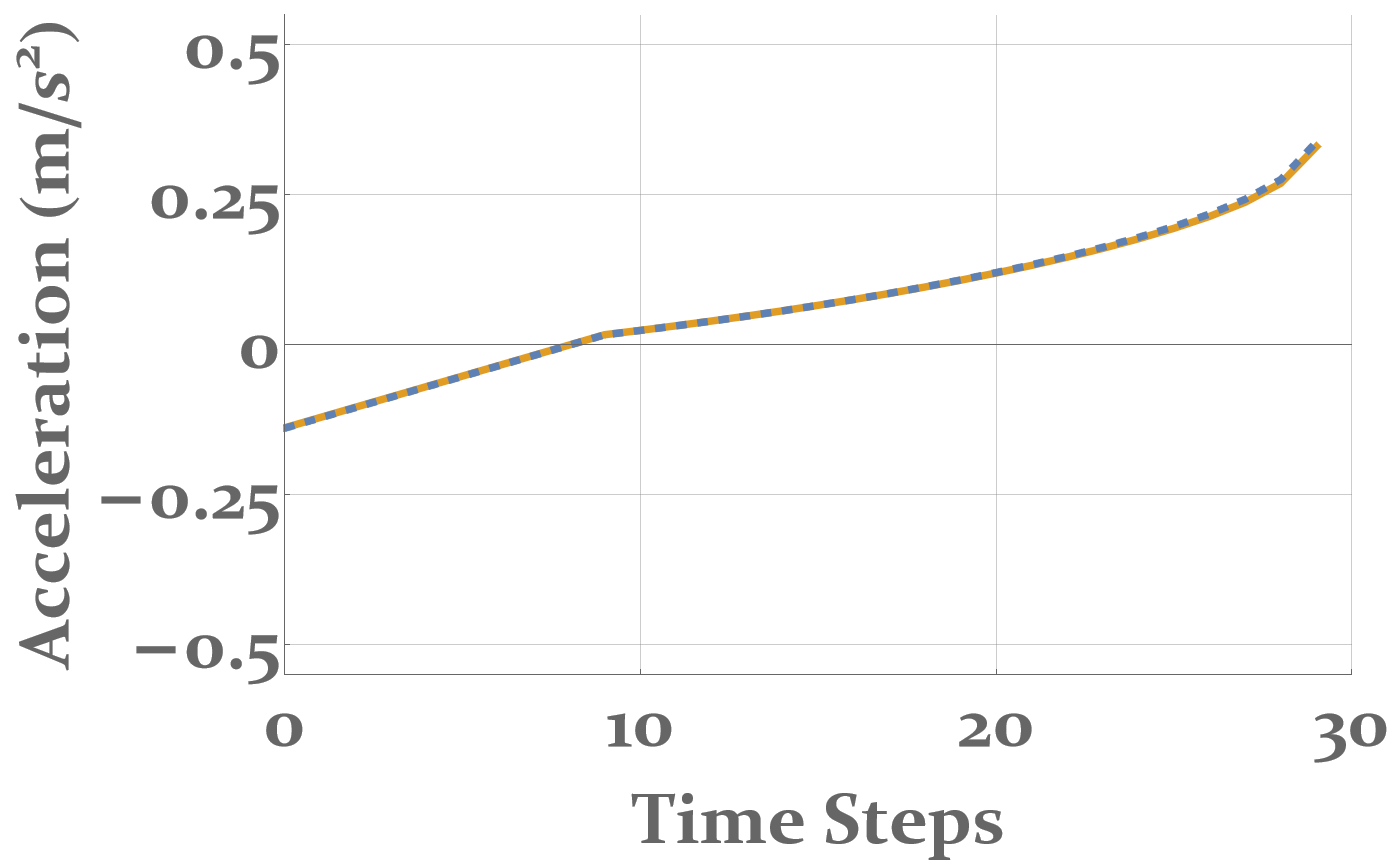}
		\caption*{Iteration 2}
	\end{subfigure}%
	\begin{subfigure}{0.33\textwidth}
		\includegraphics[width=\textwidth]{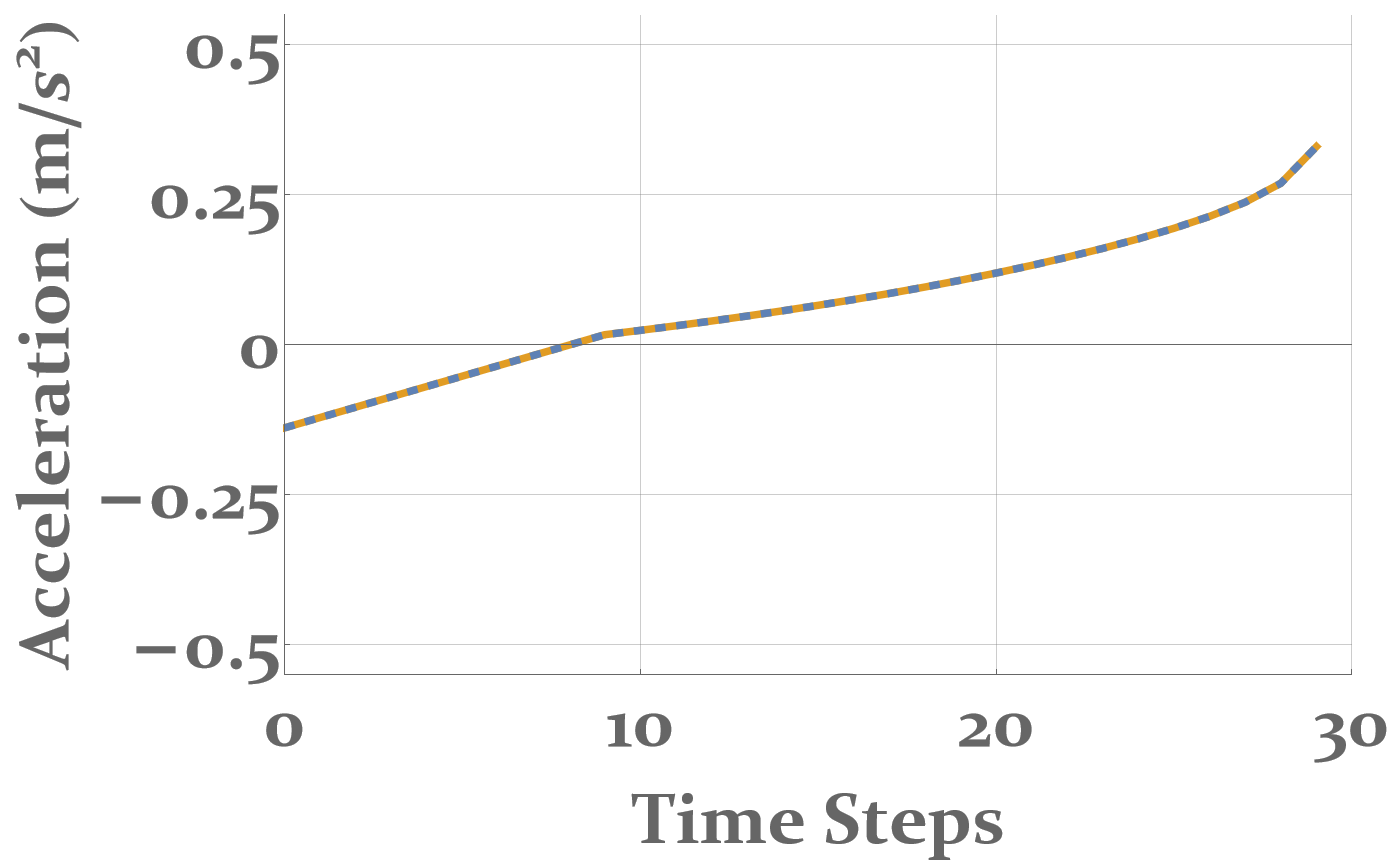}
		\caption*{Iteration 3}
	\end{subfigure}
\subcaption{ }
\end{minipage}

\begin{minipage}[b]{\textwidth}
	\begin{subfigure}{0.33\textwidth}
		\includegraphics[width=\textwidth]{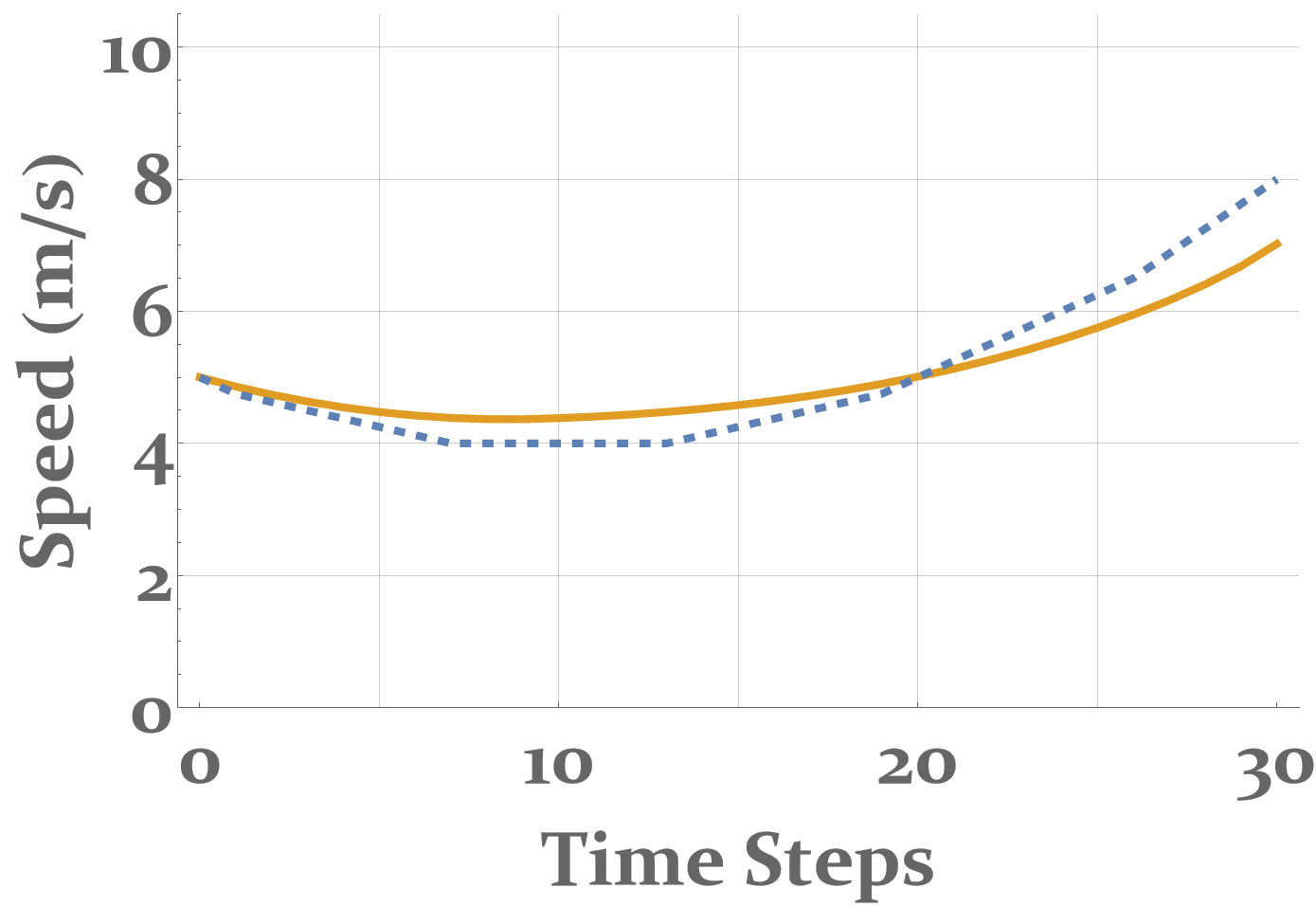}
		\caption*{Iteration 1}
	\end{subfigure}%
	\begin{subfigure}{0.33\textwidth}
		\includegraphics[width=\textwidth]{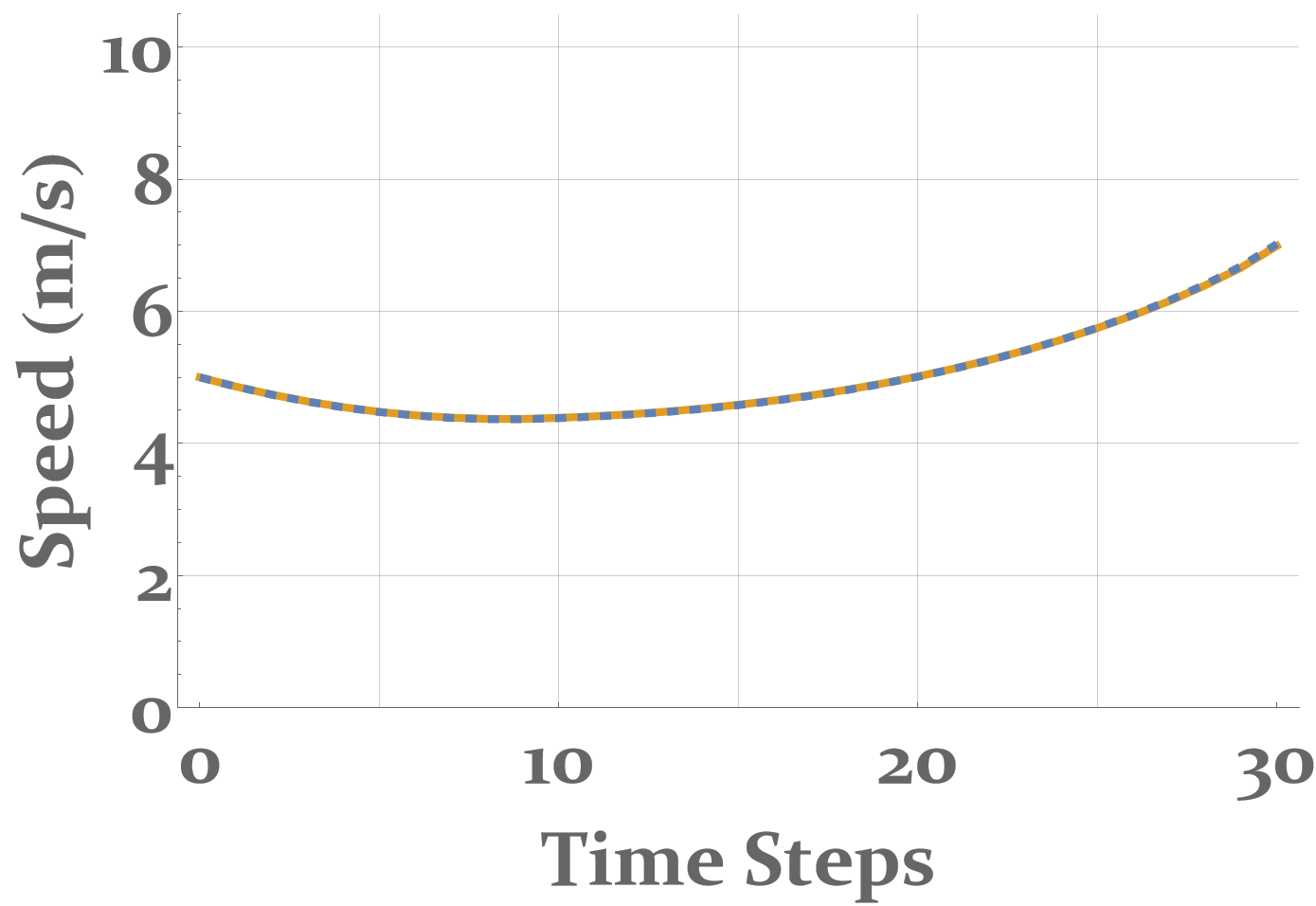}
		\caption*{Iteration 2}
	\end{subfigure}%
	\begin{subfigure}{0.33\textwidth}
		\includegraphics[width=\textwidth]{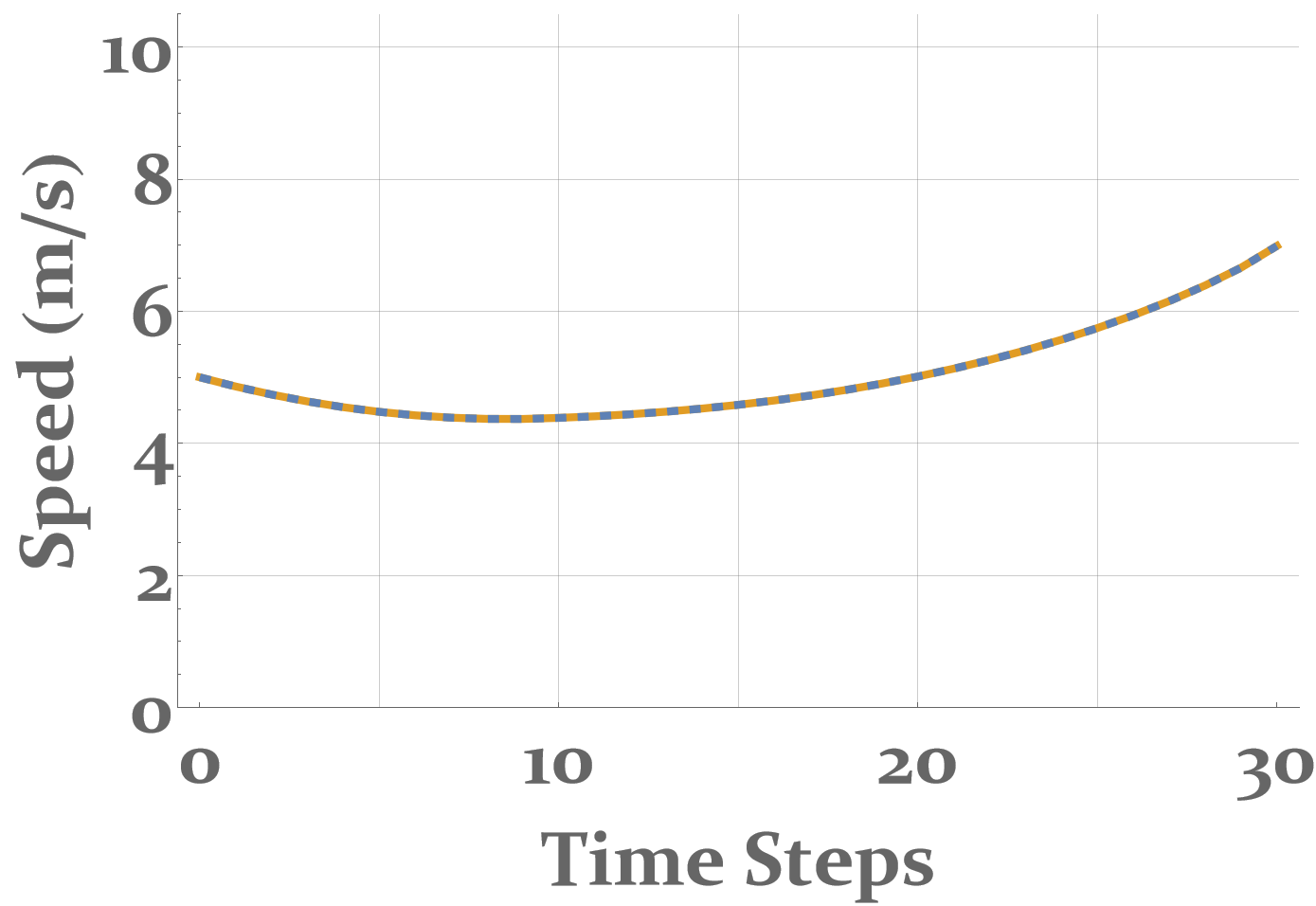}
		\caption*{Iteration 3}
	\end{subfigure}
\subcaption{ }
\end{minipage}
\caption{Received (blue dashed line) and computed (orange line) trajectories of (a) acceleration and (b) speed at each iteration of DDP algorithm for Scenario 1.}
\label{fig:ddp-s-1}
\end{figure}

Fig. \ref{fig:ddp-s-1} (Scenario 1) and Figs. \ref{fig:ddp-s-2} and \ref{fig:ddp-s-3}  (for Scenarios 2 and 3, respectively, see Appendix A) illustrate the evolution of the state (speed) and control (acceleration) trajectories produced by the DDP algorithm at each iteration. In all figures, blue dashed lines represent the initial trajectories, and solid orange lines the derived trajectories of each iteration. From the figures, it can be seen that the DDP algorithm, starting from the pessimistic deterministic GLOSA solution, converges to the overall optimal solution. It can also be observed that, at each iteration, starting from the initial (received) trajectory, the algorithm results in an improved one, which leads to an improved solution. The DDP algorithm is seen to converge very fast to the optimal solution, with the improvement being more profound at the first iterations, as the change between the nominal and the produced trajectories are more noticeable. On the other hand, in the last iterations, the improvement is smaller, as the algorithm is already very close to the optimal solution. As mentioned earlier, all intermediate DDP solutions at each iteration are continuous due to the nature of DDP, which is operating on a continuous-value basis for all involved variables.

\begin{figure}[h!t]
\centering
\begin{minipage}[b]{\textwidth}
	\begin{subfigure}{0.33\textwidth} 
		\includegraphics[width=\textwidth]{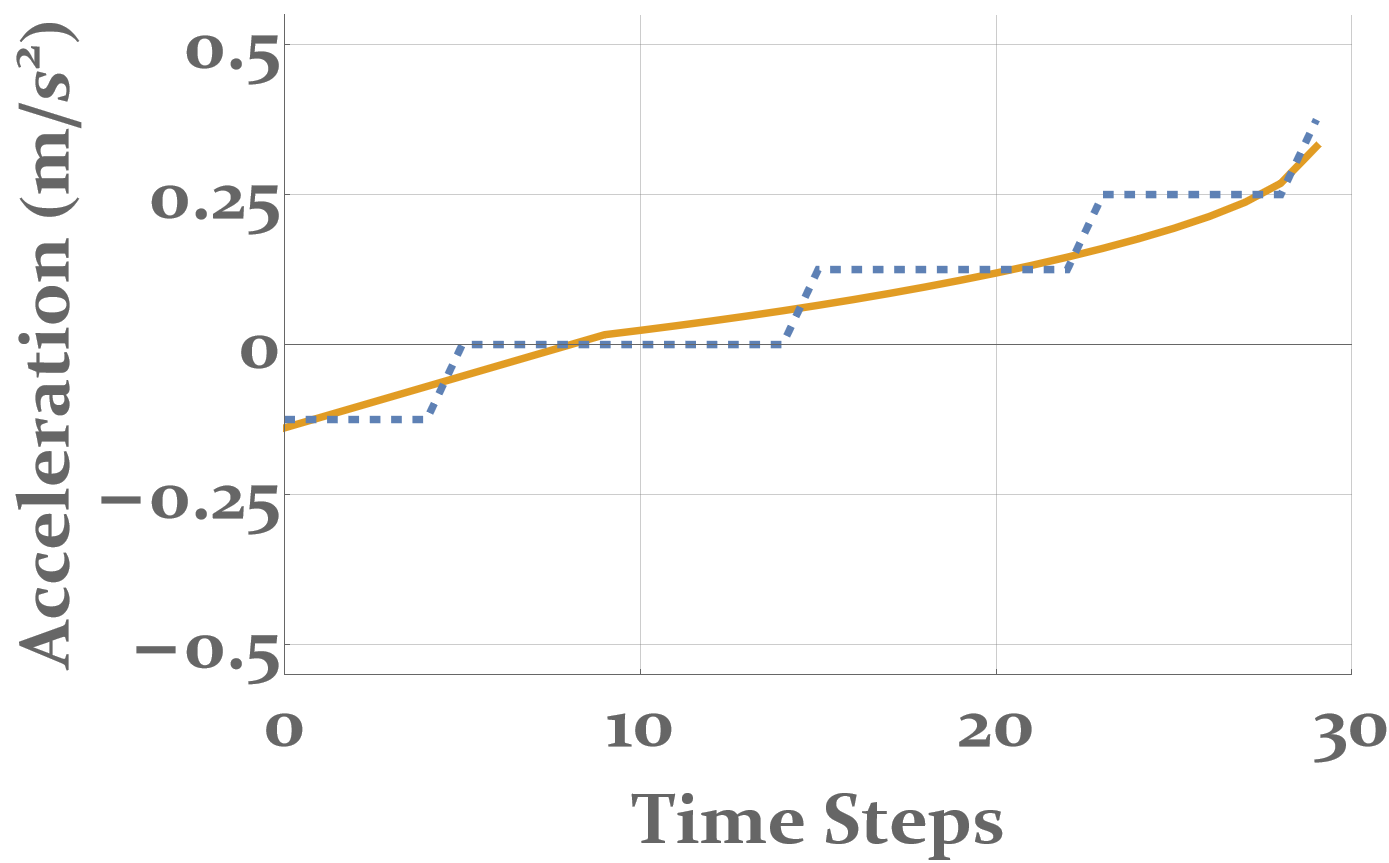}
		\caption*{Iteration 3}
	\end{subfigure}%
	\begin{subfigure}{0.33\textwidth} 
		\includegraphics[width=\textwidth]{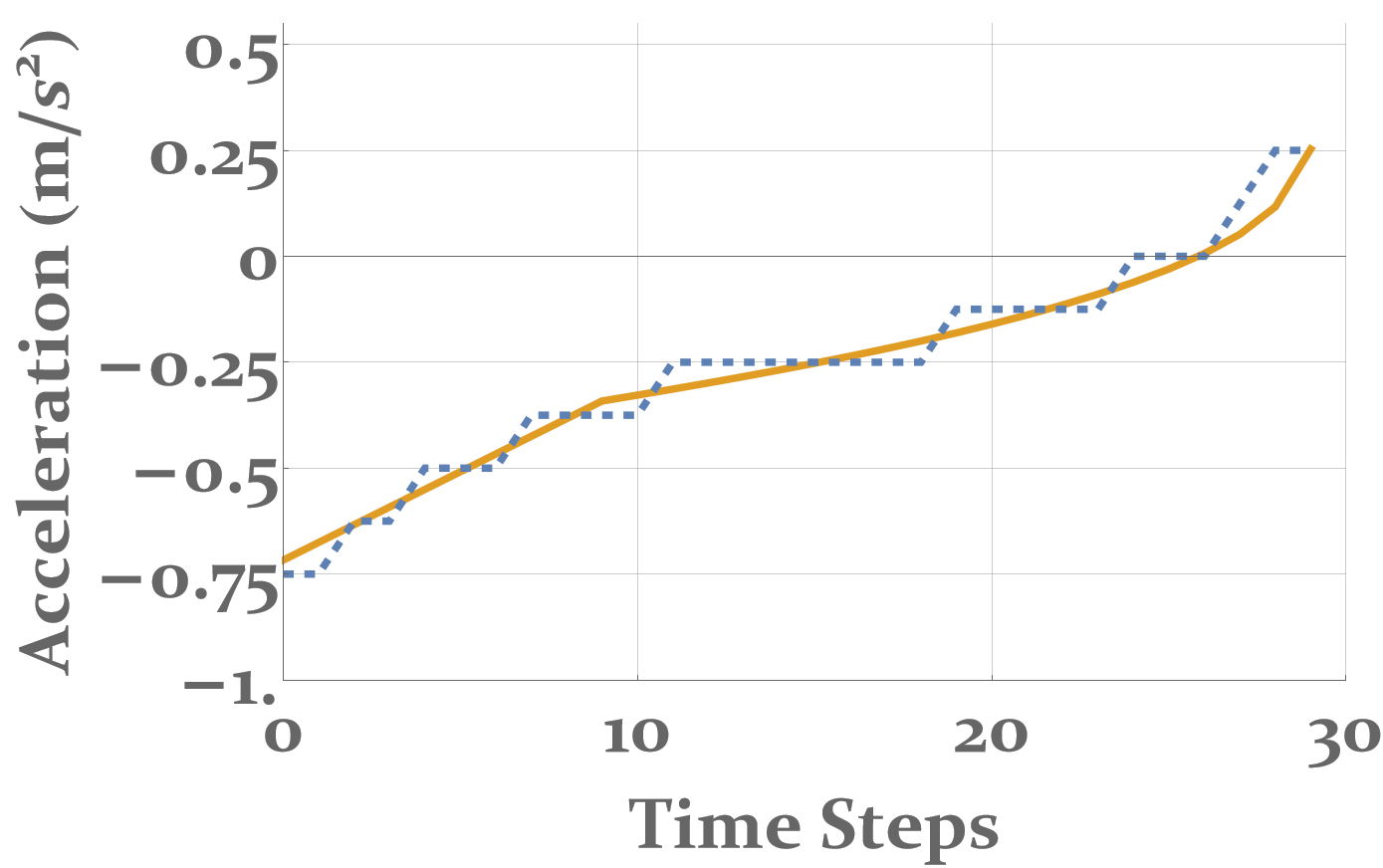}
		\caption*{Iteration 4}
	\end{subfigure}%
	\begin{subfigure}{0.33\textwidth}
		\includegraphics[width=\textwidth]{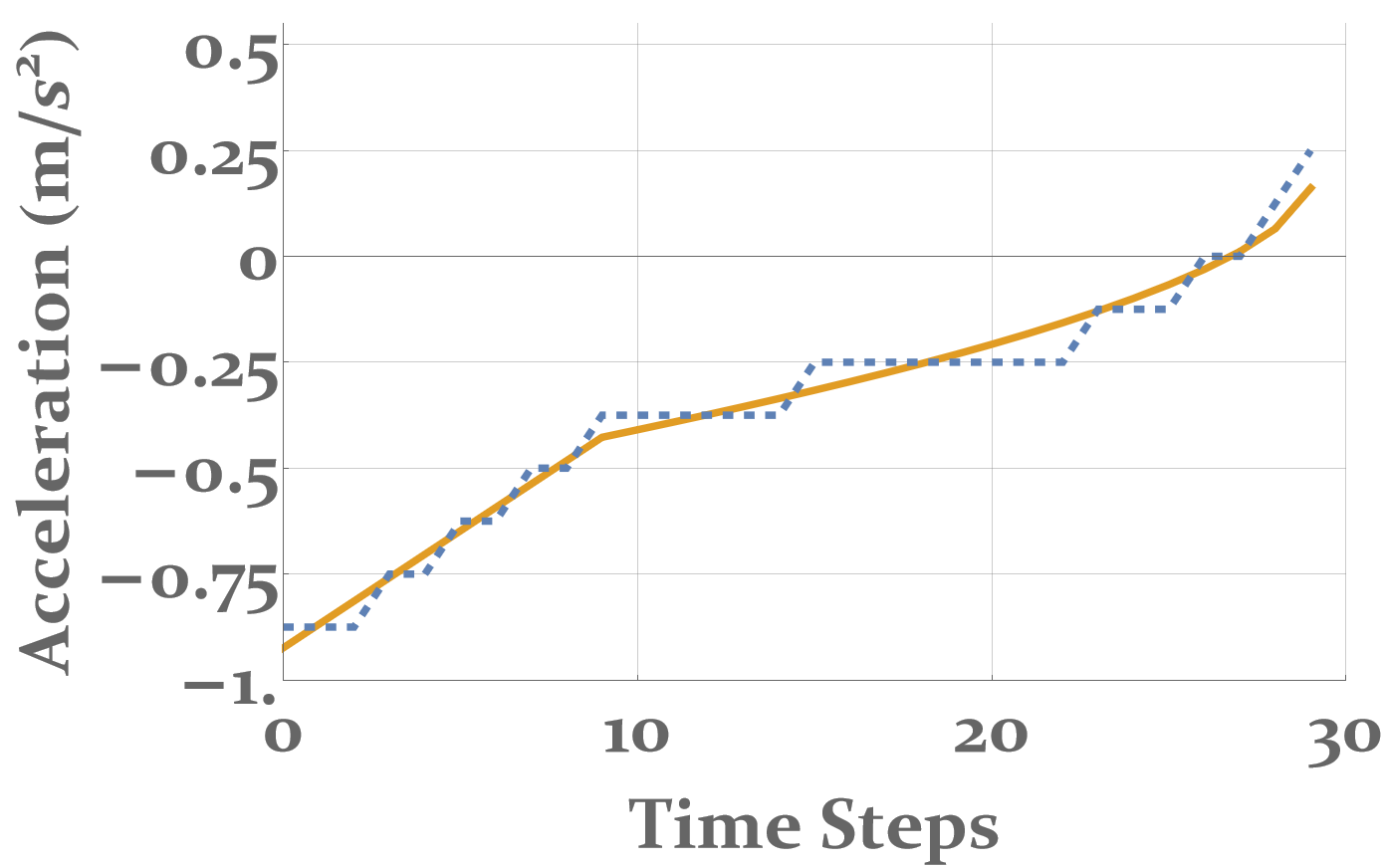}
		\caption*{Iteration 4}
	\end{subfigure}
	\begin{subfigure}{0.33\textwidth}
		\includegraphics[width=\textwidth]{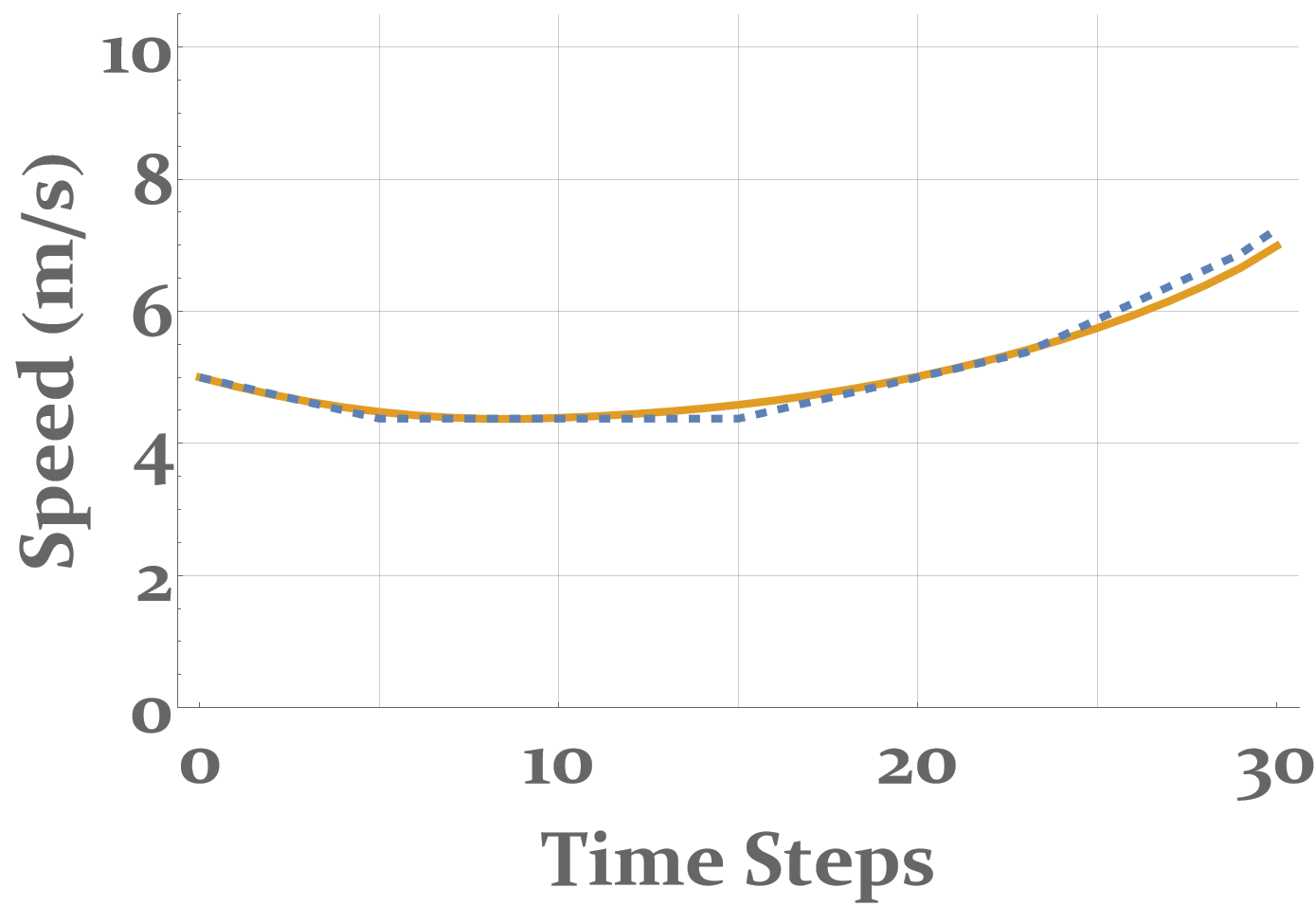}
		\caption*{Iteration 3}
		\caption{ }
	\end{subfigure}%
	\begin{subfigure}{0.33\textwidth}
		\includegraphics[width=\textwidth]{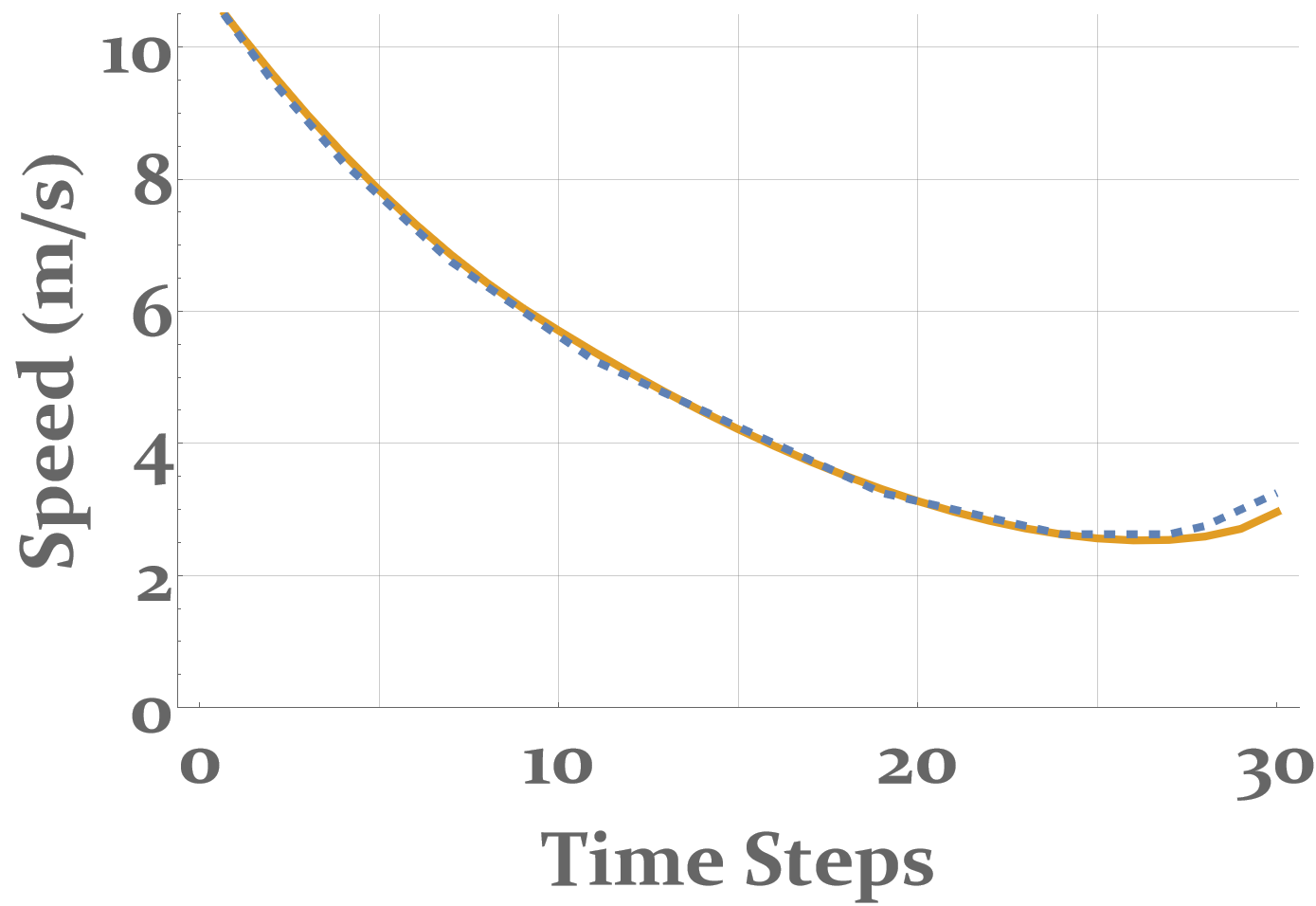}
		\caption*{Iteration 4}
		\caption{ }
	\end{subfigure}%
	\begin{subfigure}{0.33\textwidth}
		\includegraphics[width=\textwidth]{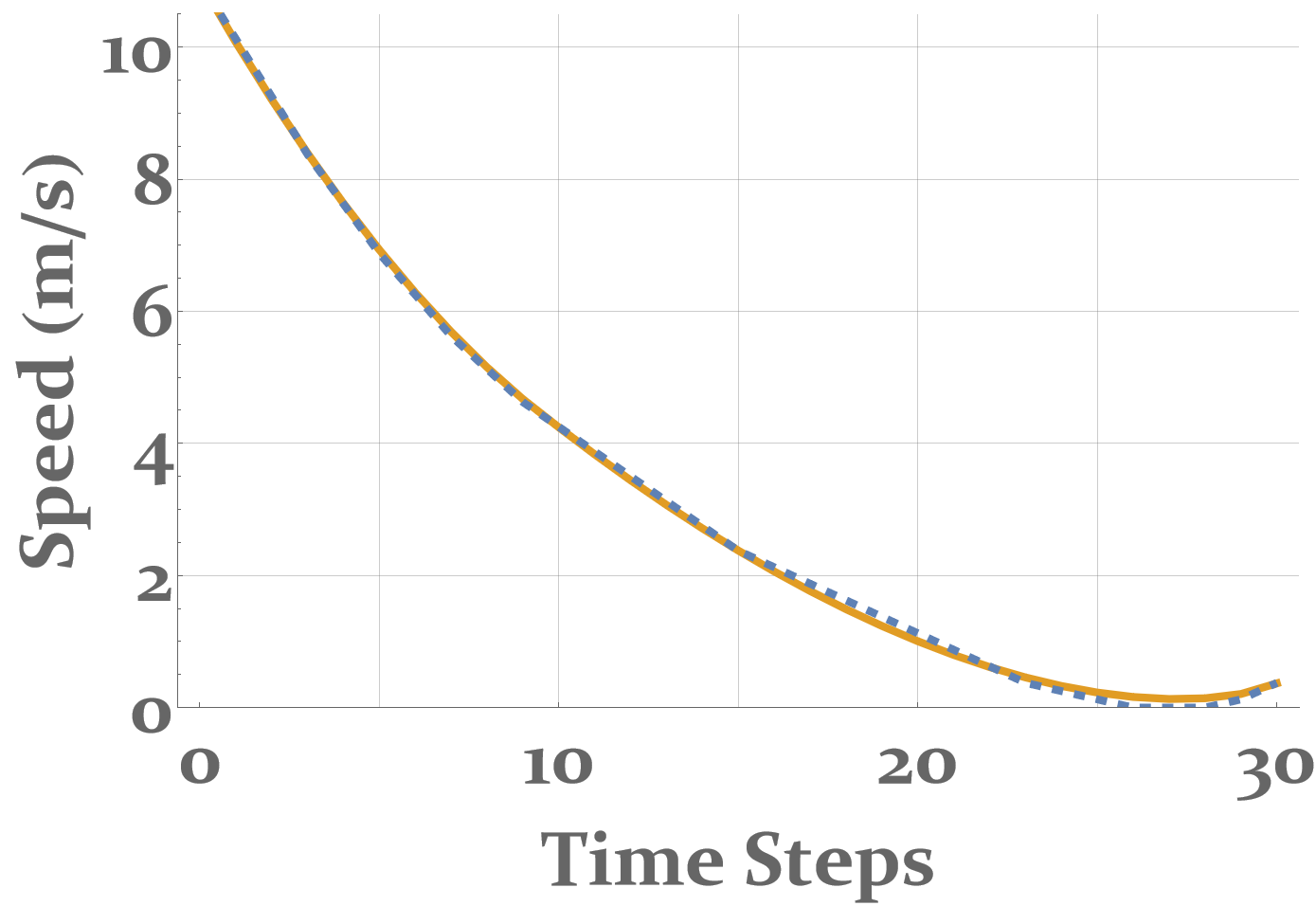}
		\caption*{Iteration 4}
		\caption{ }
	\end{subfigure}
\caption{Optimal acceleration and speed trajectories (orange lines) of DDP, compared with the corresponding optimal trajectories of the one-shot SDP (blue dashed lines) of (a) Scenario 1, (b) Scenario 2 and (c) Scenario 3.}\label{fig:ddpvssdp}
\end{minipage}
\end{figure}

Fig. \ref{fig:ddpvssdp} demonstrates the accuracy of the DDP solutions compared to the one-shot discrete SDP solutions for the three scenarios. In this figure, the optimal state (speed) and control (acceleration) trajectories obtained from DDP algorithm at convergence are contrasted to those obtained from the full-range SDP algorithm with discretization $\Delta a=0.125$ m/s$^2$. The two solutions, in all scenarios, are quite close to each other, with a slight advantage for the DDP solution, which is not bound to a fine, but still visible value-discretization of the state and control, as it is the case for the one-shot discrete SDP algorithm.

\subsection{Summarized results} \label{sec:results-sum}

Table \ref{tbl: comparison} contains the summarized information of the number of iterations until convergence, CPU time and the optimal obtained cost of the one-shot SDP, the iterative DDDP and the iterative, but non-discrete DDP algorithms for Scenarios 1, 2 and 3. Based on the results reported in the table, it can be verified that DDP's optimal solution is slightly better compared to the other two approaches. However, the most striking difference of DDP is the computational time needed to solve the problem, which is extremely small. 

\begin{table}[!htb]
\caption{Comparison of iterations, CPU times until convergence and obtained optimal cost of SDP, DDDP and DDP}
\label{tbl: comparison}
\begin{center}
\begin{tabular}{lccc||ccc||ccc}
\hline
\multicolumn{1}{c}{ } & \multicolumn{3}{c}{Scenario 1} & \multicolumn{3}{c}{Scenario 2} & \multicolumn{3}{c}{Scenario 3} \\
\hline
 		& Iter. & \begin{tabular}{@{}c@{}} CPU \\ Times (s)\end{tabular} & Cost & Iter. & \begin{tabular}{@{}c@{}} CPU \\ Times (s)\end{tabular} & Cost & Iter. & \begin{tabular}{@{}c@{}} CPU \\ Times (s)\end{tabular} & Cost\\
\hline
\begin{tabular}{@{}l@{}} 
One-shot \\ 
SDP\end{tabular} 	& -	& 631	 & 1.175 &  - & 631		& 3.906 &  - & 631 & 6.358 \\
\hline
DDDP 				& 6 & 0.692  & 1.175 & 10 & 1.246	& 3.906 &  10 & 1.128 & 6.358 \\
\hline
DDP 				& 3	& 0.0005 & 1.162 &	4 & 0.0007	& 3.892 &  4 & 0.0007 & 6.353 \\
\hline
\end{tabular}
\end{center}
\end{table}

\begin{figure}[h!t]
\centering
	\begin{subfigure}{0.45\textwidth} 
		\includegraphics[width=\textwidth]{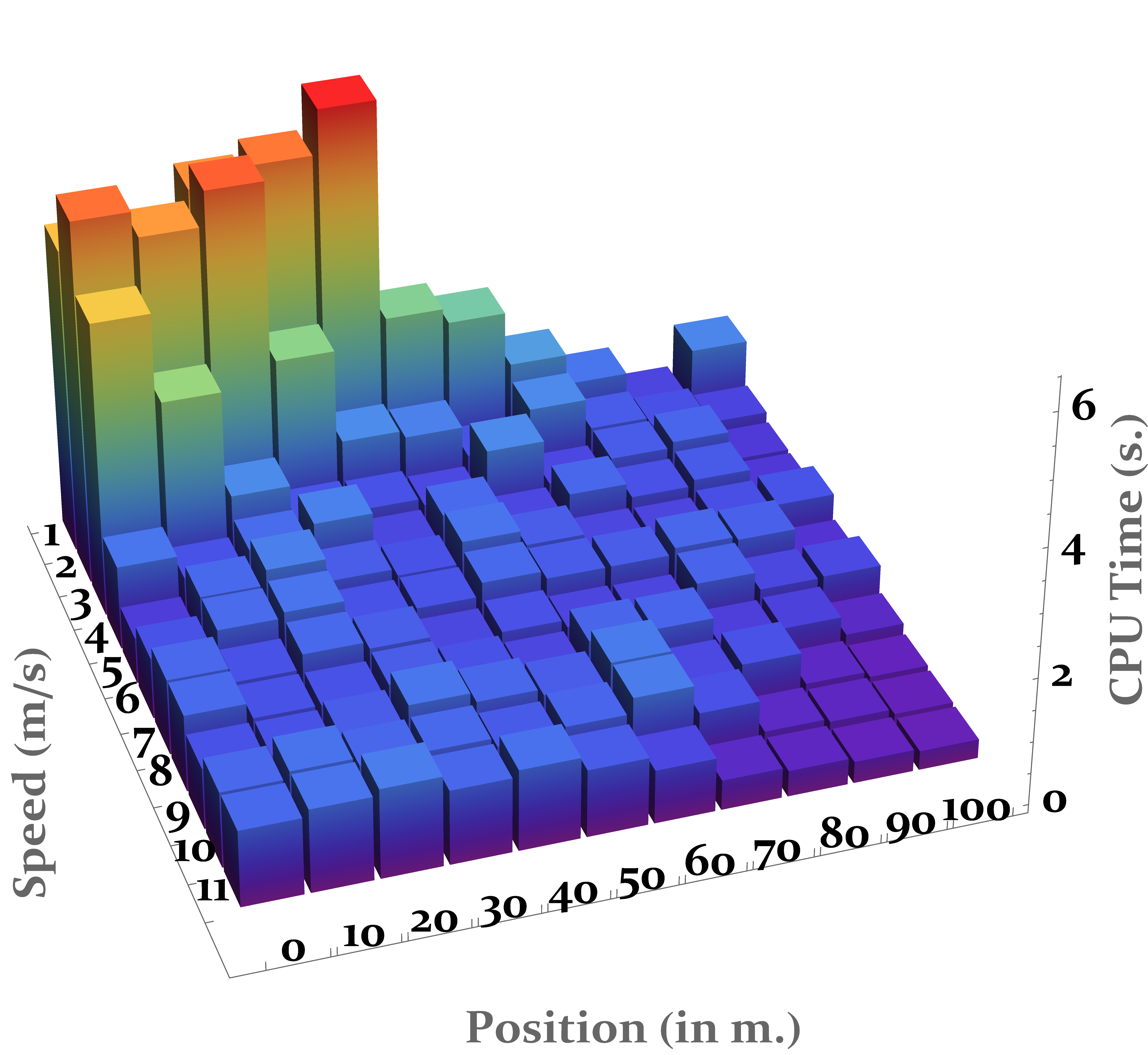}
		\caption{}
		\label{fig:cpuTimes-1}
	\end{subfigure}
	\begin{subfigure}{0.45\textwidth} 
		\includegraphics[width=\textwidth]{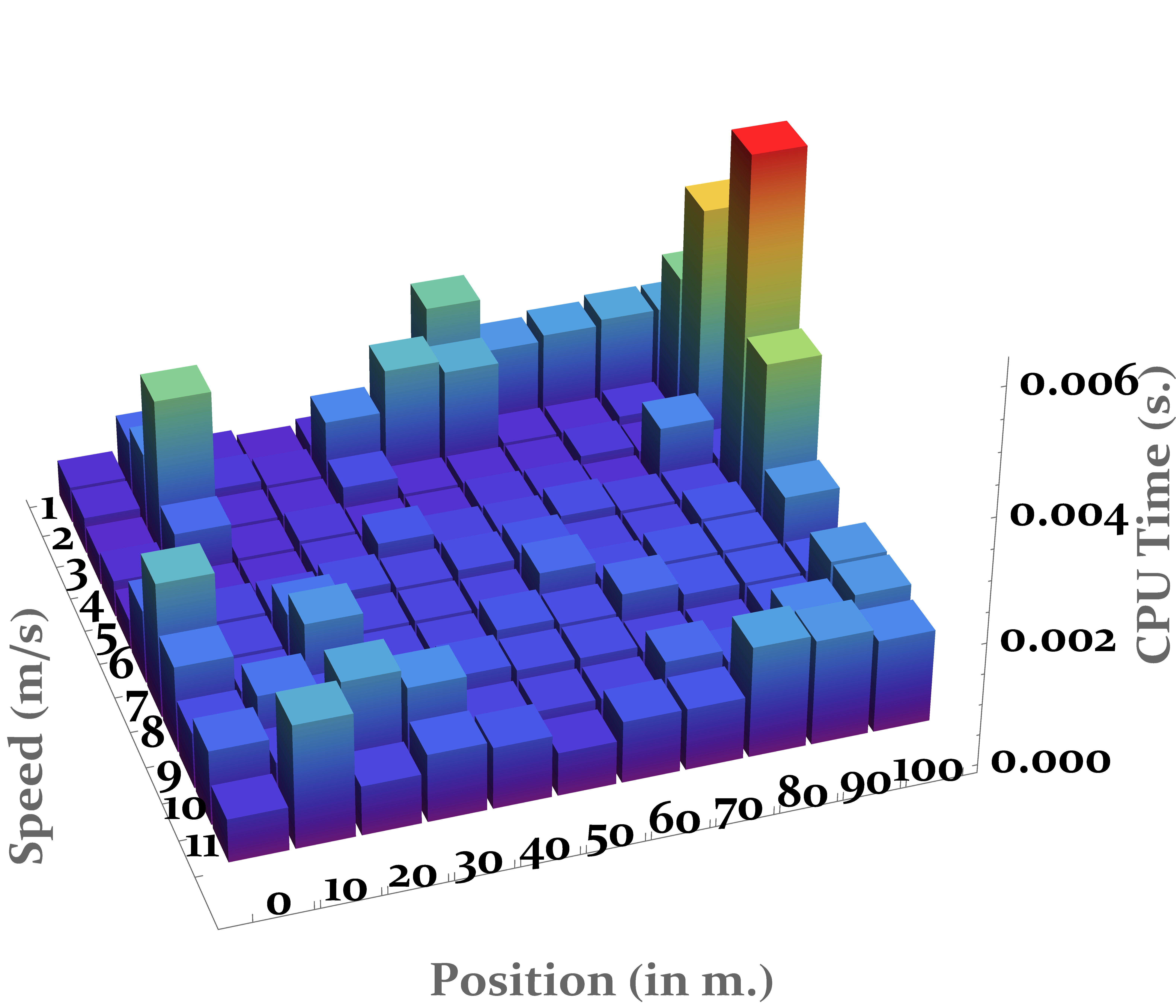}
		\caption{}
		\label{fig:cpuTimes-2}
	\end{subfigure}
\caption{Computation times until convergence of: (a) DDDP and (b) DDP algorithms for different initial conditions.}
\label{fig:cpuTimes}
\end{figure}

In order to explore the robustness of the DDDP and DDP algorithms and to further emphasize their superiority compared to the one-shot (S)DP, in terms of computation times, a more excessive analysis was conducted and is presented in Fig. \ref{fig:cpuTimes}. Specifically, Fig. \ref{fig:cpuTimes} illustrates the required CPU times until convergence for both proposed algorithms for a plethora of different initial conditions, i.e., initial speeds and positions (distance from the traffic signal) of a vehicle. From both Fig. \ref{fig:cpuTimes-1} and Fig. \ref{fig:cpuTimes-2} it can be seen that, for all evaluated initial conditions, the CPU times of both DDDP and DDP are similarly low as for the above three scenarios, with average times of 1.4 s and 0.001 s, respectively, the corresponding maximum and minimum times being 6.3 s and 0.2 s for DDDP and 0.006 s and 0.0005 s for DDP, respectively. These findings lead to the conclusion that, even though DDDP reduces the CPU time significantly compared to the one-shot SDP, the DDP algorithm achieves almost instantaneous solutions, which enables the DDP approach to be used, within an MPC framework, in the vehicle's on-board computer.

\section{Conclusions} \label{sec:concl}

In recent works \citep{typaldos2020vehicle, typaldos2021vehicle}, a stochastic GLOSA methodology was developed, by optimizing, using SDP and DDDP techniques, respectively, the vehicle kinematic trajectories subject to the intermediate stochastic traffic signal switching constraint and with a fixed final state. As an extension to these works, an iterative Differential Dynamic Programming (DDP) algorithm was employed in this paper. Within each iteration, DDP solves, for each time-step separately, a quadratic-linear approximation of the recursive Bellman equation. Demonstration and comparison results indicate that the DDP algorithm strongly outperforms both the full-range SDP and the DDDP algorithms in terms of computation time; and even delivers a slightly better optimal solutions due to lack of discretisation. This enables the DDP algorithm to be readily executable within a model predictive control framework in real time on-board approaching vehicles.

Current and future work is focused on:

\begin{itemize}
\item Generalization of the current GLOSA problem by considering uncertain switching times for both green and red phases.
\item Testing of the approach performance in presence of multiple (equipped and not equipped) vehicles, while traversing multiple subsequent intersections.
\end{itemize}

\section*{Acknowledgments}
The research leading to these results has received funding form the European Research Council under the European Union's Horizon 2020 Research and Innovation programme/ERC Grant Agreement n. [833915], project TrafficFluid.

\newpage

\section*{Appendix A}


\begin{figure}[h!t]
\centering
\begin{minipage}{\linewidth}
	\begin{subfigure}{0.3\linewidth} 
		\includegraphics[width=\textwidth]{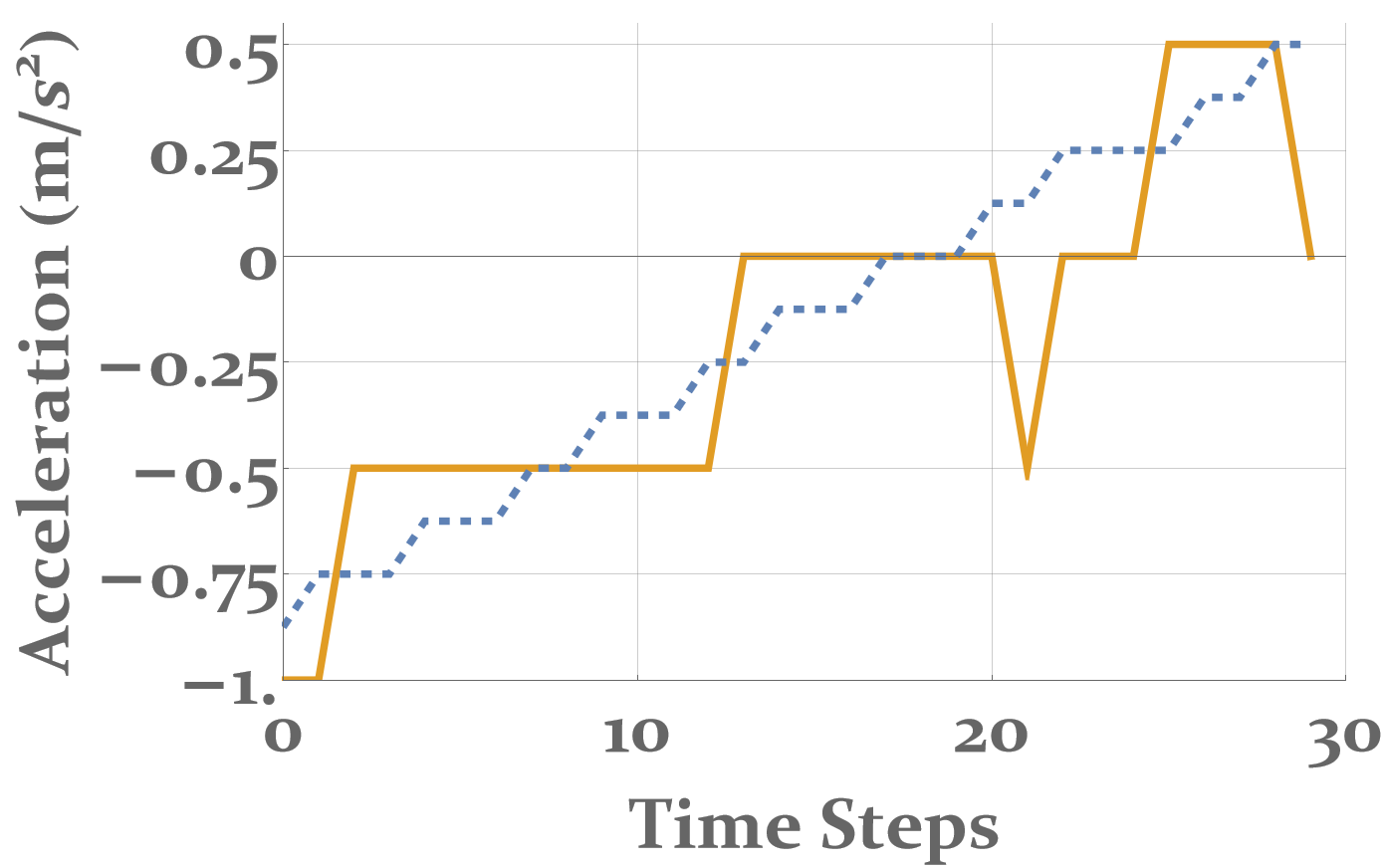}
		\caption*{Iteration 1}
	\end{subfigure}%
	\begin{subfigure}{0.3\linewidth} 
		\includegraphics[width=\textwidth]{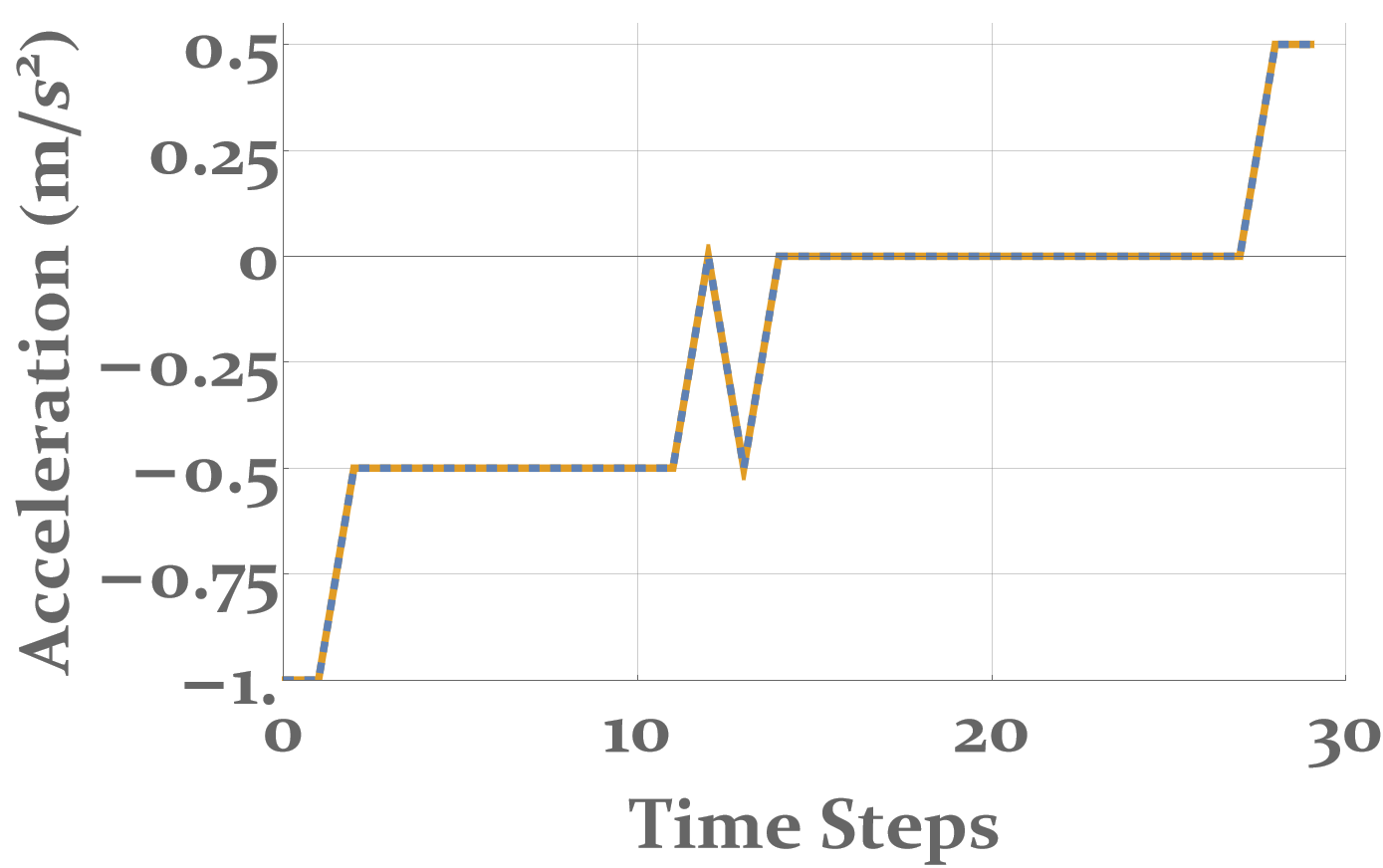}
		\caption*{Iteration 3}
	\end{subfigure}%
	\begin{subfigure}{0.3\linewidth}
		\includegraphics[width=\textwidth]{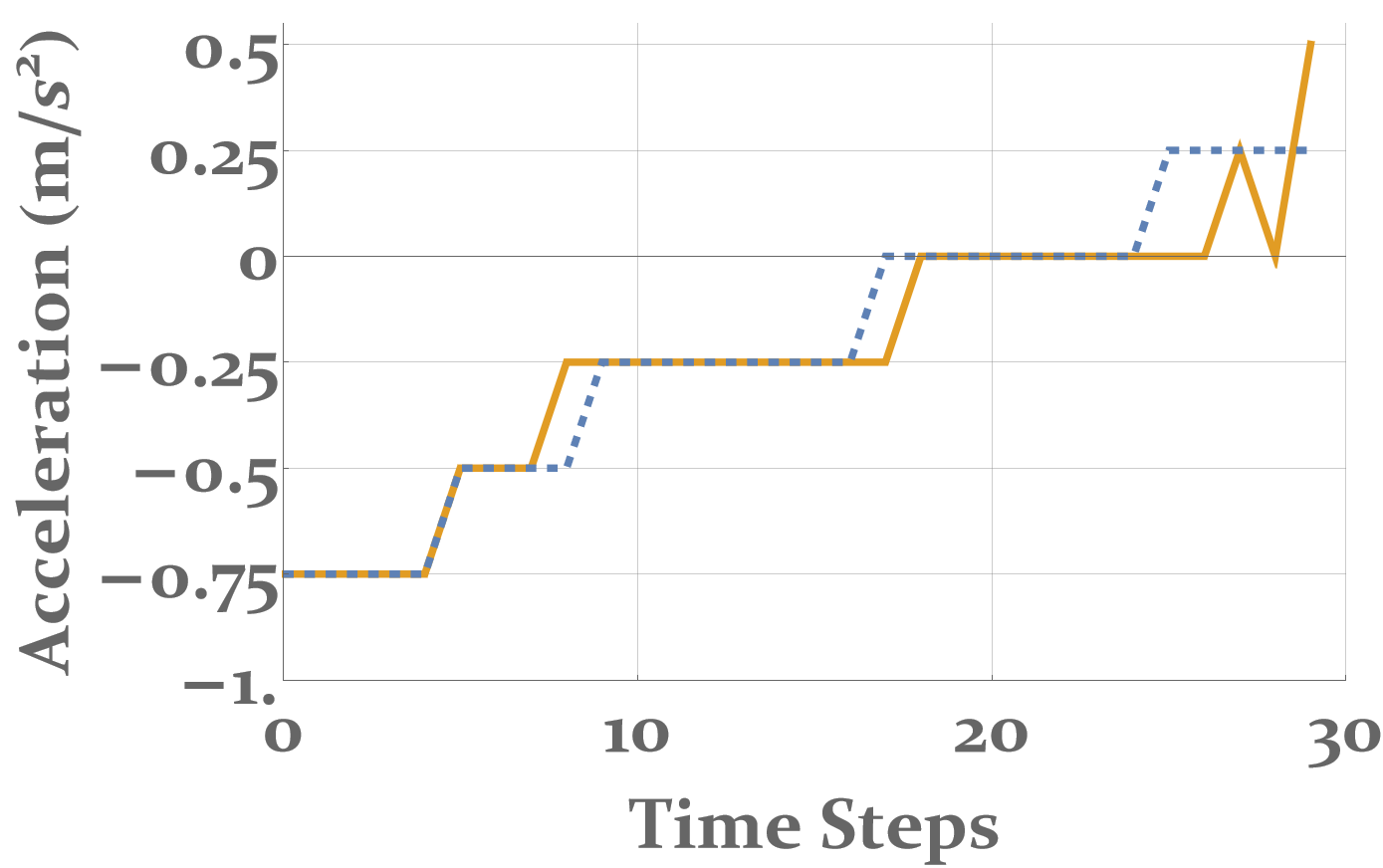}
		\caption*{Iteration 5}
	\end{subfigure}
	\begin{subfigure}{0.3\linewidth}
		\includegraphics[width=\textwidth]{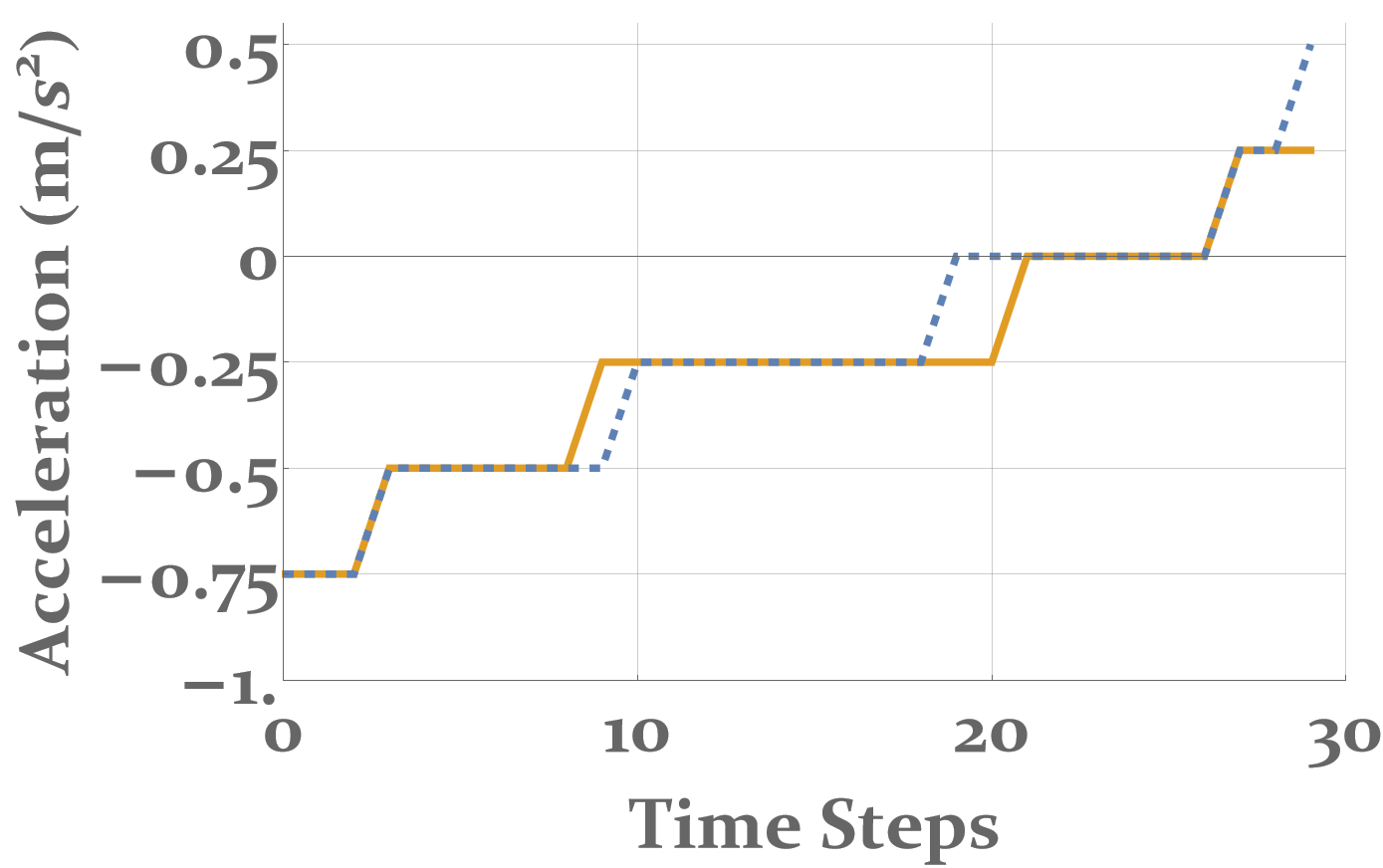}
		\caption*{Iteration 7}
	\end{subfigure}%
	\begin{subfigure}{0.3\linewidth}
		\includegraphics[width=\textwidth]{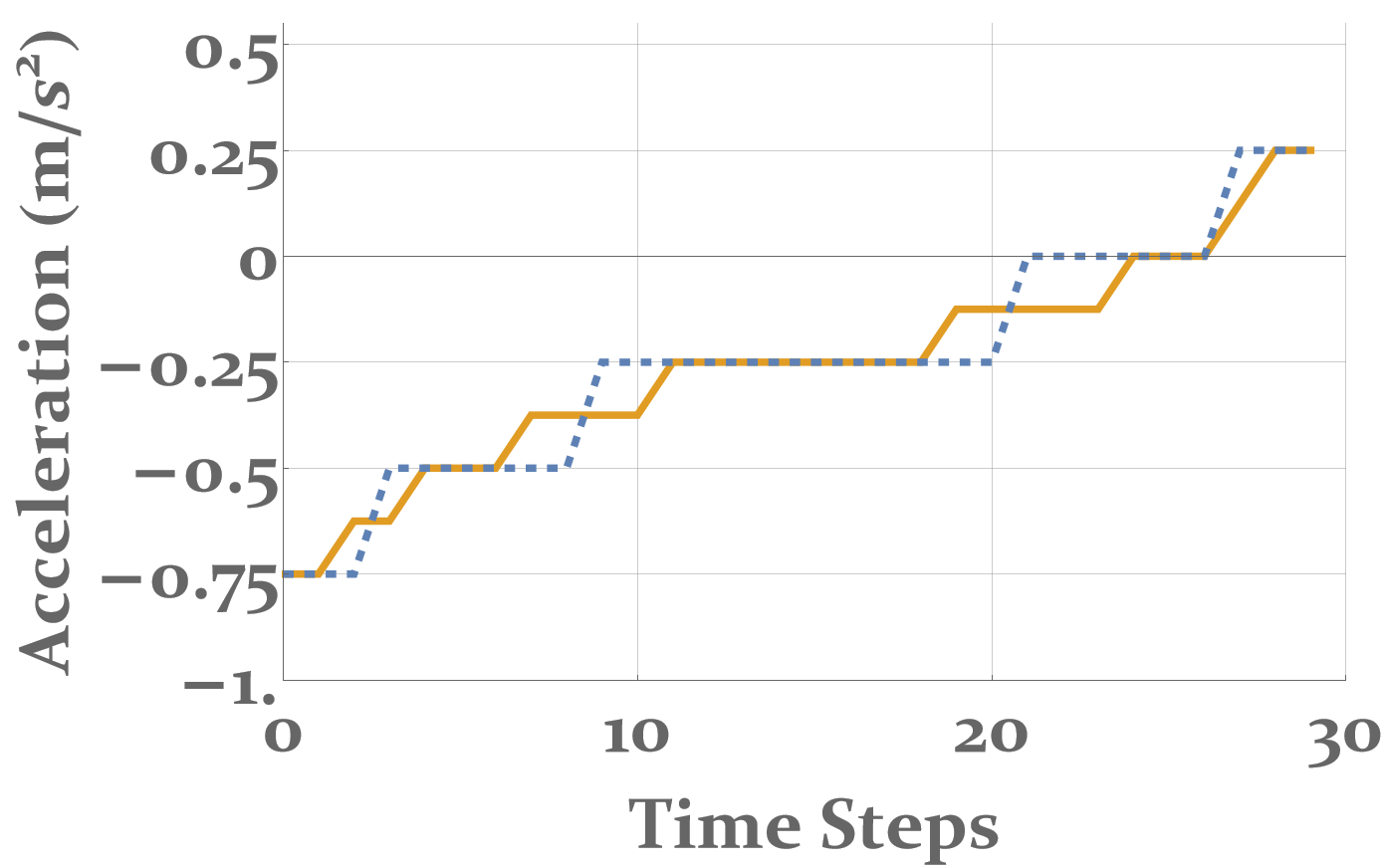}
		\caption*{Iteration 9}
	\end{subfigure}%
	\begin{subfigure}{0.3\linewidth}
		\includegraphics[width=\textwidth]{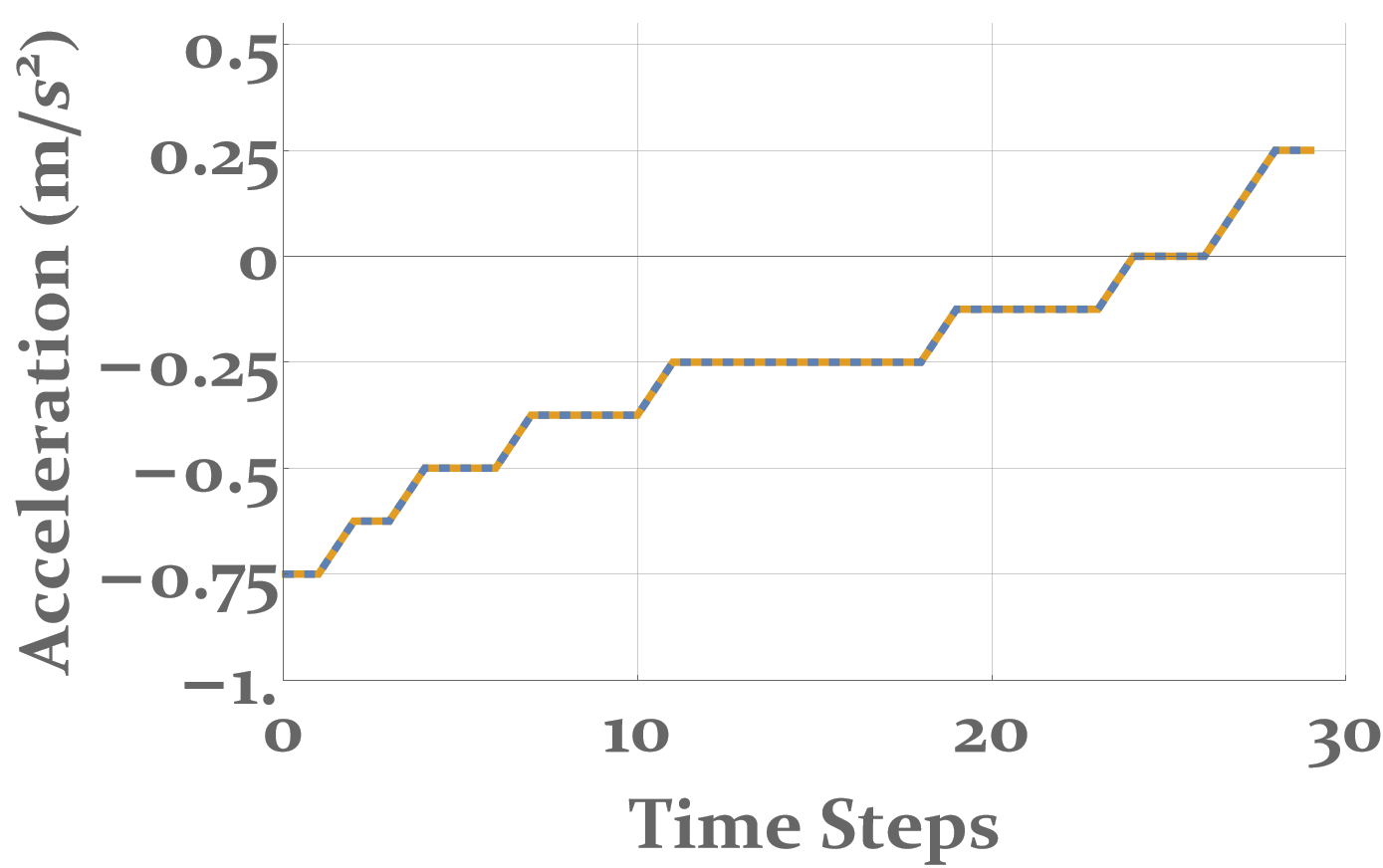}
		\caption*{Iteration 10}
	\end{subfigure}
	\subcaption{}
\end{minipage}
\begin{minipage}{\linewidth}
	\begin{subfigure}{0.3\linewidth} 
		\includegraphics[width=\textwidth]{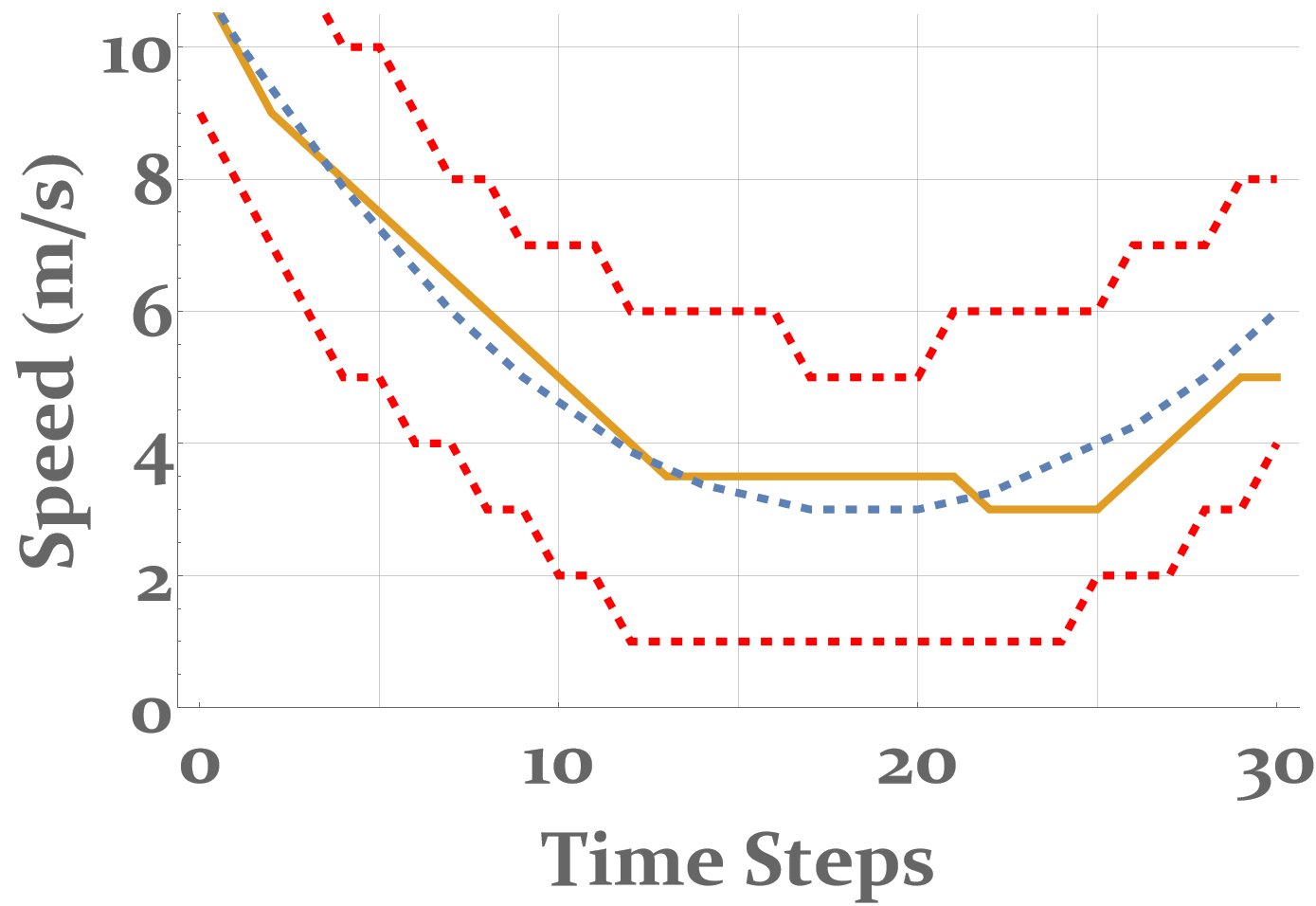}
		\caption*{Iteration 1}
	\end{subfigure}%
	\begin{subfigure}{0.3\linewidth} 
		\includegraphics[width=\textwidth]{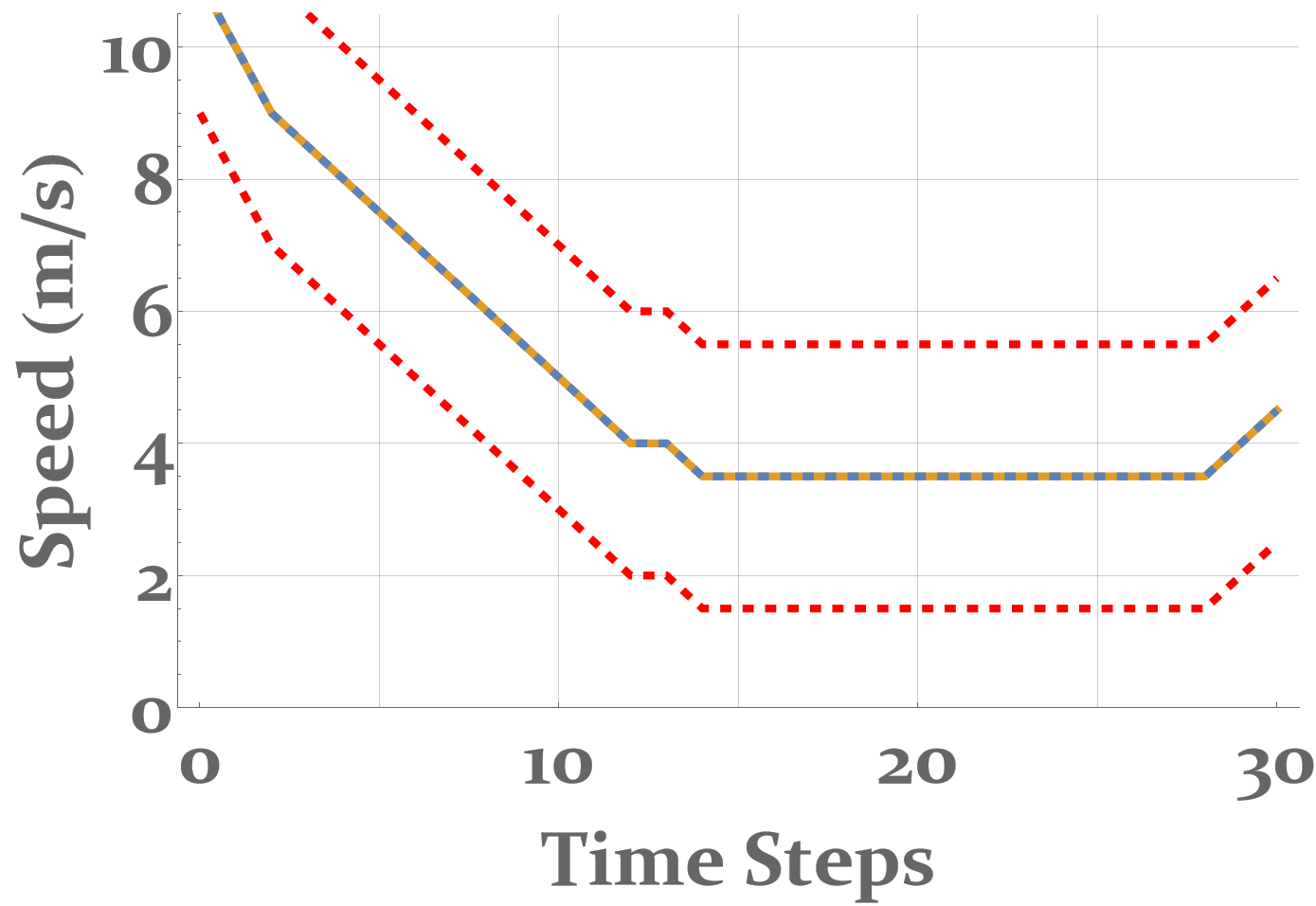}
		\caption*{Iteration 3}
	\end{subfigure}%
	\begin{subfigure}{0.3\linewidth}
		\includegraphics[width=\textwidth]{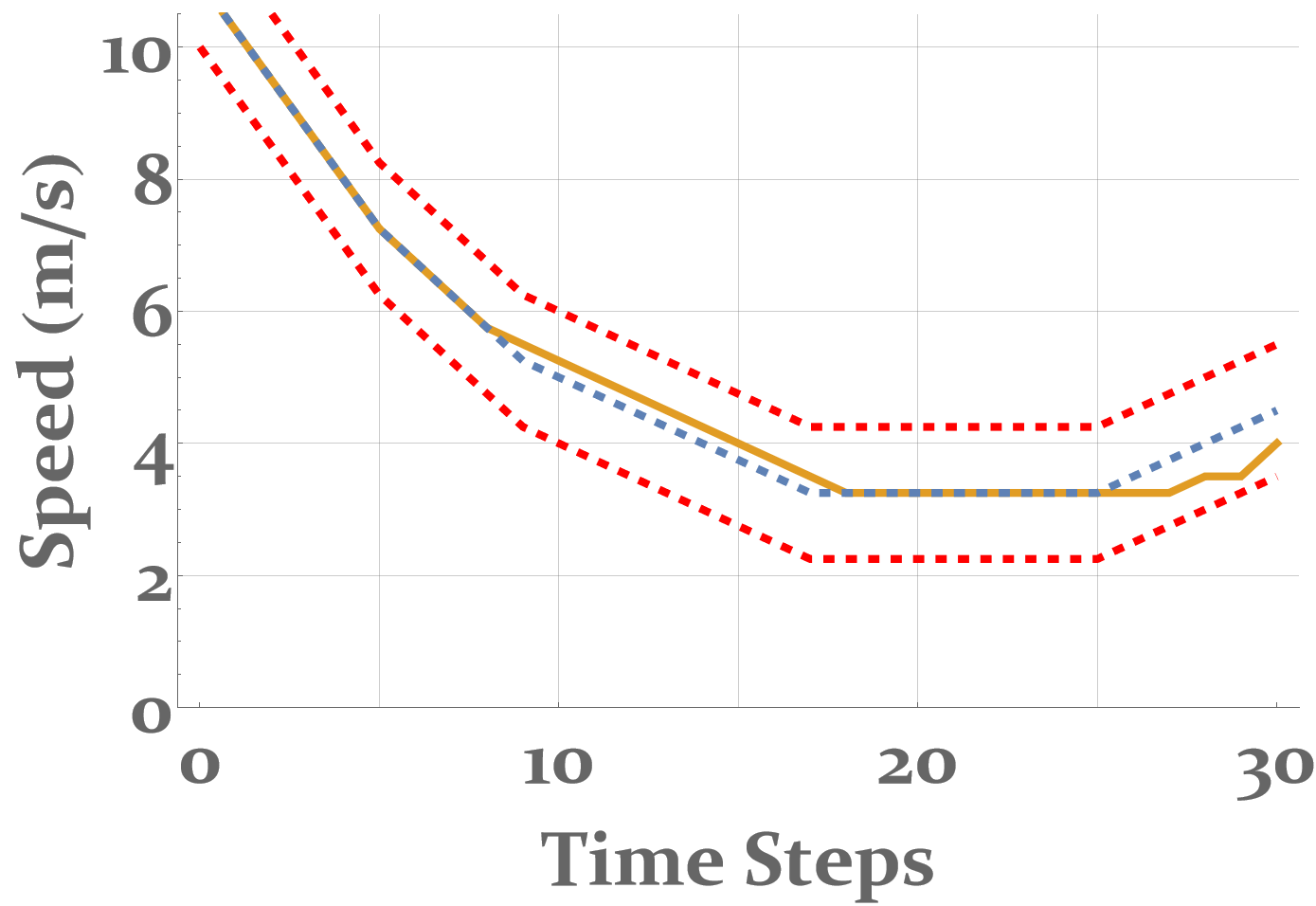}
		\caption*{Iteration 5}
	\end{subfigure}
	\begin{subfigure}{0.3\linewidth}
		\includegraphics[width=\textwidth]{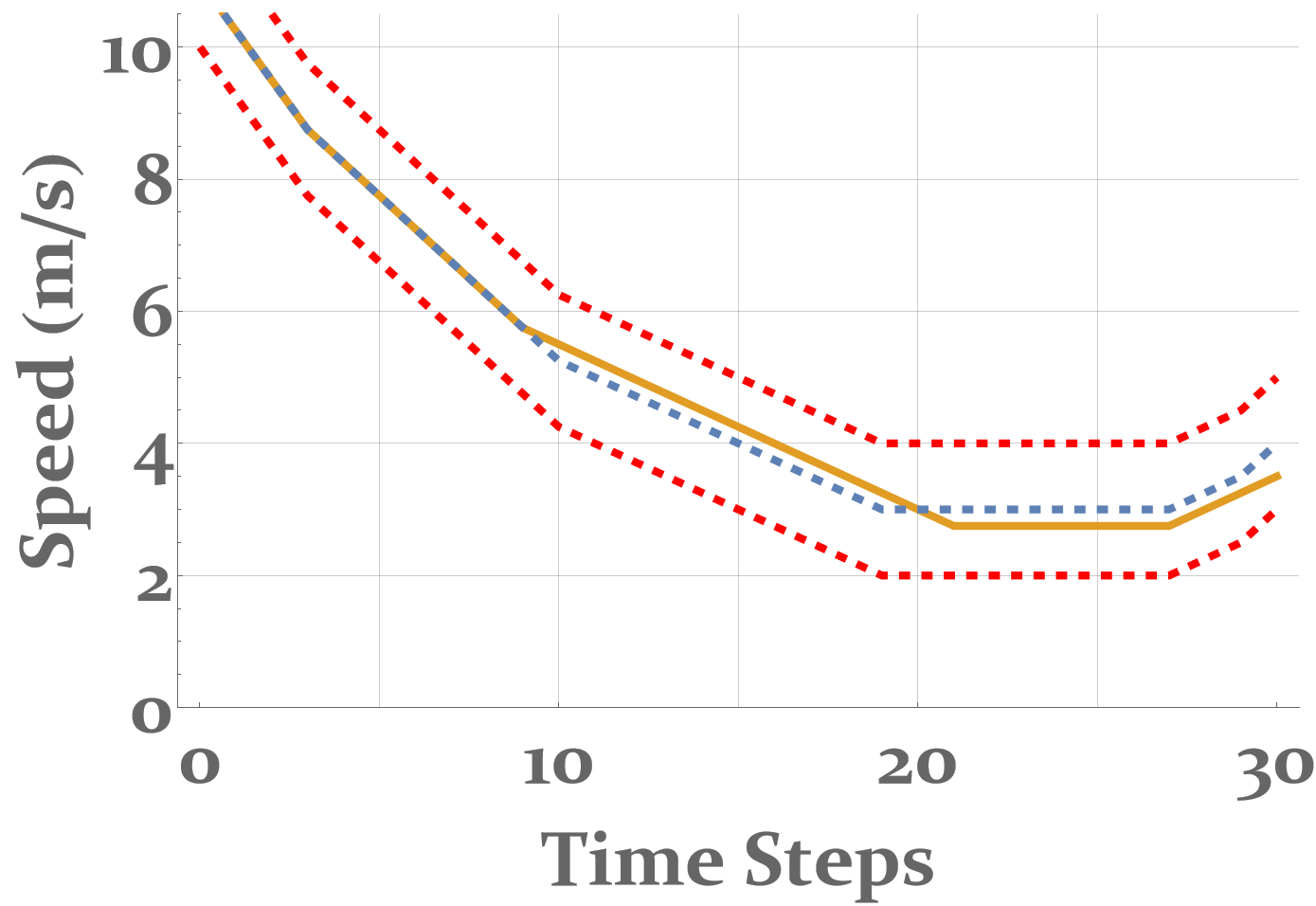}
		\caption*{Iteration 7}
	\end{subfigure}%
	\begin{subfigure}{0.3\linewidth}
		\includegraphics[width=\textwidth]{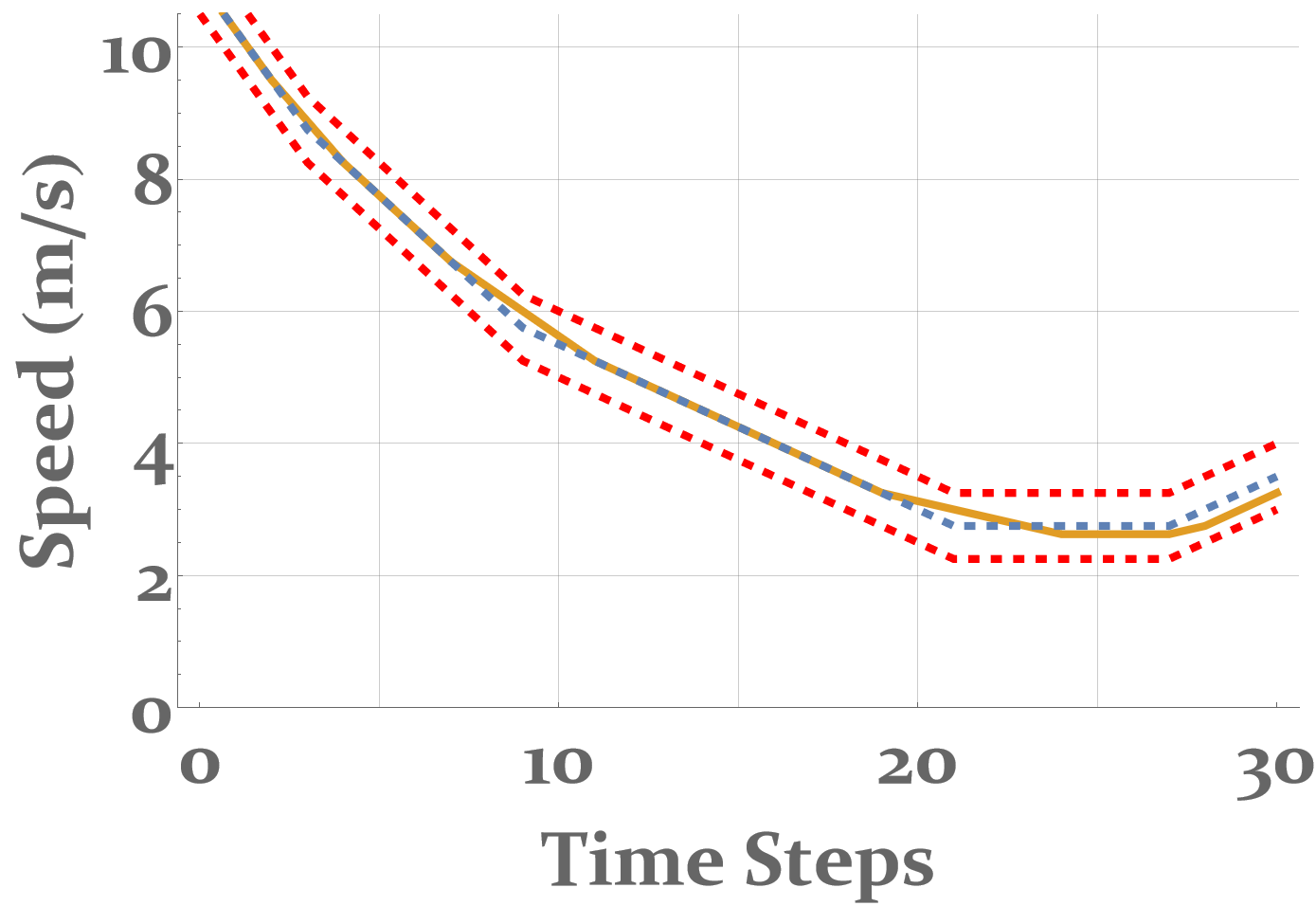}
		\caption*{Iteration 9}
	\end{subfigure}%
	\begin{subfigure}{0.3\linewidth}
		\includegraphics[width=\textwidth]{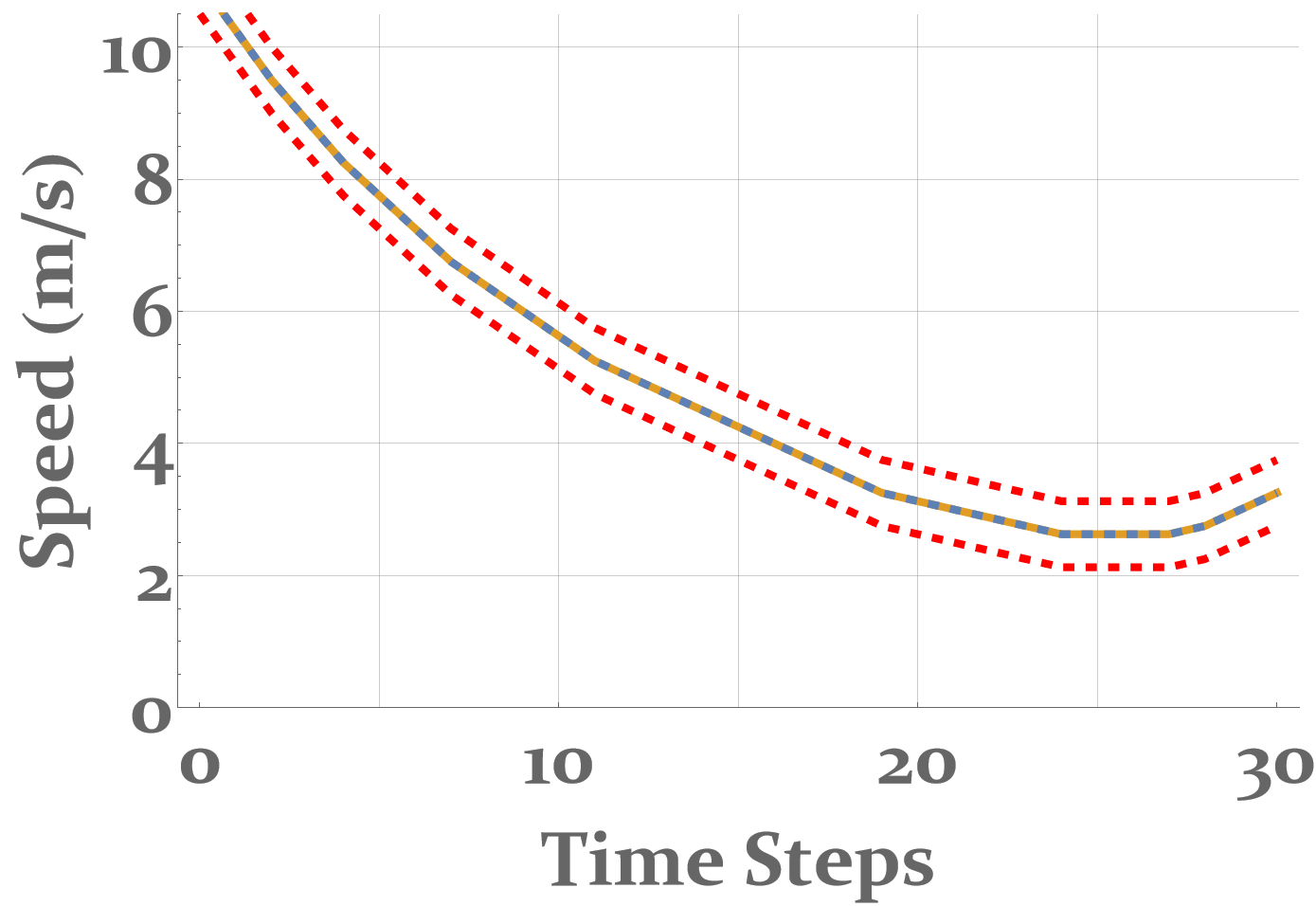}
		\caption*{Iteration 10}
	\end{subfigure}
	\subcaption{}
\end{minipage}
\caption{Received (blue dashed line) and optimal (orange line) (a) acceleration and (b) speed trajectories of DDDP algorithm in each iteration for Scenario 2.} 
\label{fig:DDDP-s-2}
\end{figure}

\vspace*{\fill}
\newpage


\begin{figure}[h!t]
\centering
\begin{minipage}{\linewidth}
	\begin{subfigure}{0.3\linewidth} 
		\includegraphics[width=\textwidth]{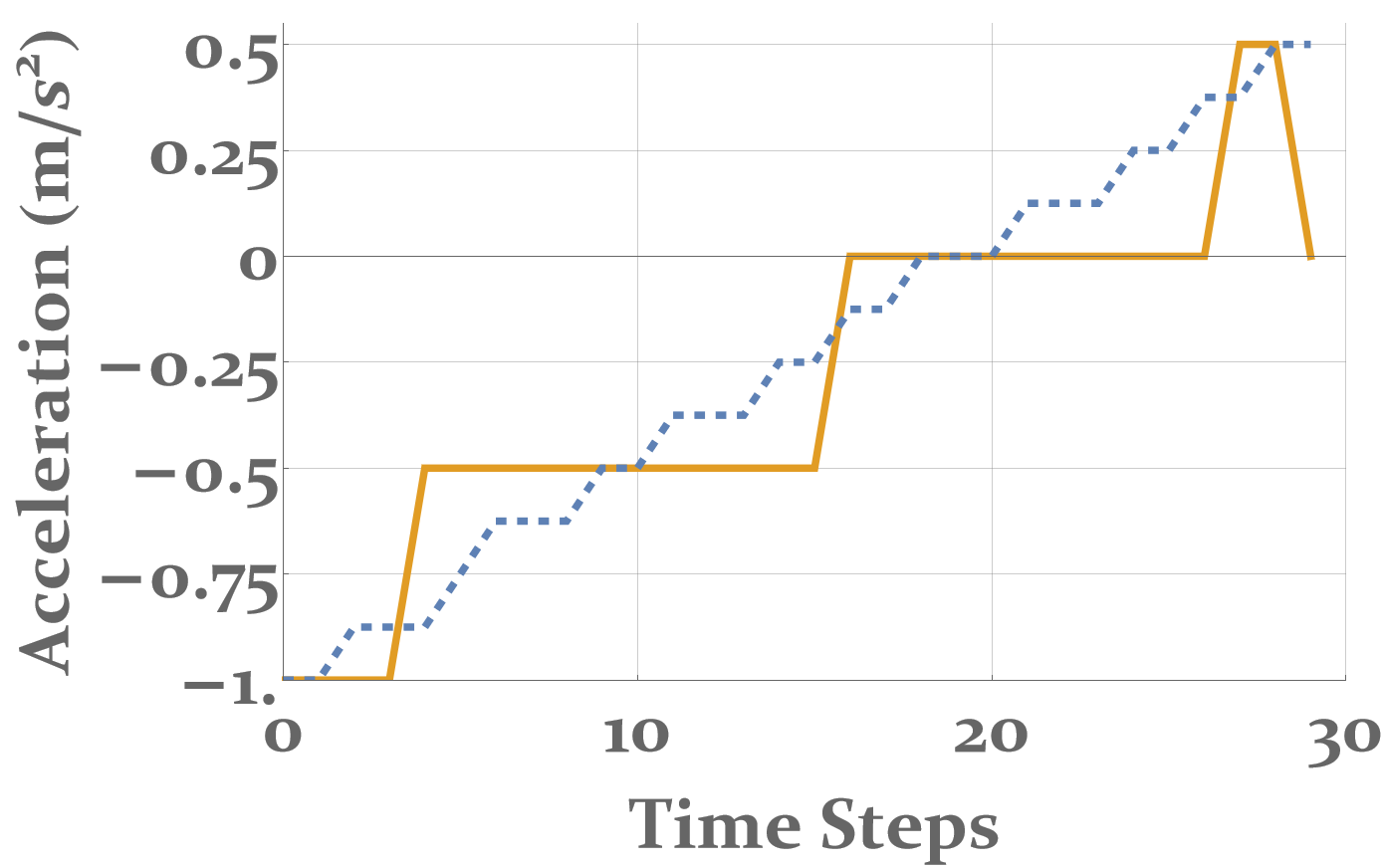}
		\caption*{Iteration 1}
	\end{subfigure}%
	\begin{subfigure}{0.3\linewidth} 
		\includegraphics[width=\textwidth]{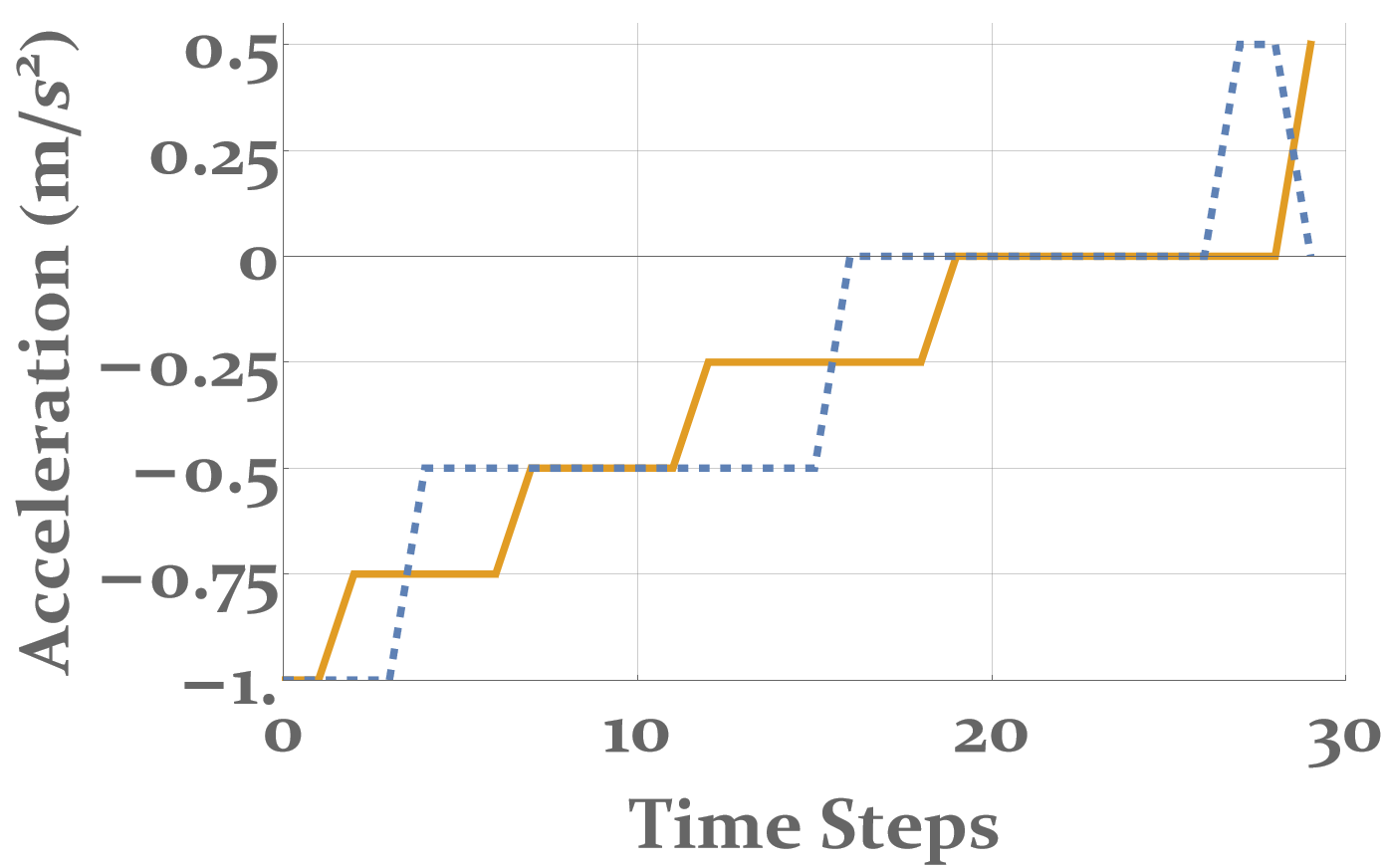}
		\caption*{Iteration 3}
	\end{subfigure}%
	\begin{subfigure}{0.3\linewidth}
		\includegraphics[width=\textwidth]{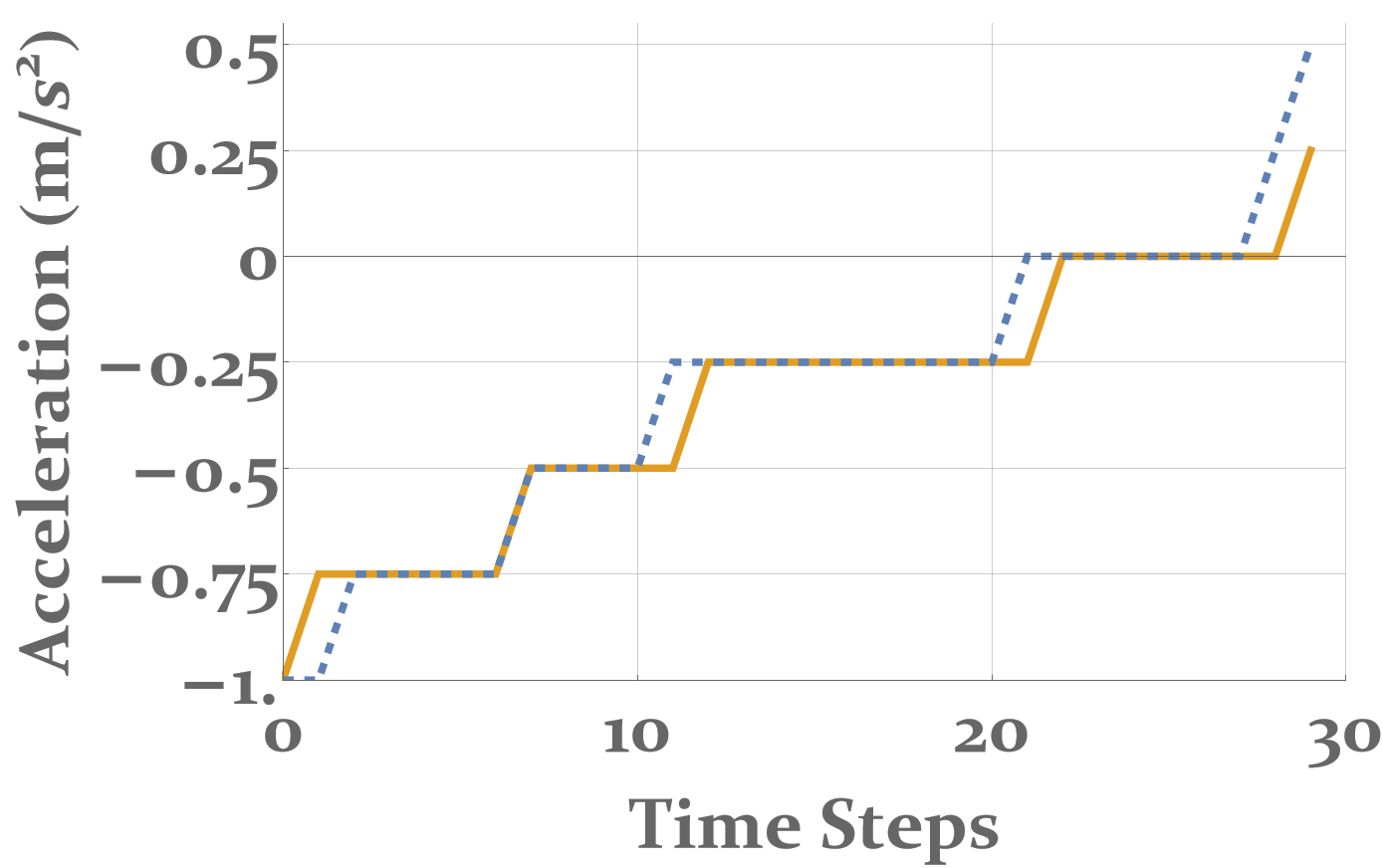}
		\caption*{Iteration 5}
	\end{subfigure}
	\begin{subfigure}{0.3\linewidth}
		\includegraphics[width=\textwidth]{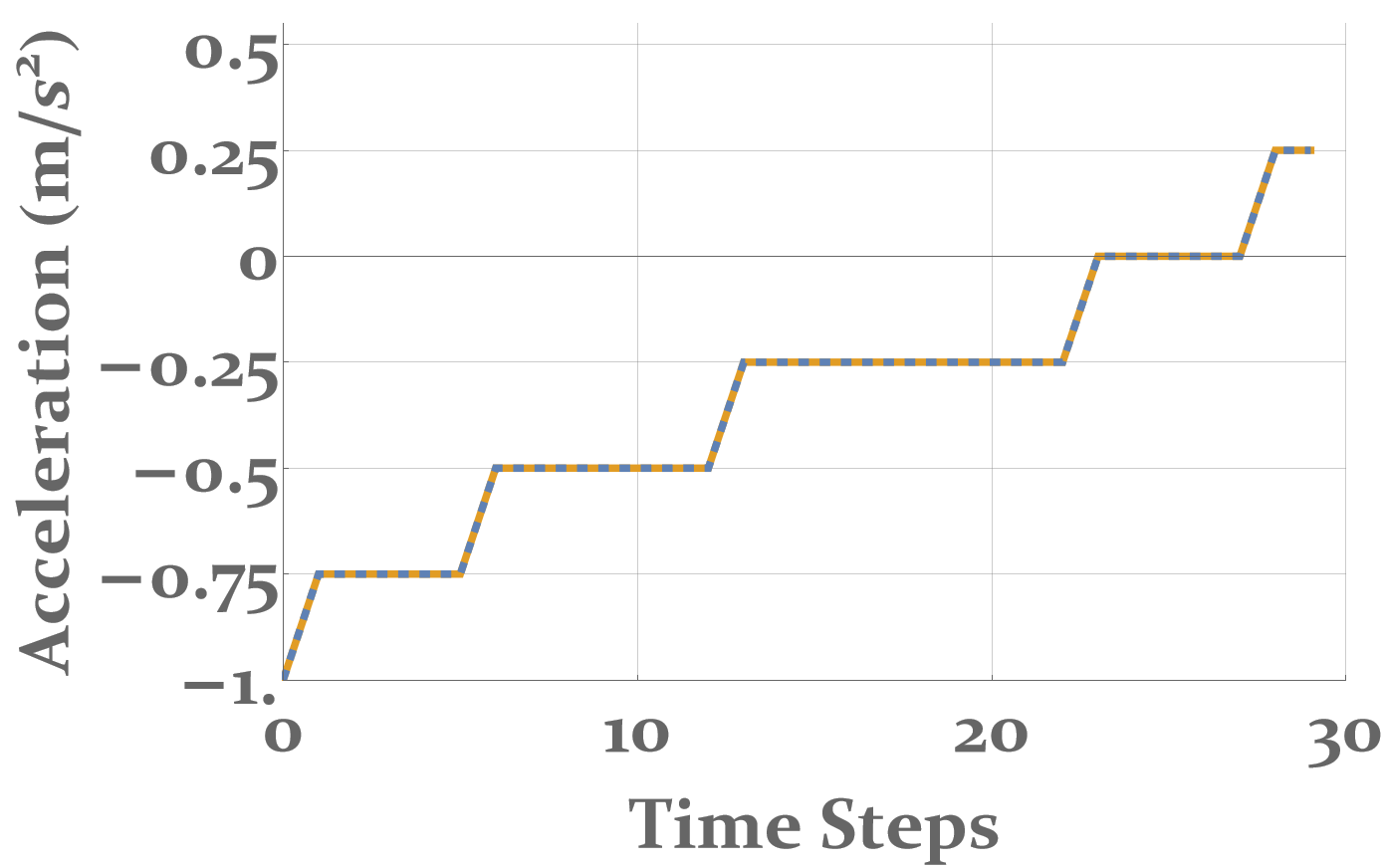}
		\caption*{Iteration 7}
	\end{subfigure}%
	\begin{subfigure}{0.3\linewidth}
		\includegraphics[width=\textwidth]{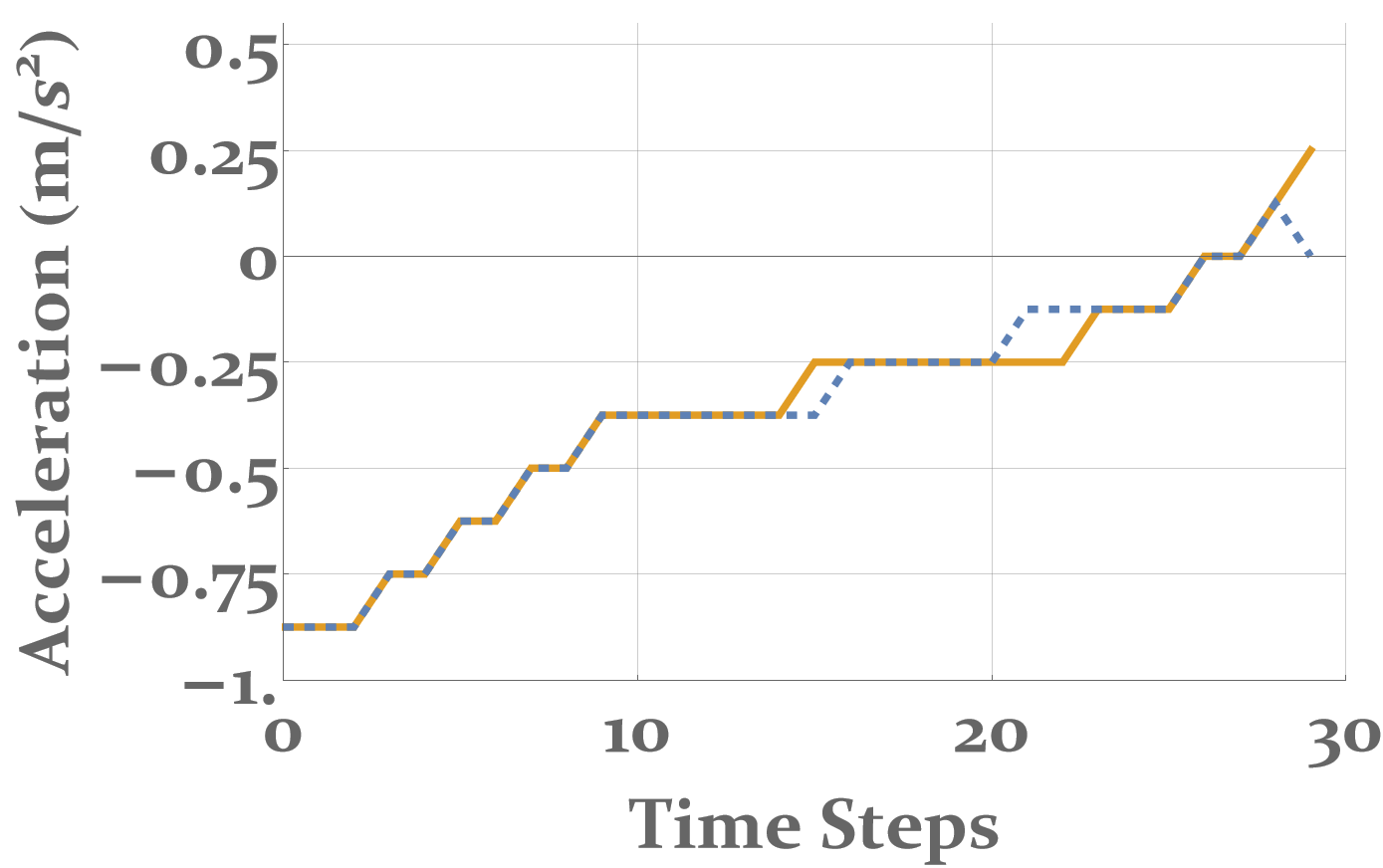}
		\caption*{Iteration 9}
	\end{subfigure}%
	\begin{subfigure}{0.3\linewidth}
		\includegraphics[width=\textwidth]{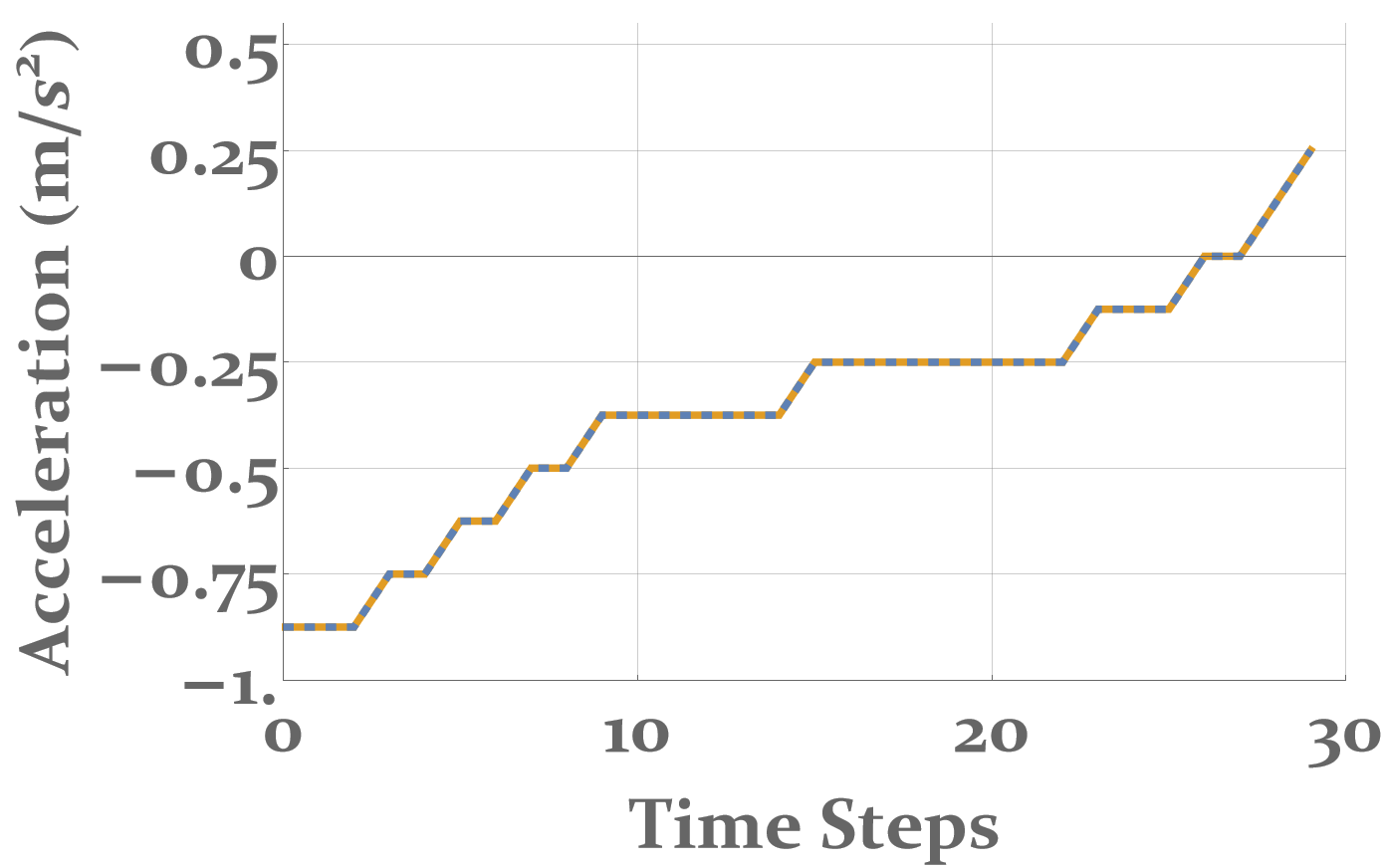}
		\caption*{Iteration 10}
	\end{subfigure}
	\subcaption{}
\end{minipage}
\begin{minipage}{\linewidth}
	\begin{subfigure}{0.3\linewidth} 
		\includegraphics[width=\textwidth]{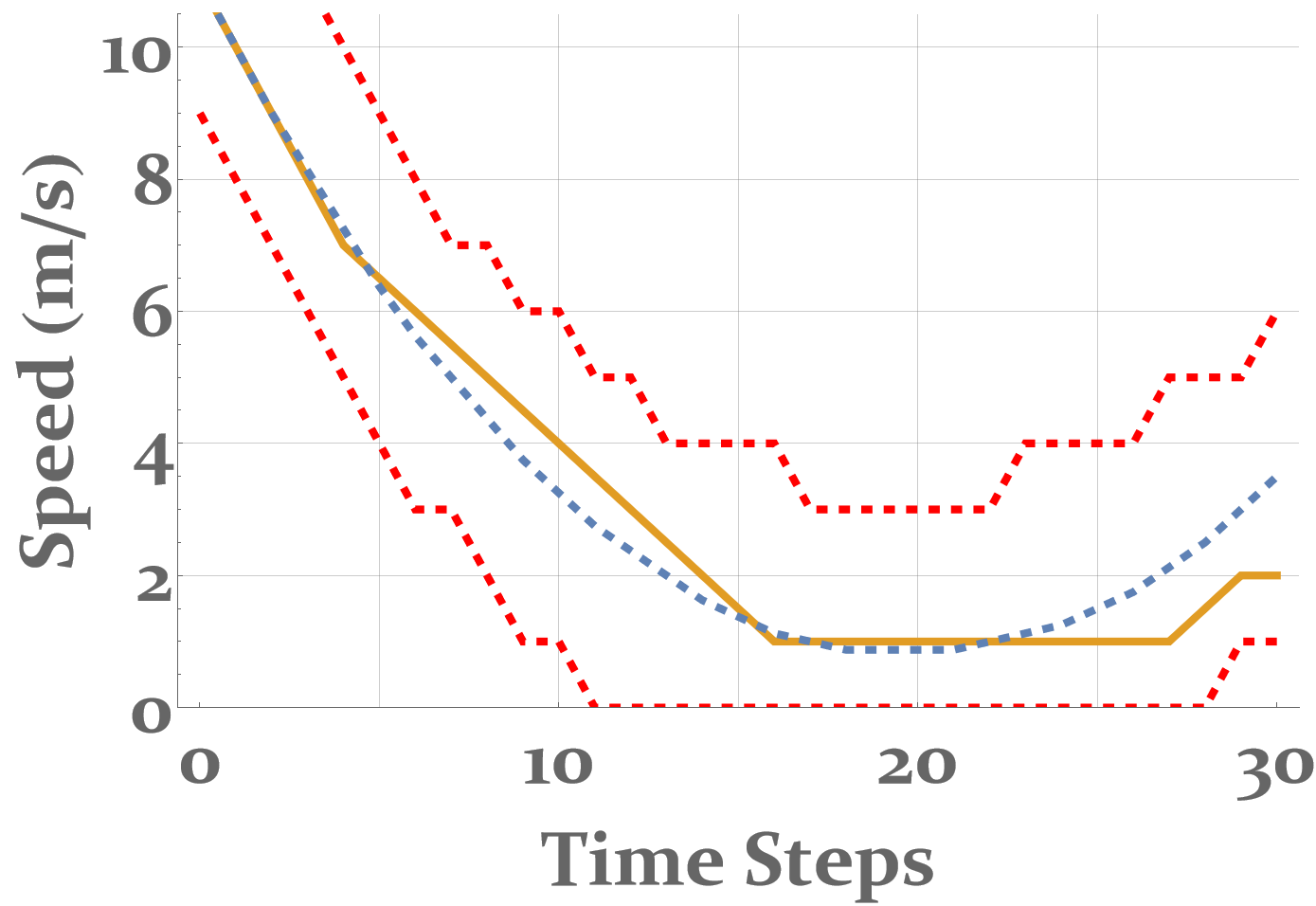}
		\caption*{Iteration 1}
	\end{subfigure}%
	\begin{subfigure}{0.3\linewidth} 
		\includegraphics[width=\textwidth]{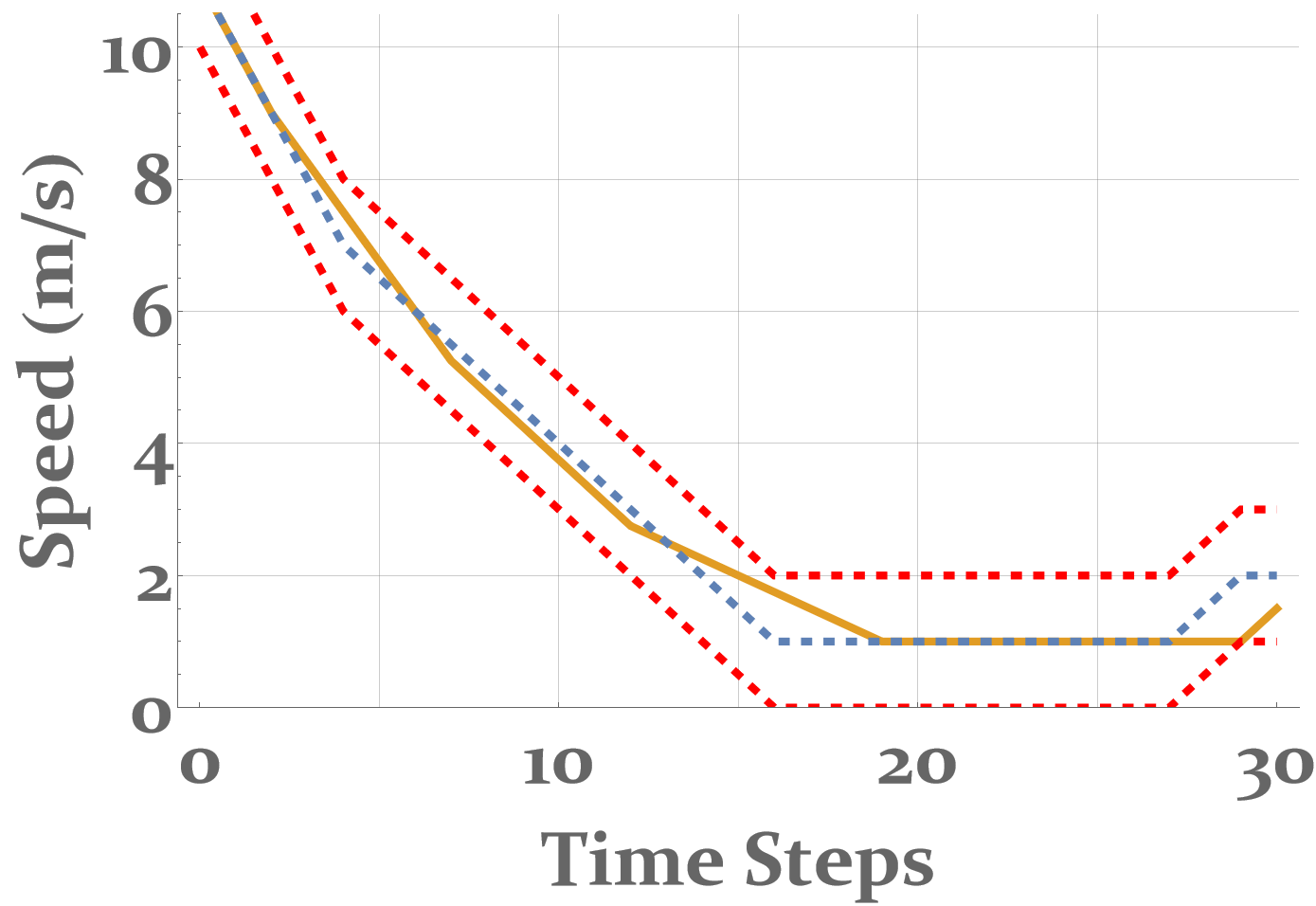}
		\caption*{Iteration 3}
	\end{subfigure}%
	\begin{subfigure}{0.3\linewidth}
		\includegraphics[width=\textwidth]{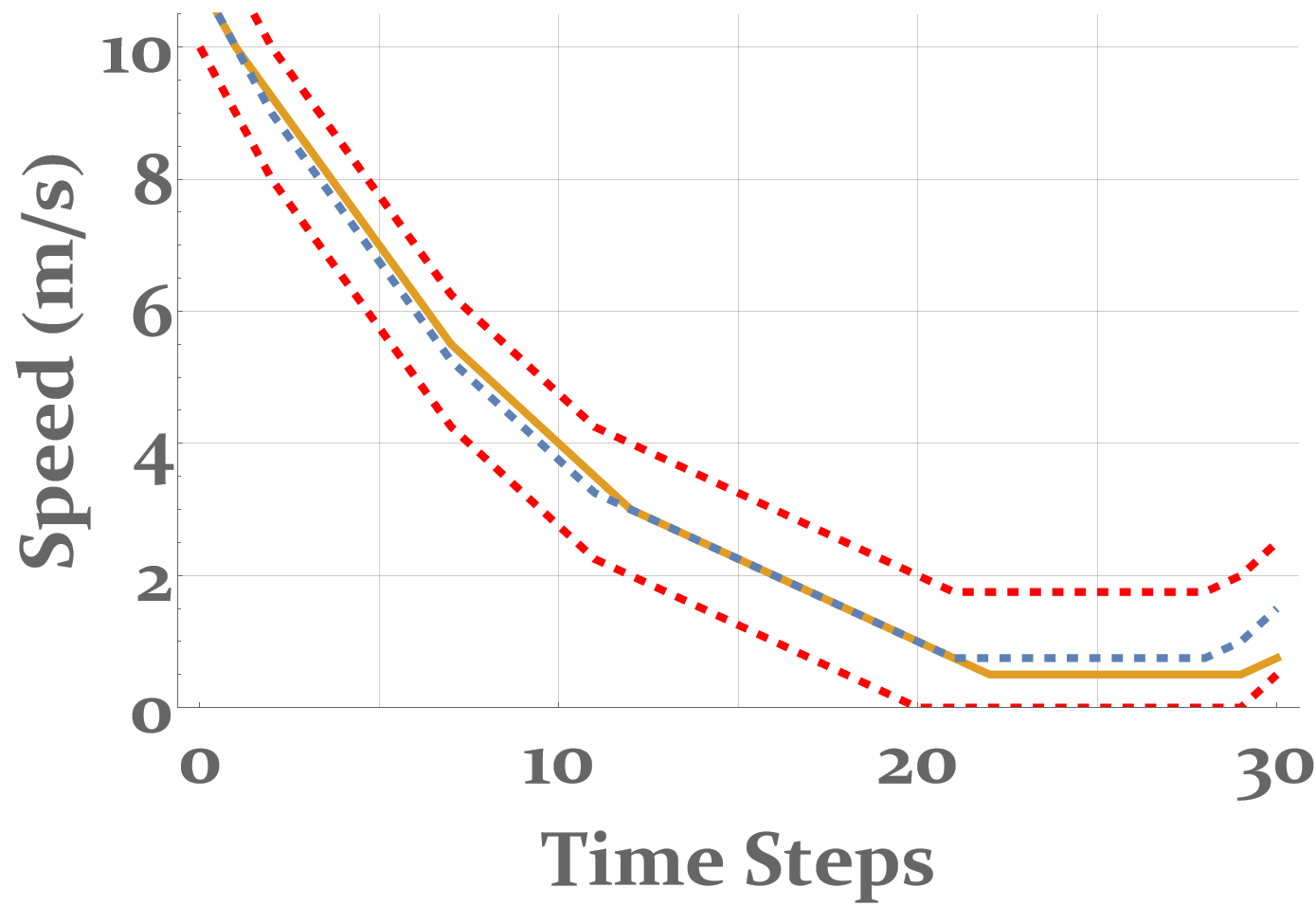}
		\caption*{Iteration 5}
	\end{subfigure}
	\begin{subfigure}{0.3\linewidth}
		\includegraphics[width=\textwidth]{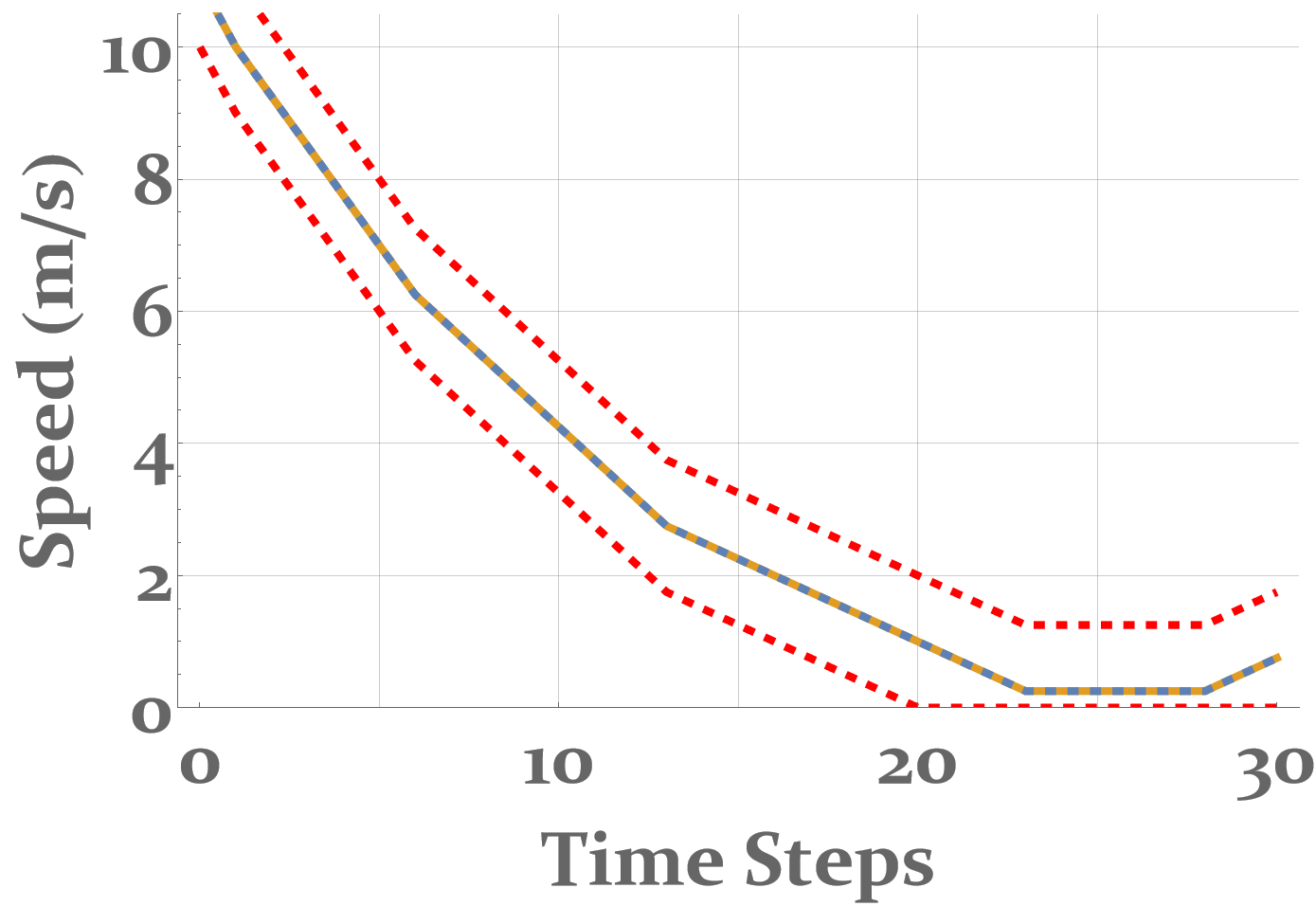}
		\caption*{Iteration 7}
	\end{subfigure}%
	\begin{subfigure}{0.3\linewidth}
		\includegraphics[width=\textwidth]{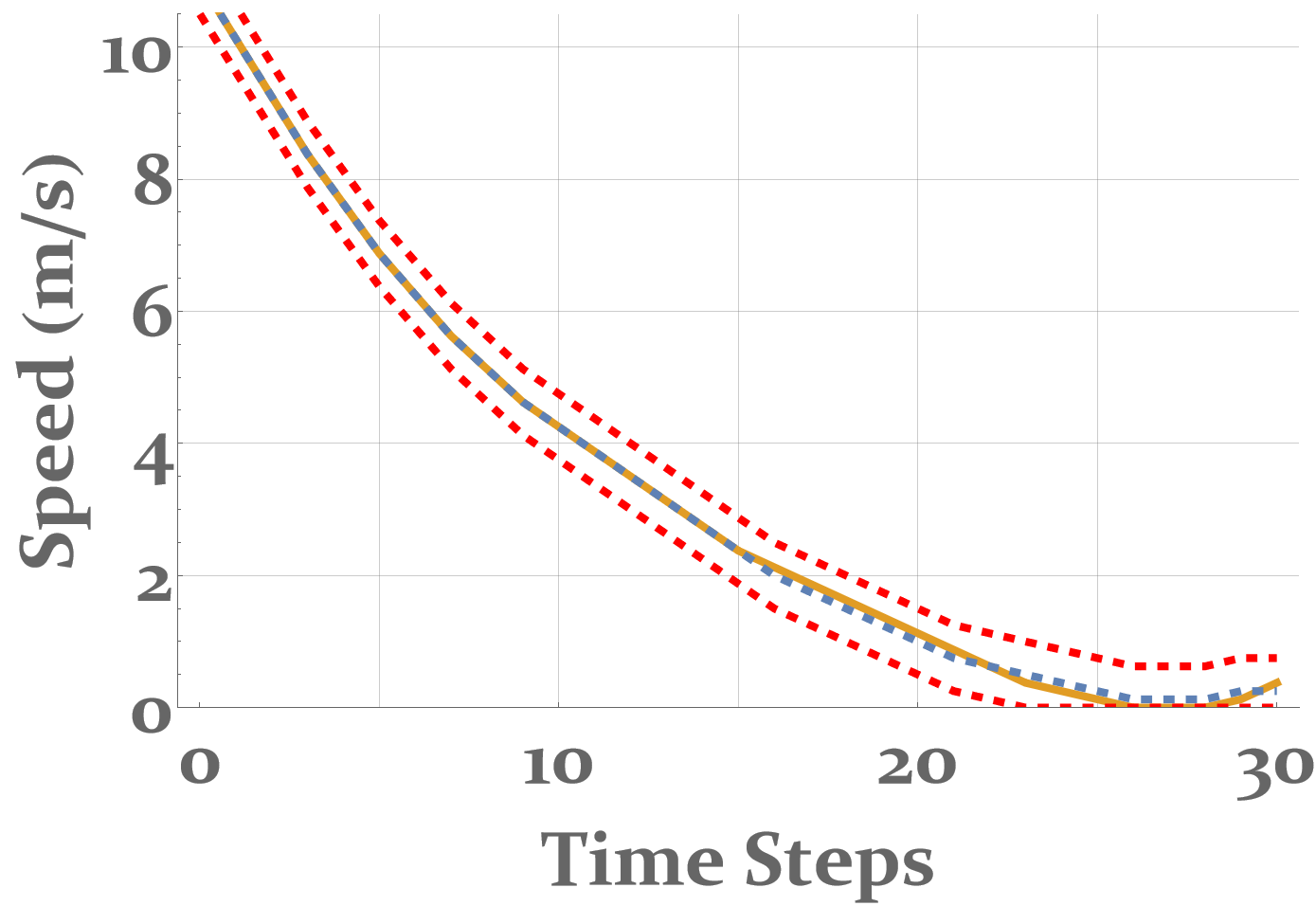}
		\caption*{Iteration 9}
	\end{subfigure}%
	\begin{subfigure}{0.3\linewidth}
		\includegraphics[width=\textwidth]{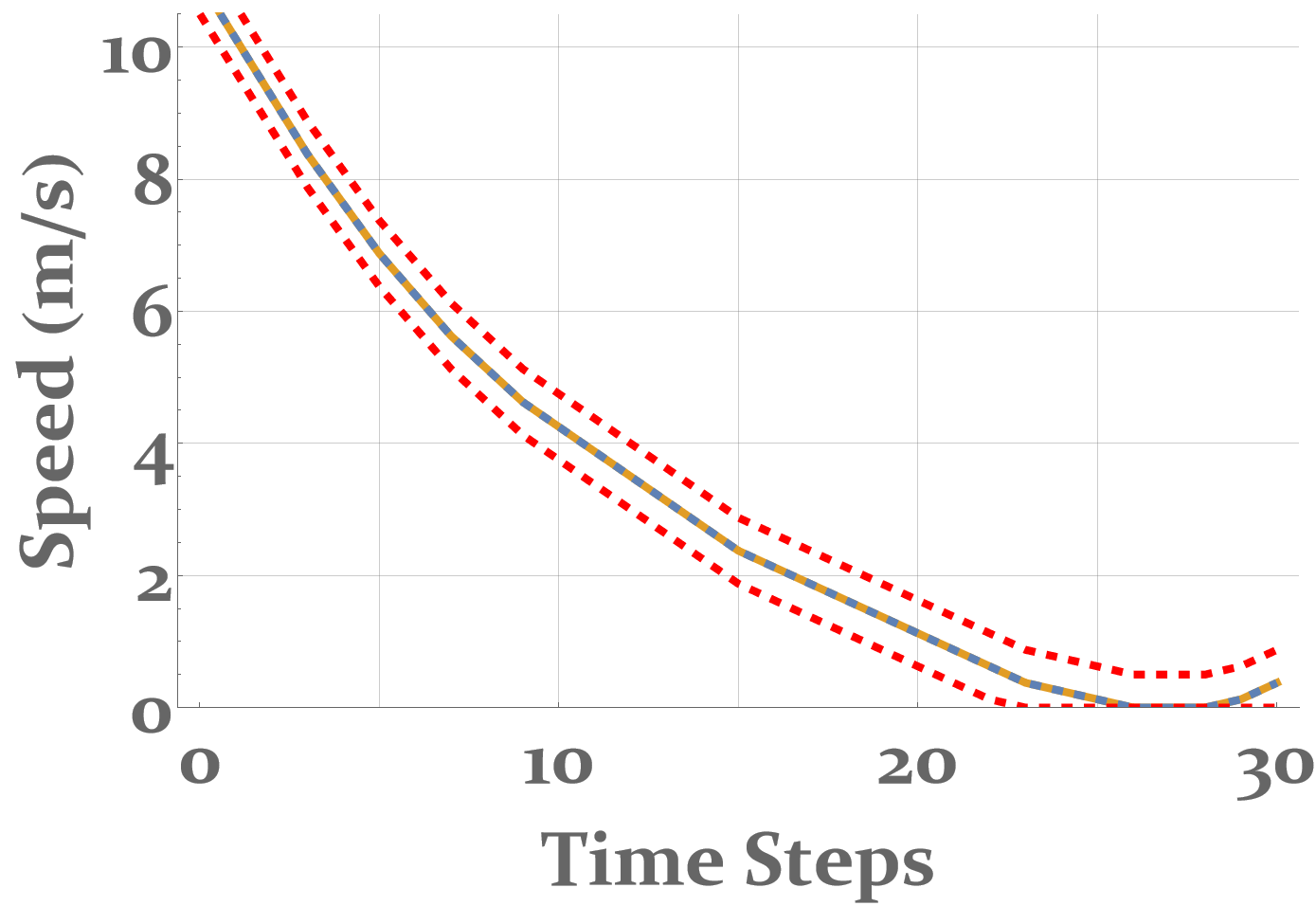}
		\caption*{Iteration 10}
	\end{subfigure}
	\subcaption{}
\end{minipage}
\caption{Received (blue dashed line) and optimal (orange line) (a) acceleration and (b) speed trajectories of DDDP algorithm in each iteration for Scenario 3.} 
\label{fig:DDDP-s-3}
\end{figure}

\vspace*{\fill}
\newpage

\begin{figure}[h!t]
\centering
\begin{minipage}[b]{\textwidth}
	\begin{subfigure}{0.25\textwidth} 
		\includegraphics[width=\textwidth]{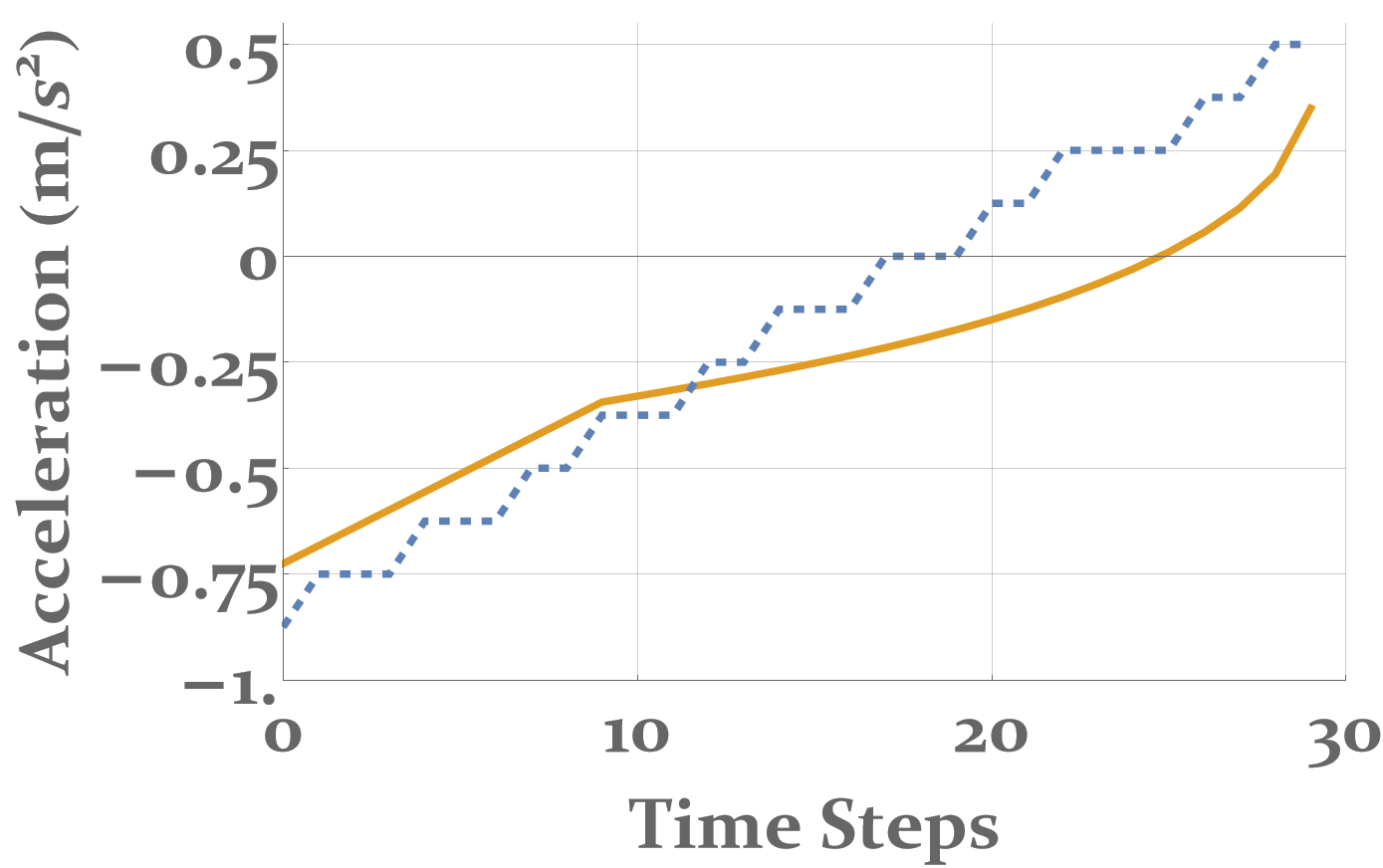}
		\caption*{Iteration 1}
	\end{subfigure}%
	\begin{subfigure}{0.25\textwidth} 
		\includegraphics[width=\textwidth]{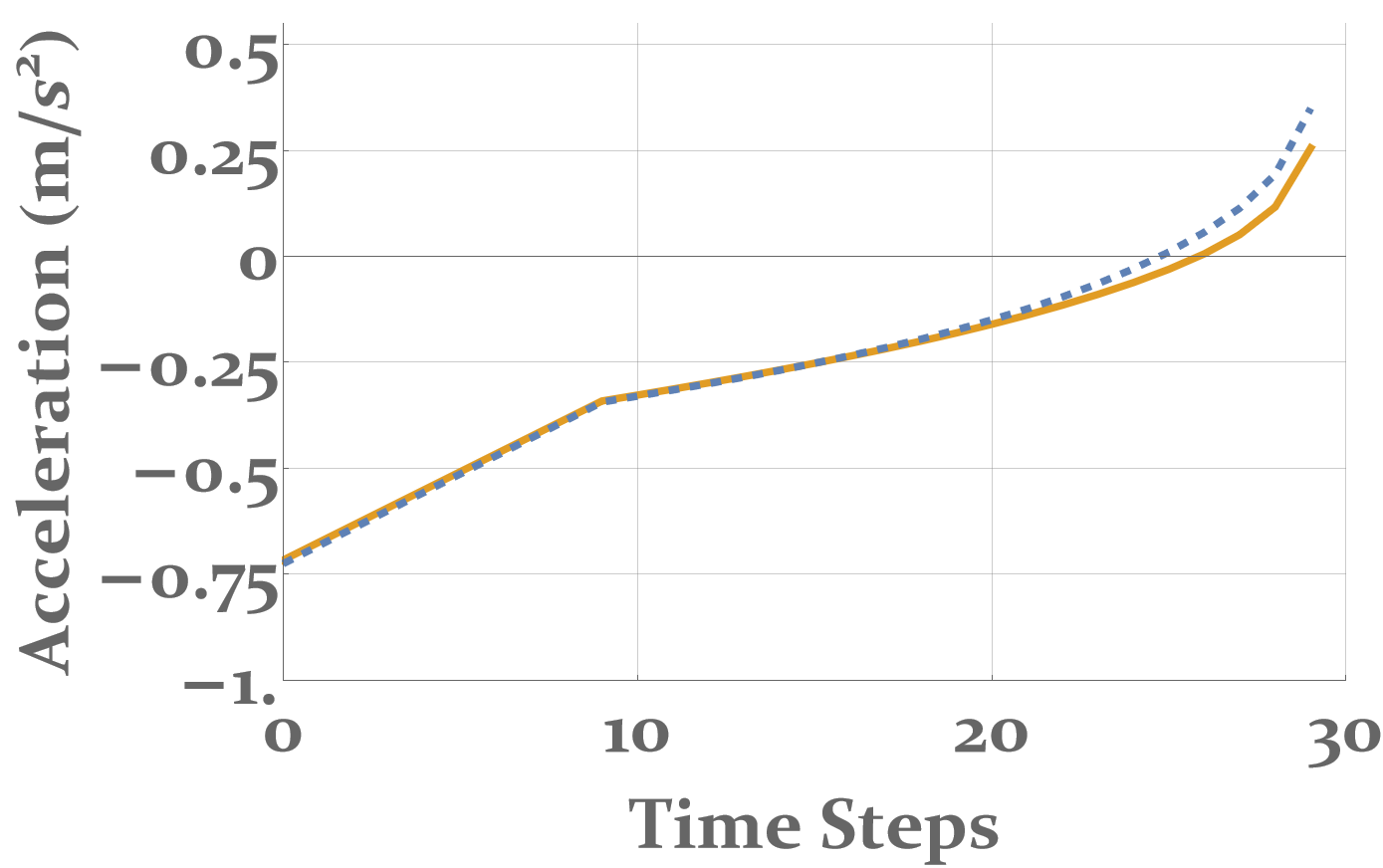}
		\caption*{Iteration 2}
	\end{subfigure}%
	\begin{subfigure}{0.25\textwidth}
		\includegraphics[width=\textwidth]{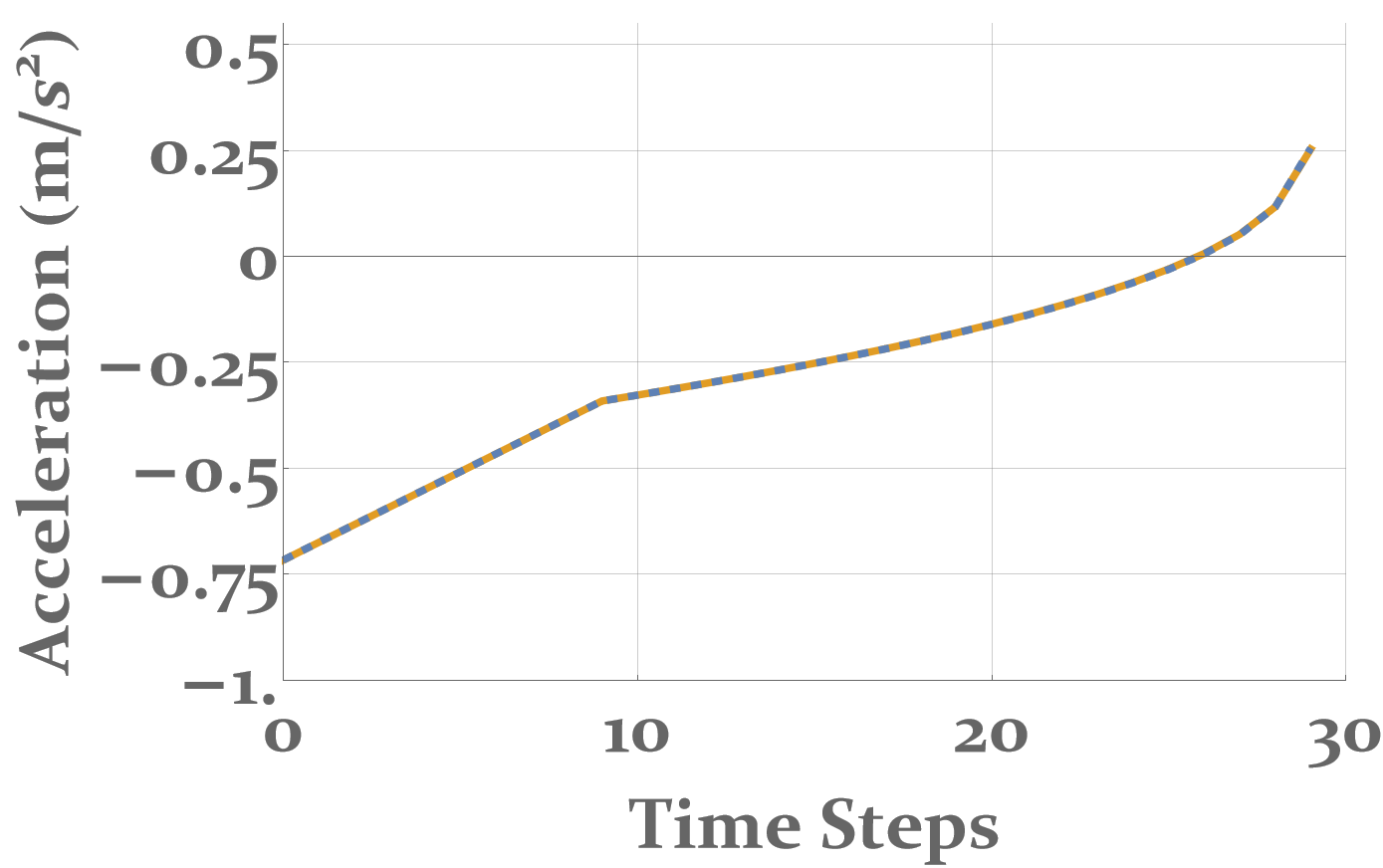}
		\caption*{Iteration 3}
	\end{subfigure}%
	\begin{subfigure}{0.25\textwidth}
		\includegraphics[width=\textwidth]{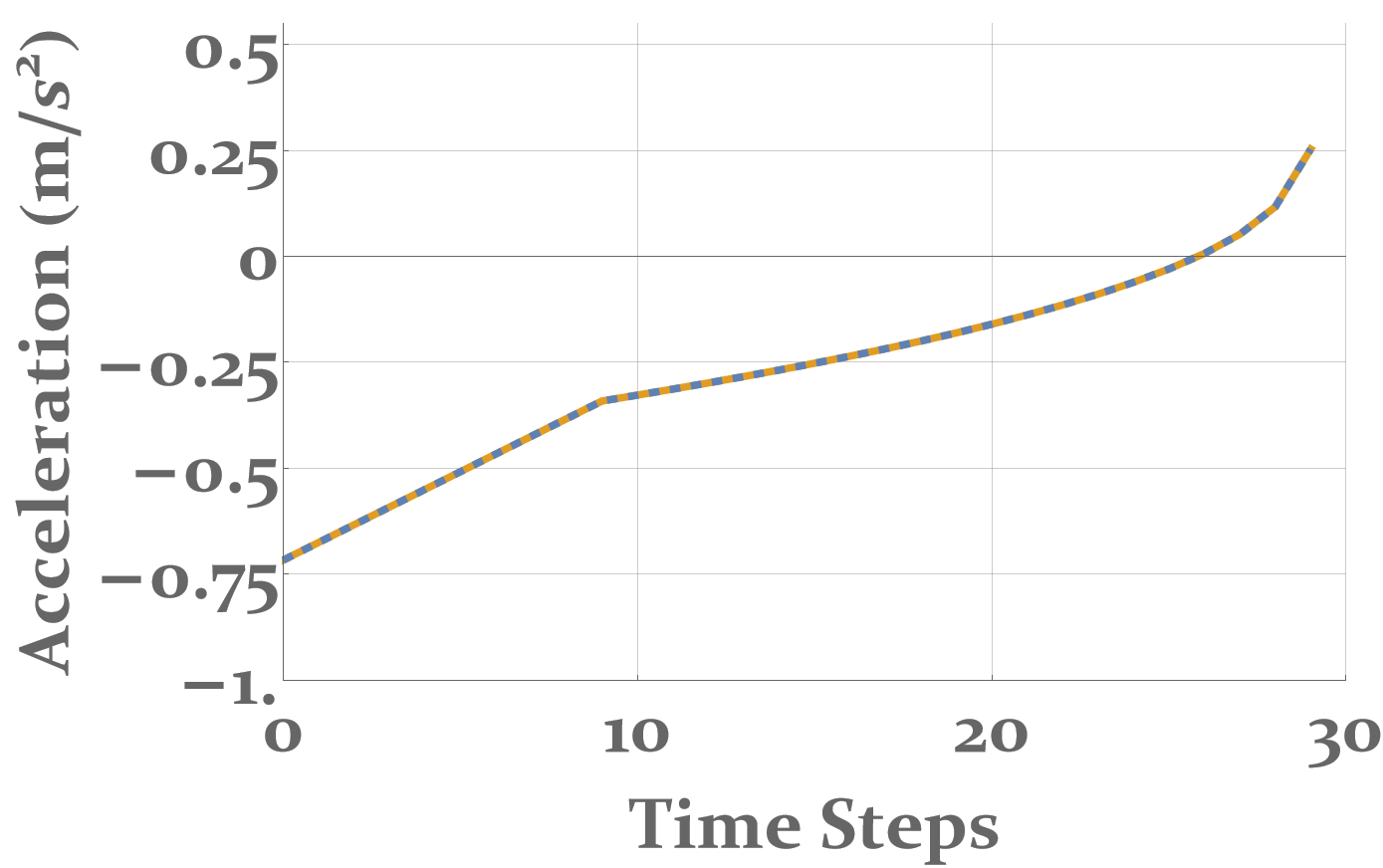}
		\caption*{Iteration 4}
	\end{subfigure}
\subcaption{ }
\end{minipage}
\hspace{0.5cm}
\begin{minipage}[b]{\textwidth}
	\begin{subfigure}{0.25\textwidth} 
		\includegraphics[width=\textwidth]{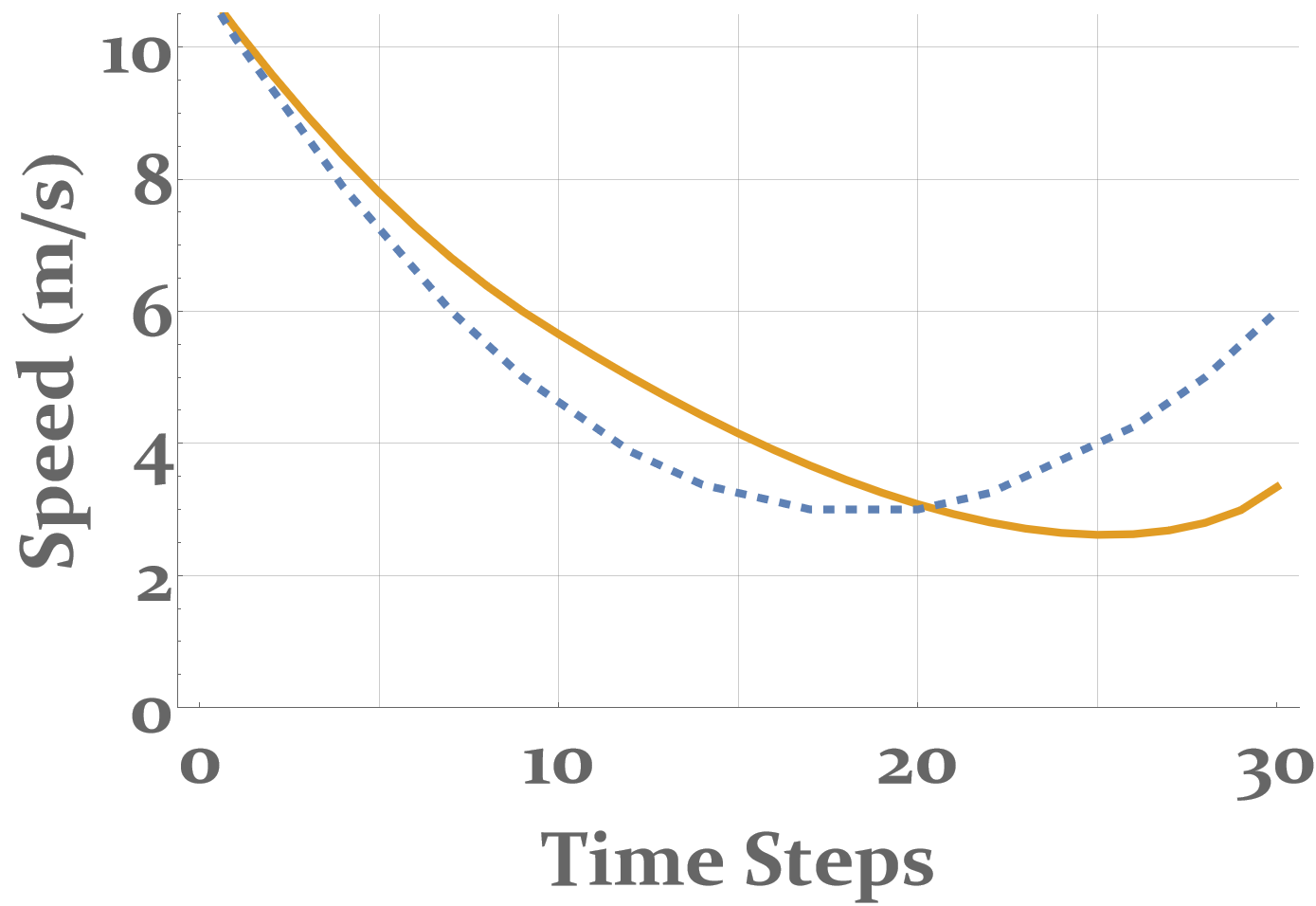}
		\caption*{Iteration 1}
	\end{subfigure}%
	\begin{subfigure}{0.25\textwidth}
		\includegraphics[width=\textwidth]{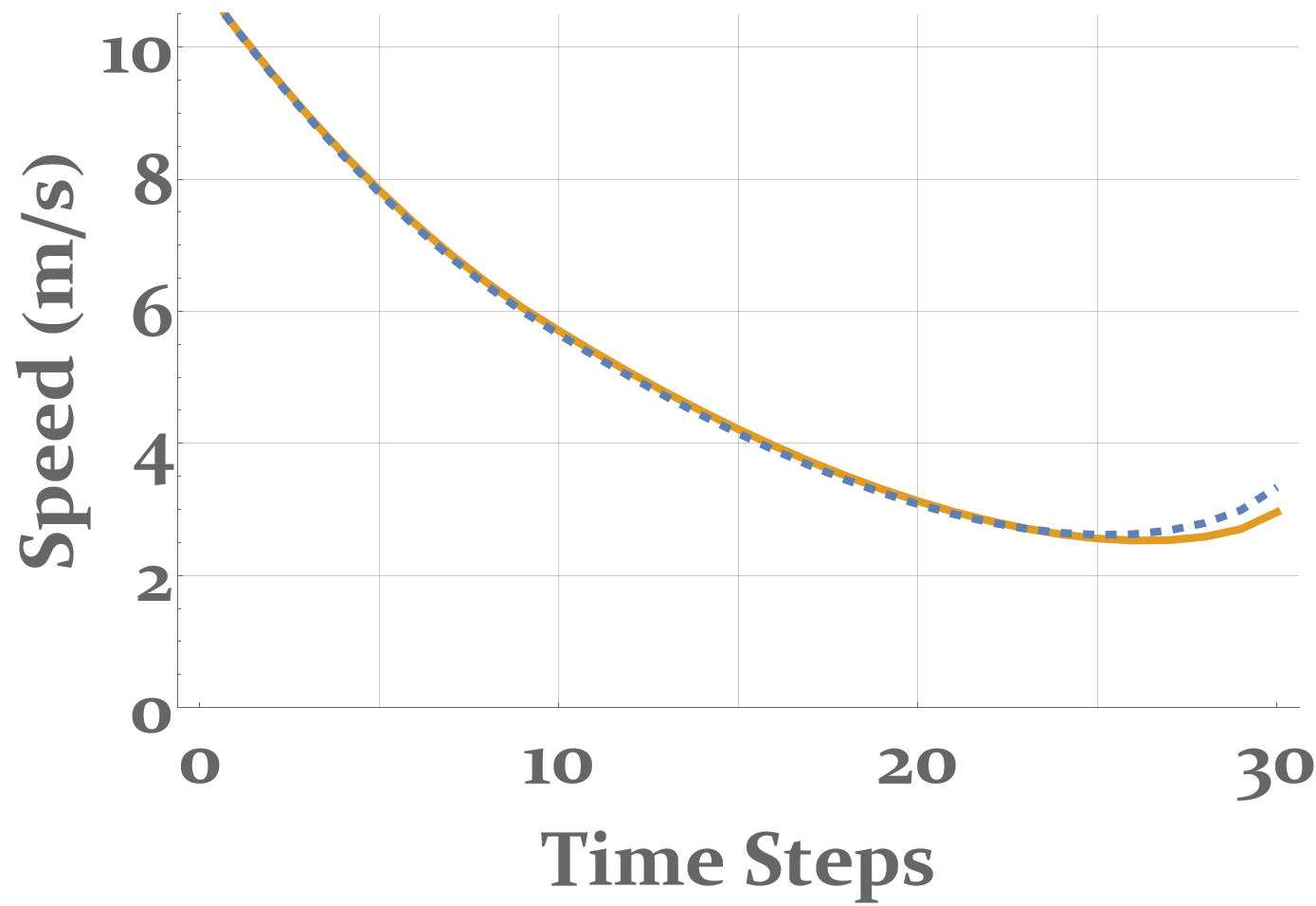}
		\caption*{Iteration 2}
	\end{subfigure}%
	\begin{subfigure}{0.25\textwidth}
		\includegraphics[width=\textwidth]{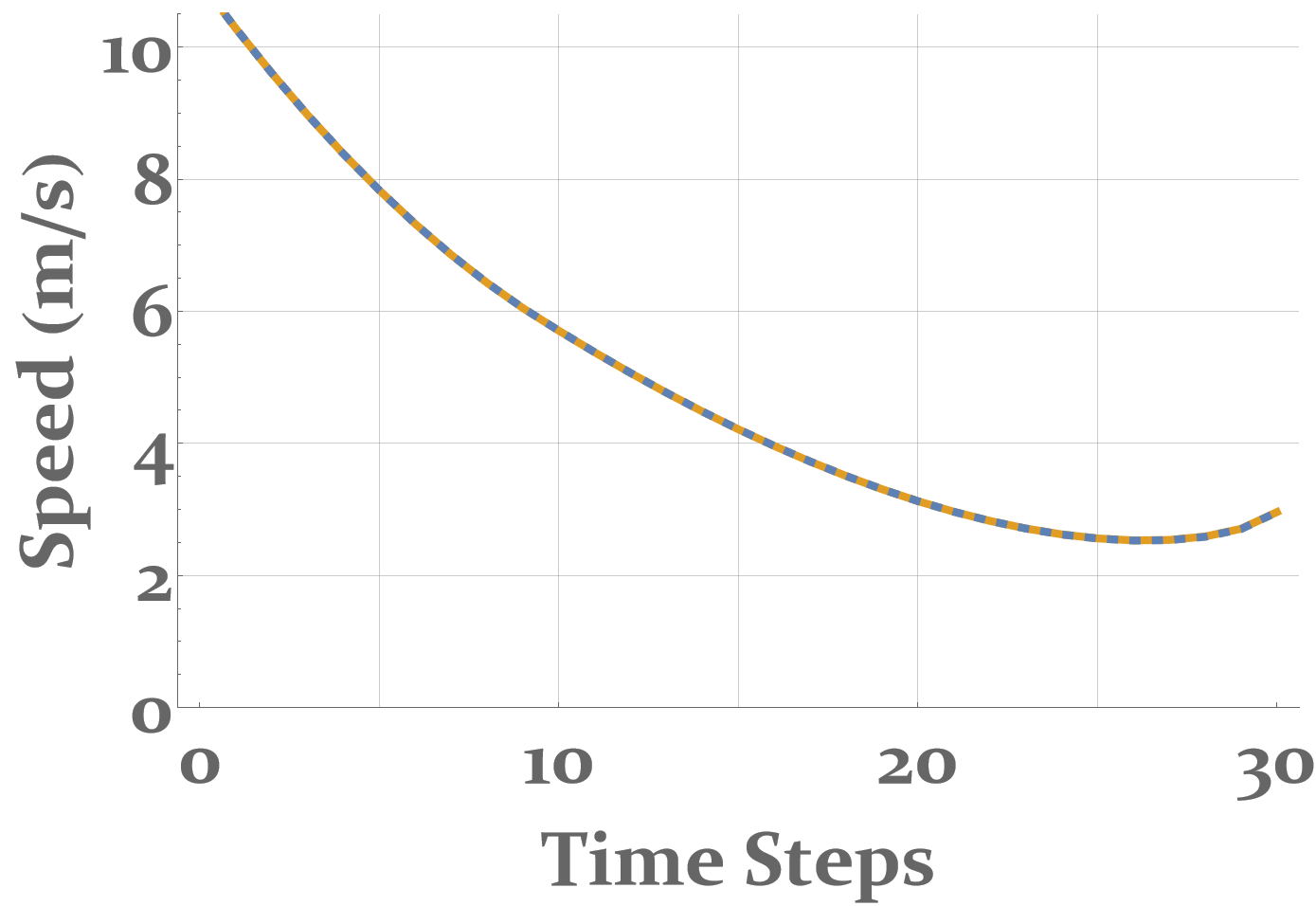}
		\caption*{Iteration 3}
	\end{subfigure}%
	\begin{subfigure}{0.25\textwidth}
		\includegraphics[width=\textwidth]{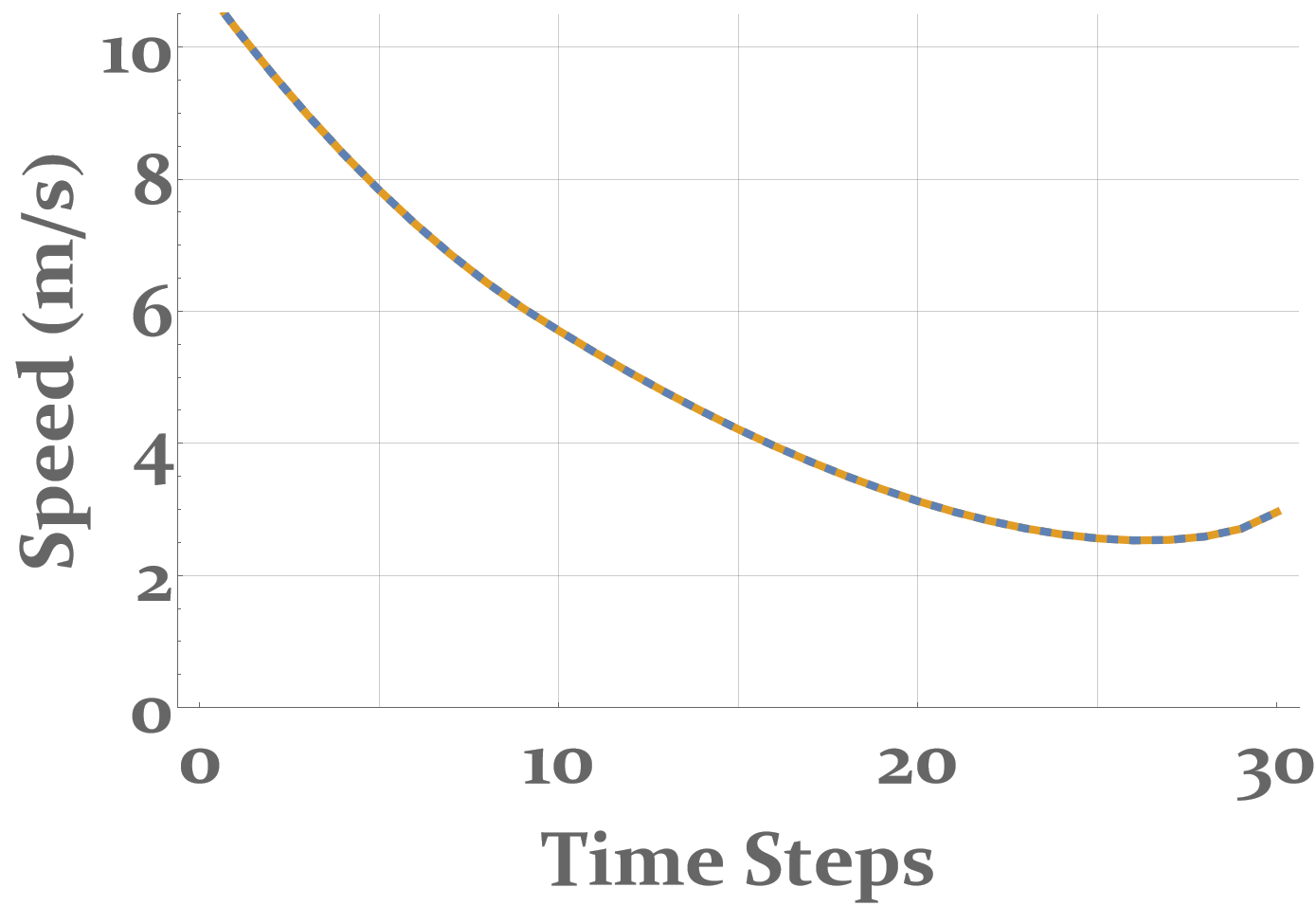}
		\caption*{Iteration 4}
	\end{subfigure}
\subcaption{ }
\end{minipage}
\caption{Received (blue dashed line) and computed (orange line) trajectories of (a) acceleration and (b) speed at each iteration of DDP algorithm for Scenario 2.}
\label{fig:ddp-s-2}
\end{figure}


\begin{figure}[h!t]
\centering
\begin{minipage}[b]{\textwidth}
	\begin{subfigure}{0.25\textwidth} 
		\includegraphics[width=\textwidth]{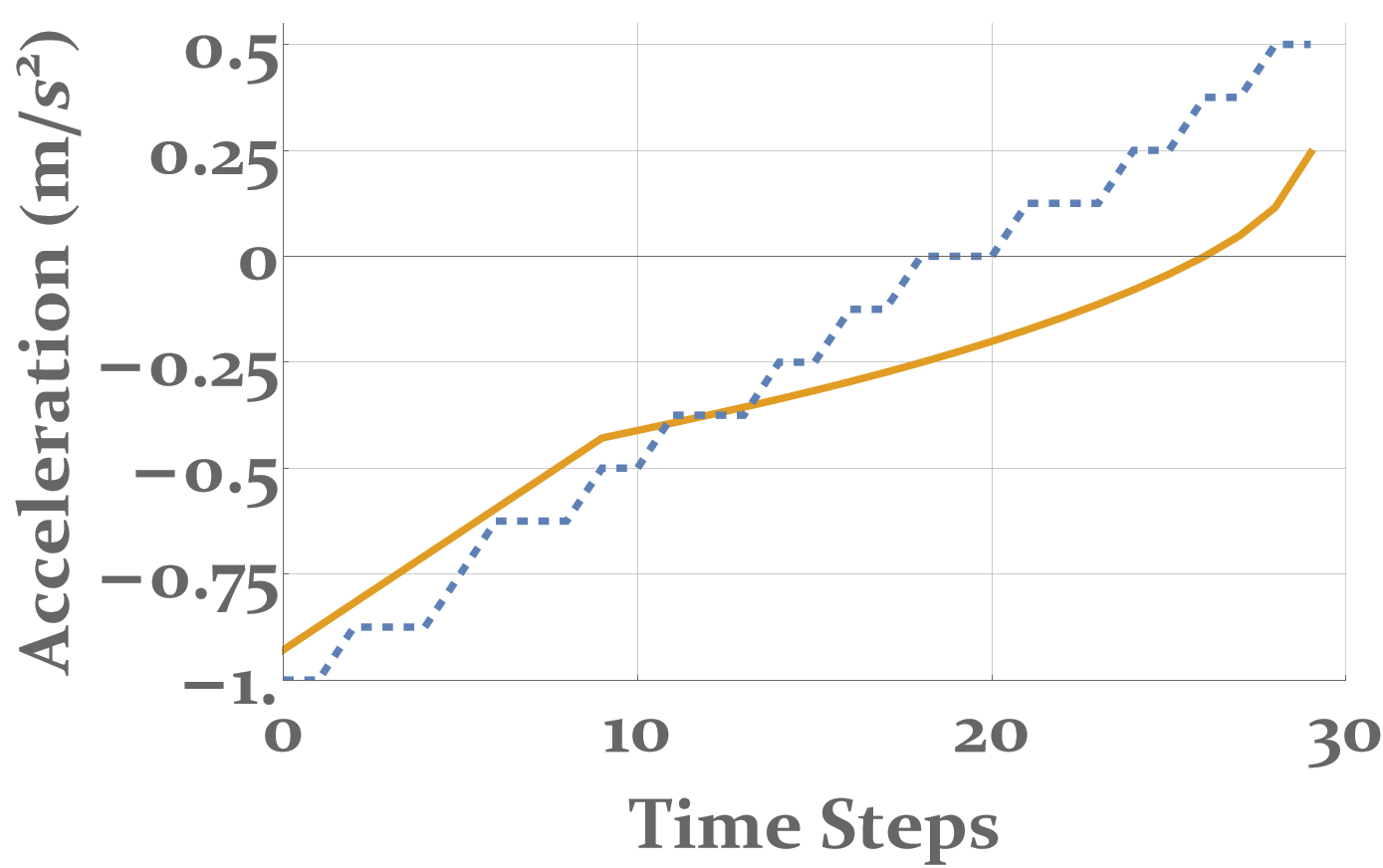}
		\caption*{Iteration 1}
	\end{subfigure}%
	\begin{subfigure}{0.25\textwidth} 
		\includegraphics[width=\textwidth]{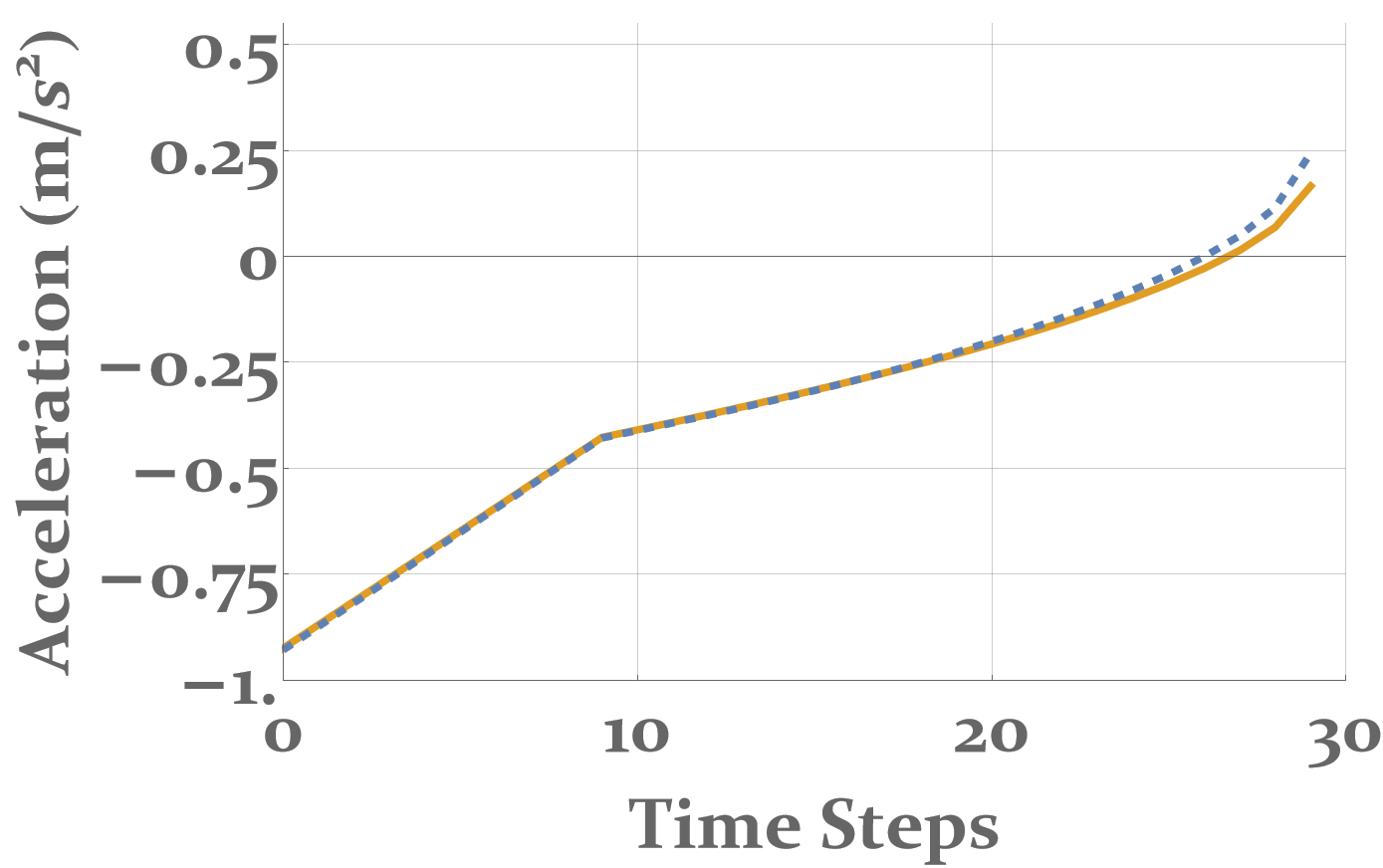}
		\caption*{Iteration 2}
	\end{subfigure}%
	\begin{subfigure}{0.25\textwidth}
		\includegraphics[width=\textwidth]{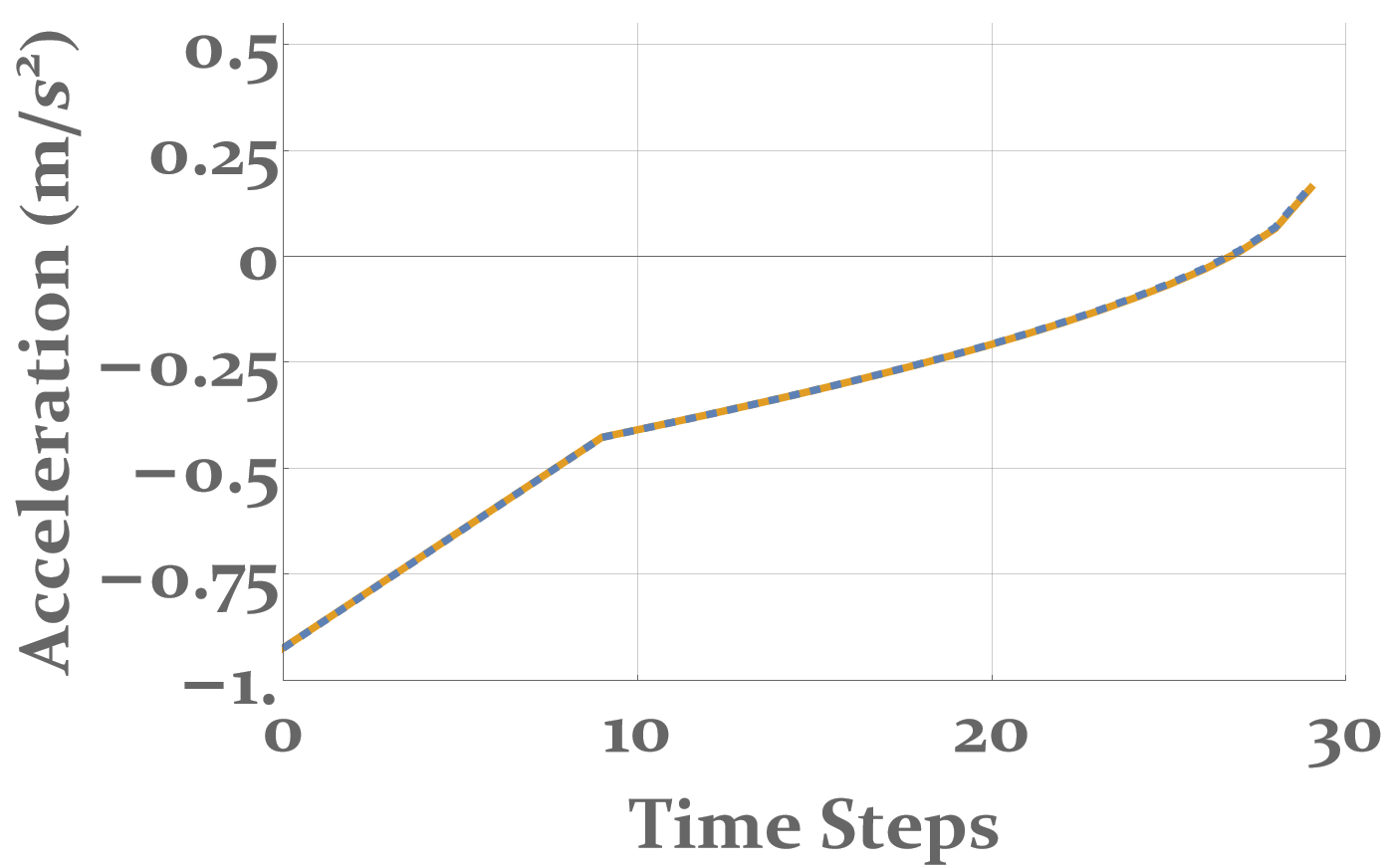}
		\caption*{Iteration 3}
	\end{subfigure}%
	\begin{subfigure}{0.25\textwidth}
		\includegraphics[width=\textwidth]{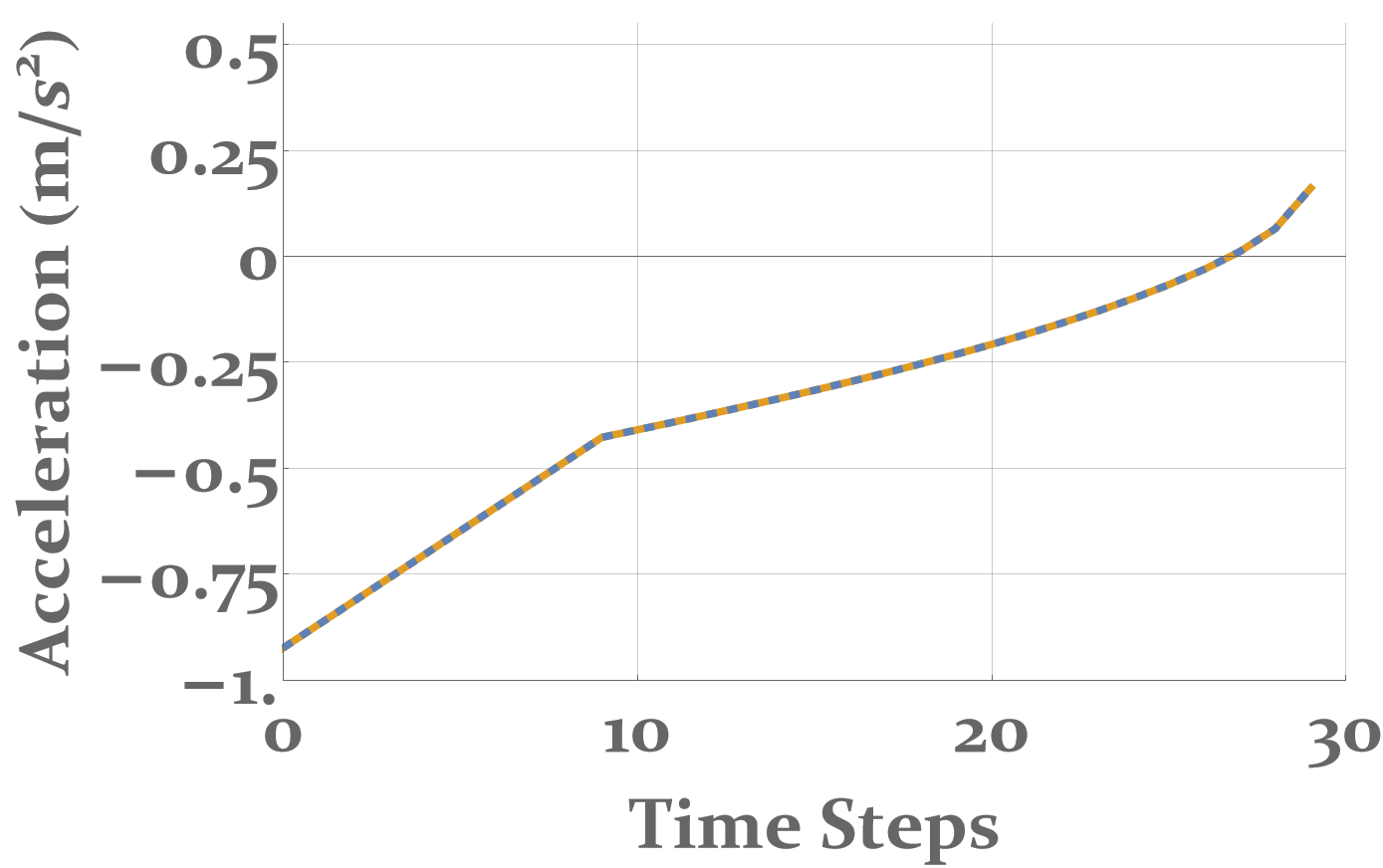}
		\caption*{Iteration 4}
	\end{subfigure}
\subcaption{ }
\end{minipage}
\hspace{0.5cm}
\begin{minipage}[b]{\textwidth}
	\begin{subfigure}{0.25\textwidth} 
		\includegraphics[width=\textwidth]{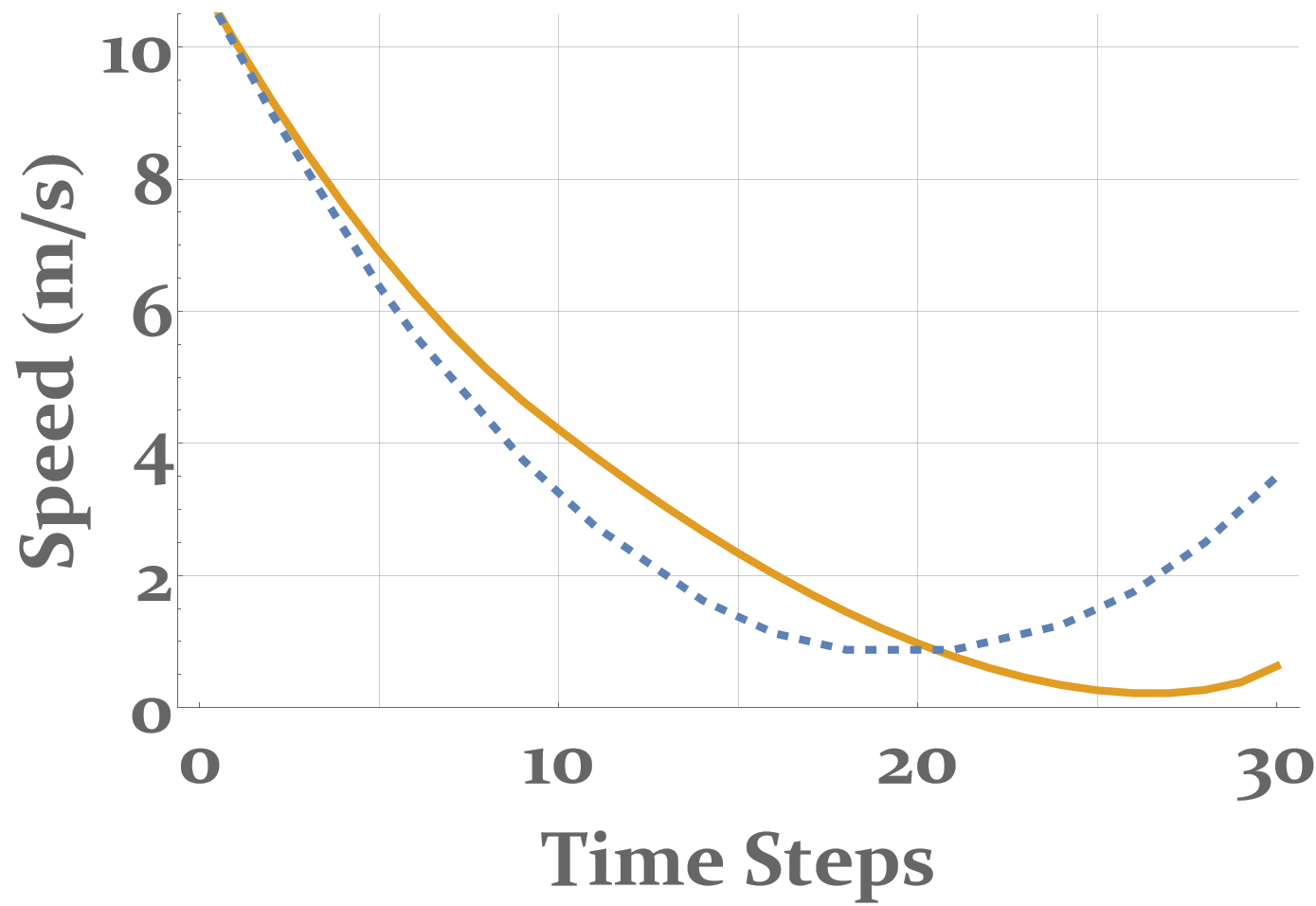}
		\caption*{Iteration 1}
	\end{subfigure}%
	\begin{subfigure}{0.25\textwidth}
		\includegraphics[width=\textwidth]{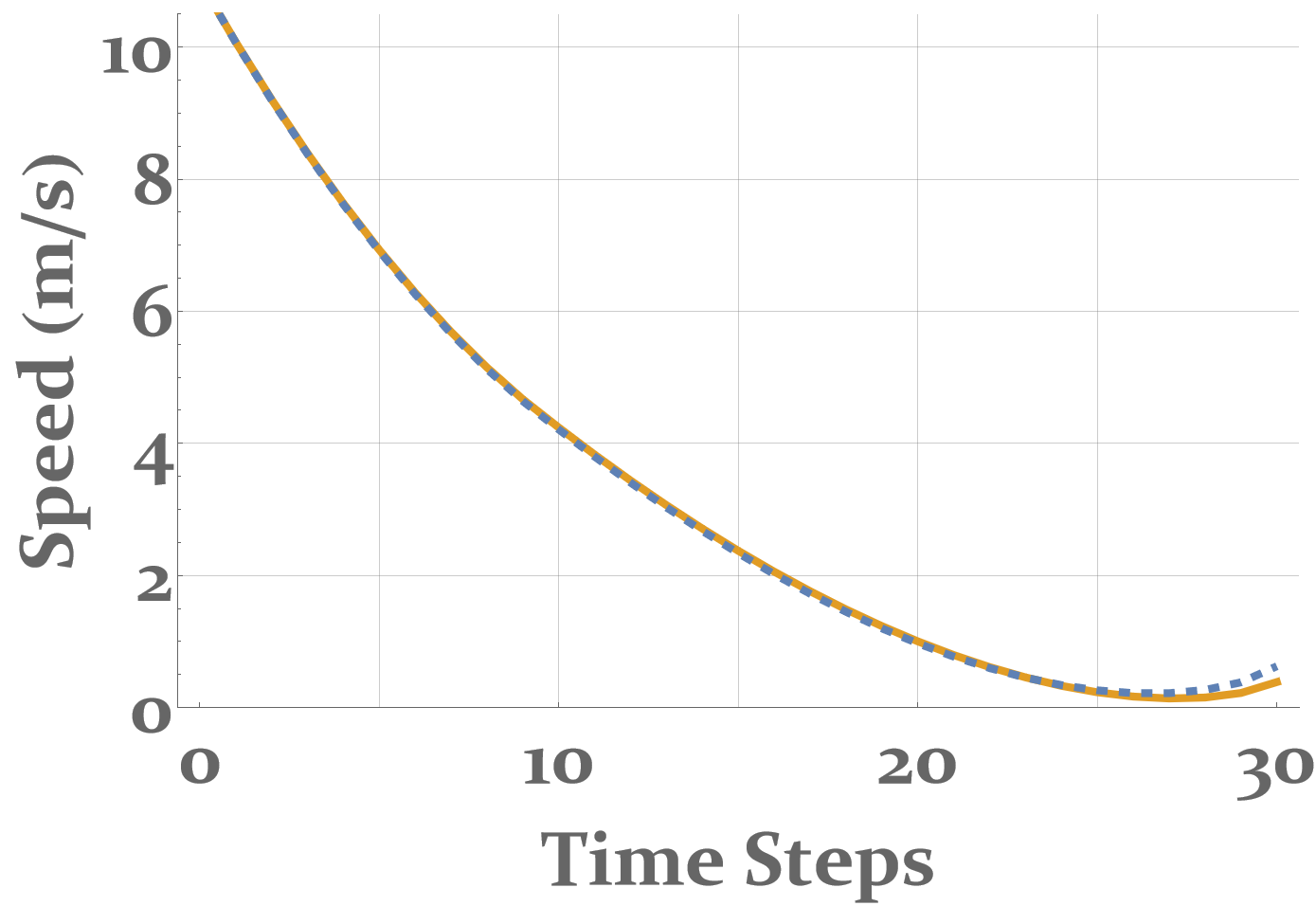}
		\caption*{Iteration 2}
	\end{subfigure}%
	\begin{subfigure}{0.25\textwidth}
		\includegraphics[width=\textwidth]{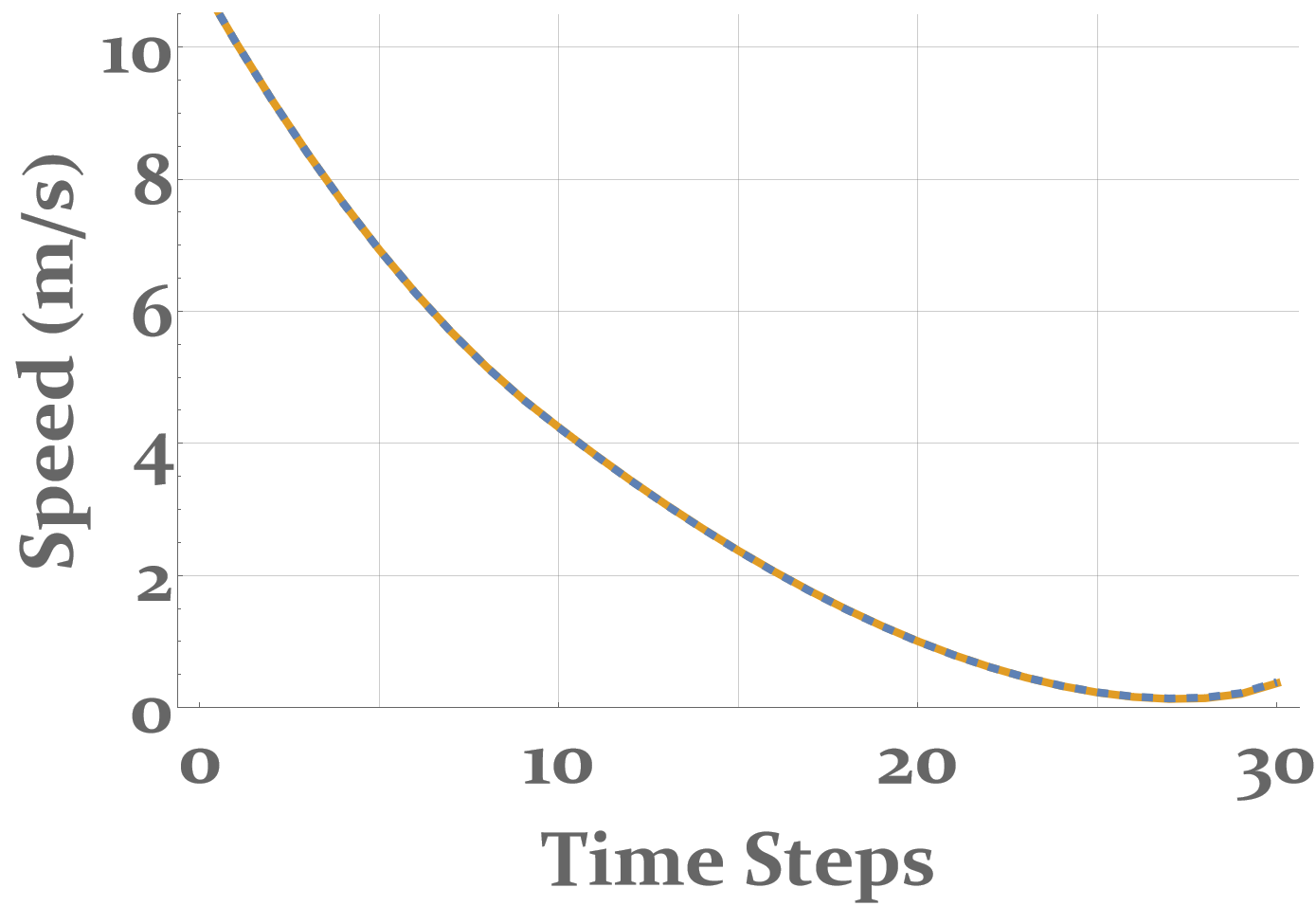}
		\caption*{Iteration 3}
	\end{subfigure}%
	\begin{subfigure}{0.25\textwidth}
		\includegraphics[width=\textwidth]{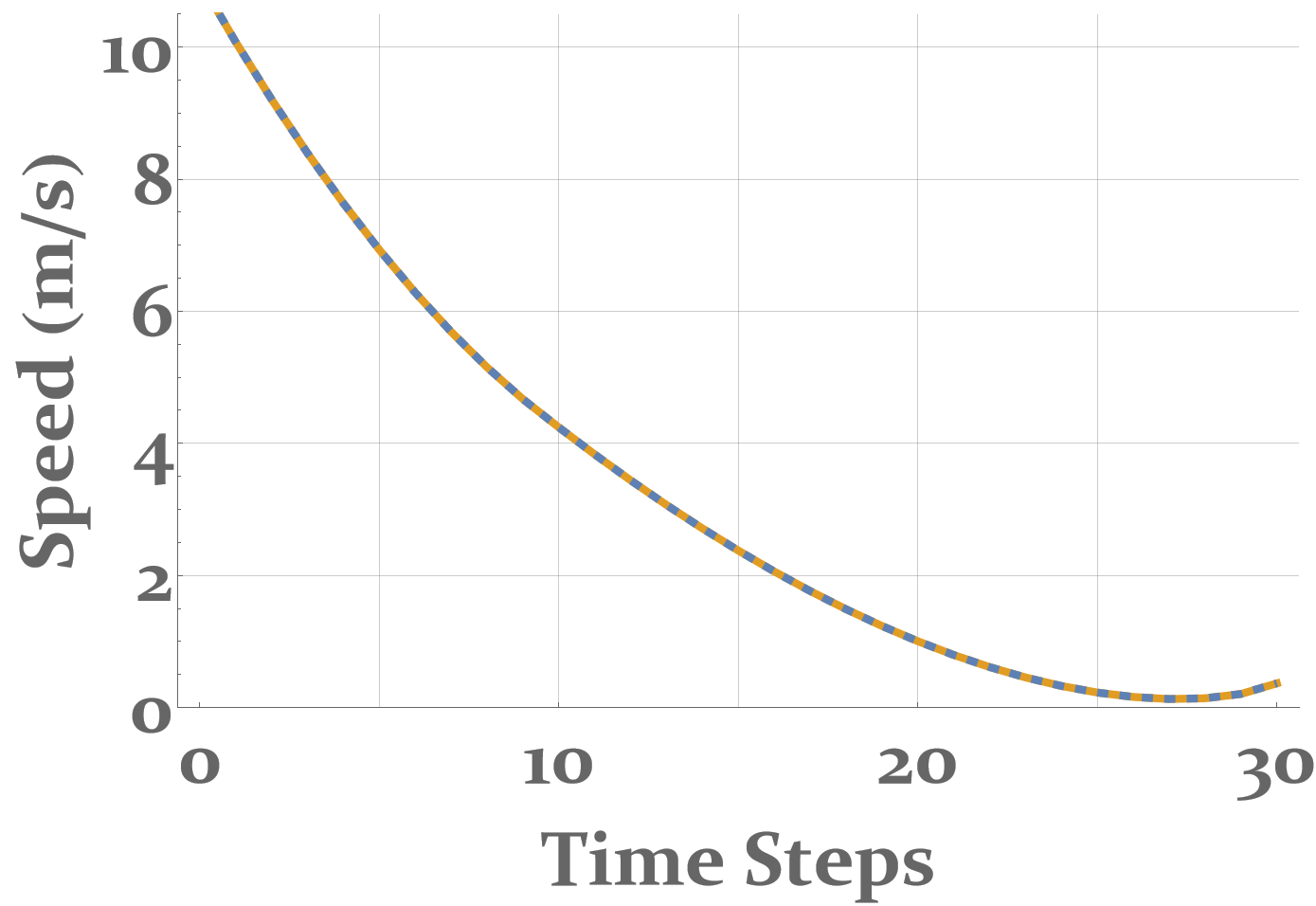}
		\caption*{Iteration 4}
	\end{subfigure}
\subcaption{ }
\end{minipage}
\caption{Received (blue dashed line) and computed (orange line) trajectories of (a) acceleration and (b) speed at each iteration of DDP algorithm for Scenario 3.}
\label{fig:ddp-s-3}
\end{figure}
\vspace*{\fill}
\newpage

\begin{table}[H]
\caption{Performance of DDDP algorithm in terms of CPU-time and optimal cost for different initial values of $\bm{C}$ and $\Delta a$.}
\label{tb:dddpC-Da}
\begin{minipage}[h]{\textwidth}
\begin{center}\small
\begin{tabular}{cccc|ccc|ccc|ccc}
\hline
 & \multicolumn{12}{c}{Scenario 1} \\\hline
 & \multicolumn{3}{c}{$\Delta a = 1.0$} & \multicolumn{3}{c}{$\Delta a = 0.5$} & \multicolumn{3}{c}{$\Delta a = 0.25$} & \multicolumn{3}{c}{$\Delta a = 0.125$}  \\\hline
$C_v$ & Iter. & \begin{tabular}{@{}c@{}} CPU \\ Times \\ (s)\end{tabular} & Cost & Iter. & \begin{tabular}{@{}c@{}} CPU \\ Times \\ (s)\end{tabular} & Cost & Iter. & \begin{tabular}{@{}c@{}} CPU \\ Times \\ (s)\end{tabular} & Cost & Iter. & \begin{tabular}{@{}c@{}} CPU \\ Times \\ (s)\end{tabular} & Cost  \\\hline
2 & 19 & 0.648 			& 1.199 	& 20 & 0.612 		& 1.193 & 10 & 0.326 		& 1.193 	& 5   & 0.280 			& 1.193 \\
3 & 19 & 1.296 			& 1.177 	& 8  & 0.552 		& 1.175 & 6  & 0.469 		& 1.179 	& 5   & 0.443 			& 1.179 \\ 
4 & 15 & $\bm{1.655}$ 	& 1.175 	& 6  & $\bm{0.692}$ & 1.175 & 7  & $\bm{0.867}$ & 1.175 	& 5   & $\bm{0.762}$ 	& 1.175 \\
5 & 12 & 1.839 			& 1.175 	& 6  & 1.052 		& 1.175 & 6  & 1.111 		& 1.175 	& 4   & 0.885 			& 1.175 \\
6 & 10 & 2.298 			& 1.175 	& 6  & 1.451 		& 1.175 & 5  & 1.291 		& 1.175 	& 3   & 0.967 			& 1.175 \\
\hline
& \multicolumn{12}{c}{Scenario 2} \\\hline
& \multicolumn{3}{c}{$\Delta a = 1.0$} & \multicolumn{3}{c}{$\Delta a = 0.5$} & \multicolumn{3}{c}{$\Delta a = 0.25$} & \multicolumn{3}{c}{$\Delta a = 0.125$}  \\\hline
$C_v$ & Iter. & \begin{tabular}{@{}c@{}} CPU \\ Times \\ (s)\end{tabular} & Cost & Iter. & \begin{tabular}{@{}c@{}} CPU \\ Times \\ (s)\end{tabular} & Cost & Iter. & \begin{tabular}{@{}c@{}} CPU \\ Times \\ (s)\end{tabular} & Cost & Iter. & \begin{tabular}{@{}c@{}} CPU \\ Times \\ (s)\end{tabular} & Cost \\\hline
2 & 29 & 1.212 			& 3.909 	& 18 & 0.790 		& 3.909 	& 16 & 0.729 		& 3.909 	& 21 & 1.014 		& 3.909 \\
3 & 22 & 1.656 			& 3.906 	& 14 & 1.093 		& 3.906 	& 12 & 0.966 		& 3.906 	& 18 & 1.568 		& 3.906 \\ 
4 & 13 & $\bm{1.578}$ 	& 3.906 	& 10 & $\bm{1.246}$ & 3.906 	& 9  & $\bm{1.149}$ & 3.906 	& 13 & $\bm{1.740}$ & 3.906 \\
5 & 11 & 1.920 			& 3.906 	& 7  & 1.273 		& 3.906 	& 8  & 1.514 		& 3.906 	& 11 & 2.178 		& 3.906 \\
6 & 14 & 3.324 			& 3.906 	& 10 & 2.471 		& 3.906 	& 7  & 1.823 		& 3.906 	& 9  & 2.457 		& 3.906 \\
\hline
& \multicolumn{12}{c}{Scenario 3} \\\hline
& \multicolumn{3}{c}{$\Delta a = 1.0$} & \multicolumn{3}{c}{$\Delta a = 0.5$} & \multicolumn{3}{c}{$\Delta a = 0.25$} & \multicolumn{3}{c}{$\Delta a = 0.125$}  \\\hline
2 & 28 & 1.152 			& 6.3583 	& 20 & 0.868 		& 6.3583	& 18 & 0.815 		& 6.3583 	& 21 & 0.998 		& 6.3644 \\
3 & 21 & 1.543 			& 6.3582 	& 16 & 1.192 		& 6.3582  	& 14 & 1.098 		& 6.3582 	& 17 & 1.419 		& 6.3639 \\
4 & 14 & $\bm{1.506}$ 	& 6.3582 	& 10 & $\bm{1.128}$ & 6.3582 	& 10 & $\bm{1.179}$ & 6.3582 	& 14 & $\bm{1.834}$ & 6.3582 \\
5 & 12 & 1.777 			& 6.3582 	& 8  & 1.253 		& 6.3582 	& 8 & 1.400 		& 6.3582 	& 11 & 2.101 		& 6.3582 \\
6 & 11 & 2.050 			& 6.3582 	& 7  & 1.368 		& 6.3582 	& 7 & 1.553 		& 6.3582 	& 9  & 2.328 		& 6.3582 \\
\hline
\end{tabular}%
\end{center}
\end{minipage}
\end{table}
\vspace*{\fill}

%
\bibliography{cas-sc-template}

\begin{thebibliography}{32}
\expandafter\ifx\csname natexlab\endcsname\relax\def\natexlab#1{#1}\fi
\providecommand{\url}[1]{\texttt{#1}}
\providecommand{\href}[2]{#2}
\providecommand{\path}[1]{#1}
\providecommand{\DOIprefix}{doi:}
\providecommand{\ArXivprefix}{arXiv:}
\providecommand{\URLprefix}{URL: }
\providecommand{\Pubmedprefix}{pmid:}
\providecommand{\doi}[1]{\href{http://dx.doi.org/#1}{\path{#1}}}
\providecommand{\Pubmed}[1]{\href{pmid:#1}{\path{#1}}}
\providecommand{\bibinfo}[2]{#2}
\ifx\xfnm\relax \def\xfnm[#1]{\unskip,\space#1}\fi
\bibitem[{Bertsekas(2012)}]{bertsekas2012dynamic}
\bibinfo{author}{Bertsekas, D.}, \bibinfo{year}{2012}.
\newblock \bibinfo{title}{Dynamic programming and optimal control: Volume I}.
  volume~\bibinfo{volume}{1}.
\newblock \bibinfo{publisher}{Athena scientific}.
\bibitem[{Blokpoel and Lu(2017)}]{blokpoel2017signal}
\bibinfo{author}{Blokpoel, R.}, \bibinfo{author}{Lu, M.}, \bibinfo{year}{2017}.
\newblock \bibinfo{title}{Signal plan stabilization to enable eco-driving}, in:
  \bibinfo{booktitle}{2017 IEEE 20th International Conference on Intelligent
  Transportation Systems (ITSC)}, \bibinfo{organization}{IEEE}. pp.
  \bibinfo{pages}{133--138}.
\bibitem[{Budhiraja et~al.(2018)Budhiraja, Carpentier, Mastalli and
  Mansard}]{budhiraja2018differential}
\bibinfo{author}{Budhiraja, R.}, \bibinfo{author}{Carpentier, J.},
  \bibinfo{author}{Mastalli, C.}, \bibinfo{author}{Mansard, N.},
  \bibinfo{year}{2018}.
\newblock \bibinfo{title}{Differential dynamic programming for multi-phase
  rigid contact dynamics}, in: \bibinfo{booktitle}{2018 IEEE-RAS 18th
  International Conference on Humanoid Robots (Humanoids)},
  \bibinfo{organization}{IEEE}. pp. \bibinfo{pages}{1--9}.
\bibitem[{Cleveland et~al.(2017)Cleveland, Grosse and
  Shyu}]{cleveland2017local}
\bibinfo{author}{Cleveland, W.S.}, \bibinfo{author}{Grosse, E.},
  \bibinfo{author}{Shyu, W.M.}, \bibinfo{year}{2017}.
\newblock \bibinfo{title}{Local regression models}, in:
  \bibinfo{booktitle}{Statistical models in S}. \bibinfo{publisher}{Routledge},
  pp. \bibinfo{pages}{309--376}.
\bibitem[{Fletcher(1971)}]{fletcher1971general}
\bibinfo{author}{Fletcher, R.}, \bibinfo{year}{1971}.
\newblock \bibinfo{title}{A general quadratic programming algorithm}.
\newblock \bibinfo{journal}{IMA Journal of Applied Mathematics}
  \bibinfo{volume}{7}, \bibinfo{pages}{76--91}.
\bibitem[{Heidari et~al.(1971)Heidari, Chow, Kokotovi{\'c} and
  Meredith}]{heidari1971discrete}
\bibinfo{author}{Heidari, M.}, \bibinfo{author}{Chow, V.T.},
  \bibinfo{author}{Kokotovi{\'c}, P.V.}, \bibinfo{author}{Meredith, D.D.},
  \bibinfo{year}{1971}.
\newblock \bibinfo{title}{Discrete differential dynamic programing approach to
  water resources systems optimization}.
\newblock \bibinfo{journal}{Water Resources Research} \bibinfo{volume}{7},
  \bibinfo{pages}{273--282}.
\bibitem[{Hounsell and McDonald(2001)}]{hounsell2001urban}
\bibinfo{author}{Hounsell, N.}, \bibinfo{author}{McDonald, M.},
  \bibinfo{year}{2001}.
\newblock \bibinfo{title}{Urban network traffic control}.
\newblock \bibinfo{journal}{Proceedings of the Institution of Mechanical
  Engineers, Part I: Journal of Systems and Control Engineering}
  \bibinfo{volume}{215}, \bibinfo{pages}{325--334}.
\bibitem[{Jacobson and Mayne(1970)}]{jacobson1970differential}
\bibinfo{author}{Jacobson, D.H.}, \bibinfo{author}{Mayne, D.Q.},
  \bibinfo{year}{1970}.
\newblock \bibinfo{title}{Differential dynamic programming}.
\newblock \bibinfo{number}{24}, \bibinfo{publisher}{Elsevier Publishing
  Company}.
\bibitem[{Jamshidnejad et~al.(2017)Jamshidnejad, Papamichail, Papageorgiou and
  De~Schutter}]{jamshidnejad2017sustainable}
\bibinfo{author}{Jamshidnejad, A.}, \bibinfo{author}{Papamichail, I.},
  \bibinfo{author}{Papageorgiou, M.}, \bibinfo{author}{De~Schutter, B.},
  \bibinfo{year}{2017}.
\newblock \bibinfo{title}{Sustainable model-predictive control in urban traffic
  networks: Efficient solution based on general smoothening methods}.
\newblock \bibinfo{journal}{IEEE Transactions on Control Systems Technology}
  \bibinfo{volume}{26}, \bibinfo{pages}{813--827}.
\bibitem[{Katsaros et~al.(2011)Katsaros, Kernchen, Dianati and
  Rieck}]{katsaros2011performance}
\bibinfo{author}{Katsaros, K.}, \bibinfo{author}{Kernchen, R.},
  \bibinfo{author}{Dianati, M.}, \bibinfo{author}{Rieck, D.},
  \bibinfo{year}{2011}.
\newblock \bibinfo{title}{Performance study of a green light optimized speed
  advisory (glosa) application using an integrated cooperative its simulation
  platform}, in: \bibinfo{booktitle}{2011 7th International Wireless
  Communications and Mobile Computing Conference},
  \bibinfo{organization}{IEEE}. pp. \bibinfo{pages}{918--923}.
\bibitem[{Kosmatopoulos et~al.(2006)Kosmatopoulos, Papageorgiou, Bielefeldt,
  Dinopoulou, Morris, Mueck, Richards and
  Weichenmeier}]{kosmatopoulos2006international}
\bibinfo{author}{Kosmatopoulos, E.}, \bibinfo{author}{Papageorgiou, M.},
  \bibinfo{author}{Bielefeldt, C.}, \bibinfo{author}{Dinopoulou, V.},
  \bibinfo{author}{Morris, R.}, \bibinfo{author}{Mueck, J.},
  \bibinfo{author}{Richards, A.}, \bibinfo{author}{Weichenmeier, F.},
  \bibinfo{year}{2006}.
\newblock \bibinfo{title}{International comparative field evaluation of a
  traffic-responsive signal control strategy in three cities}.
\newblock \bibinfo{journal}{Transportation Research Part A: Policy and
  Practice} \bibinfo{volume}{40}, \bibinfo{pages}{399--413}.
\bibitem[{Koukoumidis et~al.(2011)Koukoumidis, Peh and
  Martonosi}]{koukoumidis2011signalguru}
\bibinfo{author}{Koukoumidis, E.}, \bibinfo{author}{Peh, L.S.},
  \bibinfo{author}{Martonosi, M.R.}, \bibinfo{year}{2011}.
\newblock \bibinfo{title}{Signalguru: leveraging mobile phones for
  collaborative traffic signal schedule advisory}, in:
  \bibinfo{booktitle}{Proceedings of the 9th international conference on Mobile
  systems, applications, and services}, pp. \bibinfo{pages}{127--140}.
\bibitem[{Lawitzky et~al.(2013)Lawitzky, Wollherr and
  Buss}]{lawitzky2013energy}
\bibinfo{author}{Lawitzky, A.}, \bibinfo{author}{Wollherr, D.},
  \bibinfo{author}{Buss, M.}, \bibinfo{year}{2013}.
\newblock \bibinfo{title}{Energy optimal control to approach traffic lights},
  in: \bibinfo{booktitle}{2013 IEEE/RSJ International Conference on Intelligent
  Robots and Systems}, \bibinfo{organization}{IEEE}. pp.
  \bibinfo{pages}{4382--4387}.
\bibitem[{Lu et~al.(2020)Lu, Ferragut, Kutila and Chen}]{lu2020next}
\bibinfo{author}{Lu, M.}, \bibinfo{author}{Ferragut, J.},
  \bibinfo{author}{Kutila, M.}, \bibinfo{author}{Chen, T.},
  \bibinfo{year}{2020}.
\newblock \bibinfo{title}{Next-generation wireless networks for v2x}, in:
  \bibinfo{booktitle}{2020 IEEE 23rd International Conference on Intelligent
  Transportation Systems (ITSC)}, \bibinfo{organization}{IEEE}. pp.
  \bibinfo{pages}{1--5}.
\bibitem[{Ma et~al.(2018)Ma, Zhou, Huang, Melson, James and
  Zhang}]{ma2018hardware}
\bibinfo{author}{Ma, J.}, \bibinfo{author}{Zhou, F.}, \bibinfo{author}{Huang,
  Z.}, \bibinfo{author}{Melson, C.L.}, \bibinfo{author}{James, R.},
  \bibinfo{author}{Zhang, X.}, \bibinfo{year}{2018}.
\newblock \bibinfo{title}{Hardware-in-the-loop testing of connected and
  automated vehicle applications: a use case for queue-aware signalized
  intersection approach and departure}.
\newblock \bibinfo{journal}{Transportation Research Record}
  \bibinfo{volume}{2672}, \bibinfo{pages}{36--46}.
\bibitem[{Mahler and Vahidi(2012)}]{mahler2012reducing}
\bibinfo{author}{Mahler, G.}, \bibinfo{author}{Vahidi, A.},
  \bibinfo{year}{2012}.
\newblock \bibinfo{title}{Reducing idling at red lights based on probabilistic
  prediction of traffic signal timings}, in: \bibinfo{booktitle}{2012 American
  Control Conference (ACC)}, \bibinfo{organization}{IEEE}. pp.
  \bibinfo{pages}{6557--6562}.
\bibitem[{Murray and Yakowitz(1984)}]{murray1984differential}
\bibinfo{author}{Murray, D.}, \bibinfo{author}{Yakowitz, S.},
  \bibinfo{year}{1984}.
\newblock \bibinfo{title}{Differential dynamic programming and newton's method
  for discrete optimal control problems}.
\newblock \bibinfo{journal}{Journal of Optimization Theory and Applications}
  \bibinfo{volume}{43}, \bibinfo{pages}{395--414}.
\bibitem[{Murray and Yakowitz(1979)}]{murray1979constrained}
\bibinfo{author}{Murray, D.M.}, \bibinfo{author}{Yakowitz, S.J.},
  \bibinfo{year}{1979}.
\newblock \bibinfo{title}{Constrained differential dynamic programming and its
  application to multireservoir control}.
\newblock \bibinfo{journal}{Water Resources Research} \bibinfo{volume}{15},
  \bibinfo{pages}{1017--1027}.
\bibitem[{Papageorgiou et~al.(2003)Papageorgiou, Diakaki, Dinopoulou, Kotsialos
  and Wang}]{papageorgiou2003review}
\bibinfo{author}{Papageorgiou, M.}, \bibinfo{author}{Diakaki, C.},
  \bibinfo{author}{Dinopoulou, V.}, \bibinfo{author}{Kotsialos, A.},
  \bibinfo{author}{Wang, Y.}, \bibinfo{year}{2003}.
\newblock \bibinfo{title}{Review of road traffic control strategies}.
\newblock \bibinfo{journal}{Proceedings of the IEEE} \bibinfo{volume}{91},
  \bibinfo{pages}{2043--2067}.
\bibitem[{Papageorgiou et~al.(2015)Papageorgiou, Leibold and
  Buss}]{papageorgiou2015optimierung}
\bibinfo{author}{Papageorgiou, M.}, \bibinfo{author}{Leibold, M.},
  \bibinfo{author}{Buss, M.}, \bibinfo{year}{2015}.
\newblock \bibinfo{title}{Optimierung}. volume~\bibinfo{volume}{4}.
\newblock \bibinfo{publisher}{Springer}.
\bibitem[{Sanchez et~al.(2006)Sanchez, Cano and Kim}]{sanchez2006predicting}
\bibinfo{author}{Sanchez, M.}, \bibinfo{author}{Cano, J.C.},
  \bibinfo{author}{Kim, D.}, \bibinfo{year}{2006}.
\newblock \bibinfo{title}{Predicting traffic lights to improve urban traffic
  fuel consumption}, in: \bibinfo{booktitle}{2006 6th International Conference
  on ITS Telecommunications}, \bibinfo{organization}{IEEE}. pp.
  \bibinfo{pages}{331--336}.
\bibitem[{Stahlmann et~al.(2016)Stahlmann, M{\"o}ller, Brauer, German and
  Eckhoff}]{stahlmann2016technical}
\bibinfo{author}{Stahlmann, R.}, \bibinfo{author}{M{\"o}ller, M.},
  \bibinfo{author}{Brauer, A.}, \bibinfo{author}{German, R.},
  \bibinfo{author}{Eckhoff, D.}, \bibinfo{year}{2016}.
\newblock \bibinfo{title}{Technical evaluation of glosa systems and results
  from the field}, in: \bibinfo{booktitle}{2016 IEEE Vehicular Networking
  Conference (VNC)}, \bibinfo{organization}{IEEE}. pp. \bibinfo{pages}{1--8}.
\bibitem[{Sun et~al.(2020)Sun, Guanetti, Borrelli and Moura}]{sun2020optimal}
\bibinfo{author}{Sun, C.}, \bibinfo{author}{Guanetti, J.},
  \bibinfo{author}{Borrelli, F.}, \bibinfo{author}{Moura, S.J.},
  \bibinfo{year}{2020}.
\newblock \bibinfo{title}{Optimal eco-driving control of connected and
  autonomous vehicles through signalized intersections}.
\newblock \bibinfo{journal}{IEEE Internet of Things Journal}
  \bibinfo{volume}{7}, \bibinfo{pages}{3759--3773}.
\bibitem[{Tassa et~al.(2014)Tassa, Mansard and Todorov}]{tassa2014control}
\bibinfo{author}{Tassa, Y.}, \bibinfo{author}{Mansard, N.},
  \bibinfo{author}{Todorov, E.}, \bibinfo{year}{2014}.
\newblock \bibinfo{title}{Control-limited differential dynamic programming},
  in: \bibinfo{booktitle}{2014 IEEE International Conference on Robotics and
  Automation (ICRA)}, \bibinfo{organization}{IEEE}. pp.
  \bibinfo{pages}{1168--1175}.
\bibitem[{Typaldos et~al.(2020a)Typaldos, Kalogianni, Mountakis, Papamichail
  and Papageorgiou}]{typaldos2020vehicle}
\bibinfo{author}{Typaldos, P.}, \bibinfo{author}{Kalogianni, I.},
  \bibinfo{author}{Mountakis, K.S.}, \bibinfo{author}{Papamichail, I.},
  \bibinfo{author}{Papageorgiou, M.}, \bibinfo{year}{2020}a.
\newblock \bibinfo{title}{Vehicle trajectory specification in presence of
  traffic lights with known or uncertain switching times}.
\newblock \bibinfo{journal}{Transportation research record}
  \bibinfo{volume}{2674}, \bibinfo{pages}{53--66}.
\bibitem[{Typaldos et~al.(2020b)Typaldos, Papamichail and
  Papageorgiou}]{typaldos2020minimization}
\bibinfo{author}{Typaldos, P.}, \bibinfo{author}{Papamichail, I.},
  \bibinfo{author}{Papageorgiou, M.}, \bibinfo{year}{2020}b.
\newblock \bibinfo{title}{Minimization of fuel consumption for vehicle
  trajectories}.
\newblock \bibinfo{journal}{IEEE Transactions on Intelligent Transportation
  Systems} \bibinfo{volume}{21}, \bibinfo{pages}{1716--1727}.
\bibitem[{Typaldos et~al.(2021)Typaldos, Volakakis, Papageorgiou and
  Papamichail}]{typaldos2021vehicle}
\bibinfo{author}{Typaldos, P.}, \bibinfo{author}{Volakakis, V.},
  \bibinfo{author}{Papageorgiou, M.}, \bibinfo{author}{Papamichail, I.},
  \bibinfo{year}{2021}.
\newblock \bibinfo{title}{Vehicle-based trajectory specification in presence of
  traffic lights with stochastic switching times}.
\newblock \bibinfo{journal}{IFAC-PapersOnLine} \bibinfo{volume}{54},
  \bibinfo{pages}{298--305}.
\bibitem[{Xie et~al.(2017)Xie, Liu and Hauser}]{xie2017differential}
\bibinfo{author}{Xie, Z.}, \bibinfo{author}{Liu, C.K.},
  \bibinfo{author}{Hauser, K.}, \bibinfo{year}{2017}.
\newblock \bibinfo{title}{Differential dynamic programming with nonlinear
  constraints}, in: \bibinfo{booktitle}{2017 IEEE International Conference on
  Robotics and Automation (ICRA)}, \bibinfo{organization}{IEEE}. pp.
  \bibinfo{pages}{695--702}.
\bibitem[{Yakowitz(1986)}]{yakowitz1986stagewise}
\bibinfo{author}{Yakowitz, S.}, \bibinfo{year}{1986}.
\newblock \bibinfo{title}{The stagewise kuhn-tucker condition and differential
  dynamic programming}.
\newblock \bibinfo{journal}{IEEE transactions on automatic control}
  \bibinfo{volume}{31}, \bibinfo{pages}{25--30}.
\bibitem[{Yang et~al.(2020)Yang, Almutairi and Rakha}]{yang2020eco}
\bibinfo{author}{Yang, H.}, \bibinfo{author}{Almutairi, F.},
  \bibinfo{author}{Rakha, H.}, \bibinfo{year}{2020}.
\newblock \bibinfo{title}{Eco-driving at signalized intersections: A multiple
  signal optimization approach}.
\newblock \bibinfo{journal}{IEEE Transactions on Intelligent Transportation
  Systems} \bibinfo{volume}{22}, \bibinfo{pages}{2943--2955}.
\bibitem[{Yang et~al.(2016)Yang, Rakha and Ala}]{yang2016eco}
\bibinfo{author}{Yang, H.}, \bibinfo{author}{Rakha, H.}, \bibinfo{author}{Ala,
  M.V.}, \bibinfo{year}{2016}.
\newblock \bibinfo{title}{Eco-cooperative adaptive cruise control at signalized
  intersections considering queue effects}.
\newblock \bibinfo{journal}{IEEE Transactions on Intelligent Transportation
  Systems} \bibinfo{volume}{18}, \bibinfo{pages}{1575--1585}.
\bibitem[{Zhao et~al.(2018)Zhao, Ngoduy, Shepherd, Liu and
  Papageorgiou}]{zhao2018platoon}
\bibinfo{author}{Zhao, W.}, \bibinfo{author}{Ngoduy, D.},
  \bibinfo{author}{Shepherd, S.}, \bibinfo{author}{Liu, R.},
  \bibinfo{author}{Papageorgiou, M.}, \bibinfo{year}{2018}.
\newblock \bibinfo{title}{A platoon based cooperative eco-driving model for
  mixed automated and human-driven vehicles at a signalised intersection}.
\newblock \bibinfo{journal}{Transportation Research Part C: Emerging
  Technologies} \bibinfo{volume}{95}, \bibinfo{pages}{802--821}.

\end{thebibliography}

\end{document}